\documentclass[a4paper,DIV=12,10pt]{scrartcl}
 
\usepackage[utf8]{inputenc}
\usepackage{amsfonts,amssymb,amsmath}
\usepackage{epsfig}
\usepackage{MnSymbol}
\usepackage{booktabs}
\usepackage{float}
\usepackage{hyperref}
\usepackage[ruled,vlined,algosection,linesnumbered]{algorithm2e}
\usepackage{paralist}
\usepackage{color}
\usepackage{geometry}

\usepackage{caption}
\usepackage{subcaption}
\usepackage[binary-units=true]{siunitx-v2}

\newcommand{\ud}{\,\mathrm{d}}
\newcommand{\R}{\mathbb{R}}
\newcommand{\N}{\mathbb{N}}

\newcommand{\disc}{{\operatorname{disc}}}

\renewcommand{\vec}{\boldsymbol}

\newcommand*{\ldblbrace}{\{\mskip-5mu\{}
\newcommand*{\rdblbrace}{\}\mskip-5mu\}}

\numberwithin{equation}{section}
\numberwithin{figure}{section}
\numberwithin{table}{section}

\newtheorem{defi}{Definition}[section]

\newtheorem{lem}[defi]{Lemma}
\newtheorem{rem}[defi]{Remark}

\newtheorem{prob}[defi]{Problem}

\newenvironment{mproof}{\paragraph{Proof}}{\hfill$\blacksquare$}

\makeatletter
\let\@fnsymbol\@arabic
\makeatother

\begin{document}
	
\title{An energy-efficient GMRES--Multigrid solver for space-time finite element computation of dynamic poro- and thermoelasticity}

\author{
	Mathias Anselmann$^\ast$\thanks{anselmann@hsu-hh.de (corresponding author)}\;, Markus Bause$^\ast$,  Nils Margenberg$^\ast$, Pavel Shamko$^\ast$\\
	{\small ${}^\ast$ Helmut Schmidt University, Faculty of Mechanical and Civil Engineering, Holstenhofweg 85,}\\ 
	{\small 22043 Hamburg, Germany}
}

\date{}
\maketitle
			
\begin{abstract}
\textbf{Abstract}

\medskip \noindent 
We present families of space-time finite element methods (STFEMs) for a coupled hyperbolic-parabolic system of poro- or thermoelasticity. Well-posedness of the discrete problems is proved. Higher order approximations inheriting most of the rich structure of solutions to the continuous problem on computationally feasible grids are naturally embedded. However, the block structure and solution of the algebraic systems become increasingly complex for these members of the families. We present and analyze a robust geometric multigrid (GMG) preconditioner for GMRES iterations. The GMG method uses a local Vanka-type smoother. Its action is defined in an exact mathematical way. Due to nonlocal coupling mechanisms of unknowns, the smoother is applied on patches of elements. This ensures the damping of error frequencies. In a sequence of numerical experiments, including a challenging three-dimensional benchmark of practical interest, the efficiency of the solver for STFEMs is illustrated and confirmed. Its parallel scalability is analyzed. Beyond this study of classical performance engineering, the solver's energy efficiency is investigated as an additional and emerging dimension in the design and tuning of algorithms and their implementation on the hardware.
\end{abstract}


\clearpage
\section{Introduction}
\label{Sec:Introduction}

\subsection{Mathematical model}

We present and investigate numerically a geometric multigrid preconditioning technique, based on a local Vanka-type smoother, for solving by GMRES iterations the linear systems that arise from space-time finite element discretizations of the coupled hyperbolic-parabolic system of dynamic poroelasticity
\begin{subequations}
\label{Eq:HPS}
\begin{alignat}{3}
\label{Eq:HPS_1}
\rho \partial_t^2 \vec u - \nabla \cdot (\vec C \vec \varepsilon (\vec u)) + \alpha \vec \nabla p & = \rho \vec f\,, && \quad \text{in } \;
\Omega  \times (0,T]\,,\\[1ex]
\label{Eq:HPS_2}
c_0\partial_t p + \alpha \nabla \cdot \partial_t \vec u  -  \nabla \cdot (\vec K \nabla p)  & = g\,, && \quad \text{in } \; \Omega \times (0,T]\,,\\[1ex]
\label{Eq:HPS_3}
\vec u (0) = \vec u_0\,, \quad \partial_t \vec u (0) = \vec u_1\,, \quad p(0) & = p_0\,, && \quad \text{in } \; \Omega\times \{0\} \,,\\[1ex]
\label{Eq:HPS_4}
\vec u & = \vec u_D\,, && \quad \text{on } \; \Gamma_{\vec u}^{D} \times (0,T]\,,\\[1ex]
\label{Eq:HPS_5}
-(\vec C\vec \varepsilon(\vec u) - \alpha p\vec E) \vec n  & = \vec t_N\,,  && \quad \text{on } \; \Gamma_{\vec u}^{{N}} \times (0,T]\,,\\[1ex]
\label{Eq:HPS_6}
p & = p_D\,, && \quad \text{on } \; \Gamma_p^{{D}} \times (0,T]\,,\\[1ex]
\label{Eq:HPS_7}
- \vec K \nabla p \cdot \vec n  & = p_N\,, && \quad \text{on } \; \Gamma_p^{{N}} \times (0,T]\,.
\end{alignat}
\end{subequations}
In \eqref{Eq:HPS}, $\Omega \subset \R^d$, with $d\in \{2,3\}$, is an open bounded Lipschitz domain with outer unit normal vector $\vec n$ to the boundary $\partial \Omega$ and $T>0$ is the final time point. We let $\partial \Omega = {\Gamma_{\vec u}^{D}} \cup {\Gamma_{\vec u}^{N}}$ and $\partial \Omega = {\Gamma_{p}^{D}}\cup {\Gamma_{p}^{N}}$ with closed portions $\Gamma_{\vec u}^{D}$ and $\Gamma_{p}^{D}$ of non-zero measure. Important applications of the model \eqref{Eq:HPS}, that is studied as a prototype system, arise in poroelasticity; cf.~\cite{S00} and \cite{B41,B55,B72}. In poroelasticity, Eqs.~\eqref{Eq:HPS} are referred to as the dynamic Biot model. The system \eqref{Eq:HPS} is used to describe flow of a slightly compressible viscous fluid through a deformable porous matrix. The small deformations of the matrix are described by the Navier equations of linear elasticity, and the diffusive fluid flow is described by Duhamel’s equation. The unknowns are the effective solid phase displacement $\vec u$ and the effective fluid pressure $p$. The quantity $\vec \varepsilon (\vec u):= (\nabla \vec u + (\nabla \vec u)^\top)/2$ denotes the symmetrized gradient or strain tensor and $\vec E\in \R^{d,d}$ is the identity matrix. Further, $\rho$ is the effective mass density, $\vec C$ is Gassmann’s fourth order effective elasticity tensor, $\alpha$ is Biot’s pressure-storage coupling tensor, $c_0$ is the specific storage coefficient and $\vec K$ is the permeability field. For brevity, the positive quantities $\rho>0$, $\alpha>0$ and $c_0 >0$ as well as the tensors $\vec C$ and $\vec K$ are assumed to be constant in space and time. The tensors $\vec C$ and $\vec K$ are assumed to be symmetric and positive definite,
\begin{subequations}
\label{Eq:PosDef}
\begin{alignat}{2}
\label{Eq:PosDefC}
\exists k_0>0 \; \forall  \vec \xi = \vec \xi^\top \in \R^{d,d}:  &\quad \sum_{i,j,k,l=1}^d \xi_{ij} C_{ijkl} \xi_{kl} \geq k_0 \sum_{j,k=1}^d |\xi_{jk}|^2\,,\\[1ex]
\label{Eq:PosDefK}
 \exists k_1>0 \; \forall \vec \xi \in \R^d:  & \quad \sum_{i,j,=1}^d \xi_{i} K_{ij} \xi_{j} \geq k_1 \sum_{i=1}^d |\xi_{i}|^2\,.
\end{alignat}
\end{subequations}
 In order to enhance physical realism in poroelasticity, generalizations of the model \eqref{Eq:HPS} have been developed and investigated in, e.g., \cite{BKNR22,MW12}. We note that the system \eqref{Eq:HPS} is also formally equivalent to the classical coupled thermoelasticity system which describes the flow of heat through an elastic structure; cf.~\cite{JR18} and \cite{C72,L86}). In that context, $p$ denotes the temperature, $c_0$ is the specific heat of the medium, and $\vec K$ is the conductivity. Then, the quantity $\alpha \nabla p$ arises from the thermal stress in the structure, and the term $\alpha \nabla \cdot \partial_t \vec u$ models the internal heating due to the dilation rate. Well-posedness of \eqref{Eq:HPS} is ensured; cf., e.g., \cite{JR18,S89,STW22}. This can be shown by different mathematical techniques, by semigroup methods \cite[Thm.\ 2.2]{JR18}, Rothe's method \cite[Thm.\ 4.18 and Cor.\ 4.33]{S89} and Picard's theorem \cite[Thm.\ 6.2.1]{STW22}. In these works, boundary conditions different to the ones in \eqref{Eq:HPS_4} to   \eqref{Eq:HPS_4} are partly used. 
\subsection{Space-time finite element and multigrid techniques}
The coupled hyperbolic-parabolic structure of the system \eqref{Eq:HPS} of partial differential equations adds an additional facet of complexity onto its numerical simulation. A natural and promising approach for the numerical approximation of coupled systems is given by the application of space-time finite element methods that are based on an uniform treatment of the space and time variables by variational formulations. Therein, the discrete unknown functions are defined on the entire space-time domain $\Omega\times [0,T]$ and are represented in terms of basis functions. This facilitates the discretization of even complex coupling terms that may involve combinations of temporal and spatial derivatives as in \eqref{Eq:HPS_2}. Moreover, space-time finite element methods enable the natural construction of higher order members in the families of schemes. Higher order discretizations offer the potential to achieve accurate results on computationally feasible grids with a minimum of numerical costs. In particular for second-order hyperbolic problems, higher order time discretization schemes have shown their superiority, measured in terms of accuracy over temporal mesh resolution and compute time; cf., e.g., \cite{ABBM20}. Finally, space-time adaptivity based on a-posteriori error control by duality concepts and multi-rate in time approaches become feasible. For this we refer to \cite{BR03} for the presentation of the abstract framework of duality based error control and, e.g., to \cite{BGR10,BKB22} for its application to physical problems. For the variational time discretization of second-order hyperbolic and parabolic equations, various families have been proposed and investigated recently; cf., e.g., \cite{ABBM20,BKR22,KM04,FTW19,SZ22}. They differ by the approximation order and regularity in time of the discrete solution as well as stability and energy conservation properties. In this work, we use a \textbf{discontinuous Galerkin time discretization of arbitrary polynomial order} and focus on the efficient solution of the resulting complex algebraic block systems. Higher order (variational) time discretizations for three-dimensional applications put a facet of complexity on the algebraic solver and have rarely been studied yet in the literature. For the discretization of the space variables, numerous finite element approaches have been proposed and studied, in particular, for the quasi-static Biot system of flow in deformable porous media. The latter model results from \eqref{Eq:HPS} by neglecting the acceleration term $\rho \partial_t^2 \vec u$ in \eqref{Eq:HPS_1}. For the quasi-static Biot model, three-field formulations with the fluid flux $\vec q = - \vec K \nabla p$ as an additional explicit variable and the application of mixed finite element methods for the spatial discretization of $p$ and $\vec q$ attracted strong interest of researchers; cf., e.g., \cite{BRK17,HK18,KR18,PW07,PW08} and the references therein. In this work, we use an approximation of the variables $\vec u$ and $p$ by higher order inf-sup stable pairs of finite element spaces with continuous or discontinuous discrete variable $p$. Precisely, either the well-known Taylor--Hood family of finite element pairs for $\vec u$ and $p$ or a continuous approximation of $\vec u$ combined with a discontinuous approximation of $p$ is used; cf.~\eqref{Eq:DefVhQh}. In the literature, the latter option has shown its superiority for the approximation of incompressible viscous flow problems \cite{JM01}. Structure-preserving approximations with the fluid flux and stress tensor as additional variables, that are based on a reformulation of the coupled problem \eqref{Eq:HPS} as a first-order in space and time system (cf.~\cite[Subsec.~7.1]{FTW19}), remain as a work for the future. 

\begin{figure}[t]
	\begin{center}
		\includegraphics[width=8cm]{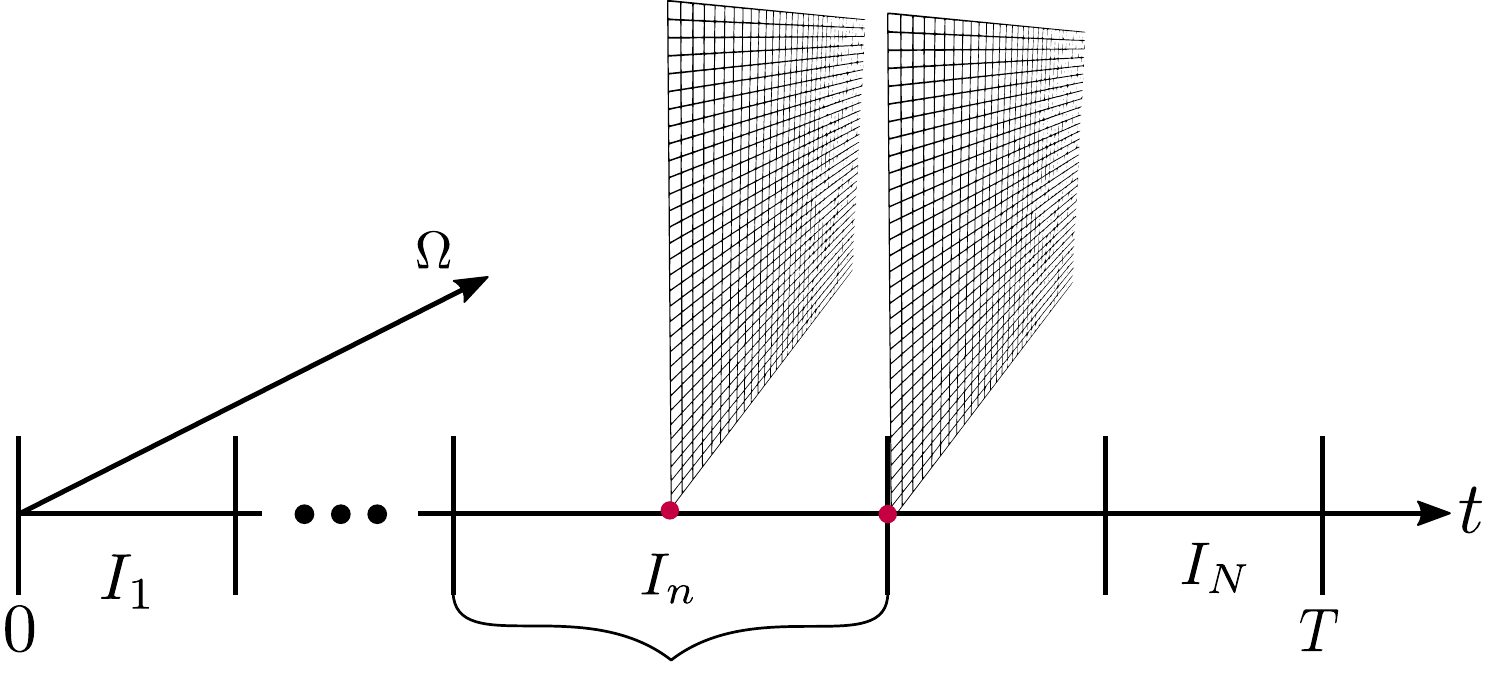}
	\end{center}
	\caption{Space-time mesh for a piecewise linear ($k=1$)  discontinuous Galerkin time discretization and a Lagrange basis w.r.t.\ the $k+1$ Gauss-Radau quadrature points of $I_n$.}	
	\label{Fig:dG1}
\end{figure}
In our space-time finite element  approach, we choose a discontinuous temporal test basis such a time-marching process is obtained. On each subinterval $I_n=(t_{n-1},t_n]$ of the time mesh $\mathcal{M}_\tau := \{I_1,\ldots, I_N\}$, we use a Lagrangian basis with respect to the $(k+1)$ Gauss--Radau quadrature points of $I_n$ for the time interpolation. Two discrete systems differing by the treatment of the term $\nabla \cdot \partial_t \vec u$ in \eqref{Eq:HPS_2} are then proposed. Dirichlet boundary conditions are implemented in a weak form by means of Nitsche's method and/or a discontinuous (spatial) approximation of the scalar variable $p$. Well-posedness of the discrete problems is carefully proved for arbitrary polynomial degrees of the discretization. The degrees of freedom on $I_n$ are associated with the $(k+1)$ Gauss--Radau quadrature nodes of this subinterval; cf.~Fig.\ \ref{Fig:dG1}. On each subinterval $I_n$, the resulting algebraic problem is a linear system of equations with a block matrix of $(k+1)\times (k+1)$ subsystems where each of the subsystems exhibits a complex structure itself. Here, $k$ denotes the piecewise polynomial order of the time discretization; cf.~\cite{AB22,AB21,HST14}. The solution of these algebraic problems with complex block structure demands for efficient and robust iteration schemes. Geometric multigrid (GMG) methods are known as the most efficient iterative techniques for the solution of large linear systems arising from the discretization of partial differential equations, particularly of elliptic type; cf.\ \cite{BL11,H85,TOS01} for a general overview.  In contrast to other iterative solvers (like for instance BiCGStab or GMRES), multigrid solvers converge with a rate which is independent of the (spatial) mesh size and require computational costs that are only linearly dependent on the number of unknowns. Massively parallel implementations of GMG methods on modern architectures show excellent scalability properties and their high efficiency has been recognized, e.g., in \cite{LLRS94,GHJRW16,GHJRW16}. Today they are widely used and are employed in many variants; cf.\ \cite{DJRWZ18,JT00,T99} for flow and saddle point problems. Nowadays, GMG methods are employed as preconditioner in Krylov subspace iterations, for instance GMRES iterations, to enhance their robustness, which is also done here.  Analyses of GMG methods (cf., e.g., \cite{DJRWZ18,HKXZ16,M06}) have been done for linear systems in saddle point form, with matrix $\vec{\mathcal A}= \begin{pmatrix} \vec A & \vec B^\top \\ \vec B & - \vec C \end{pmatrix}$ and symmetric and positive definite submatrices $\vec A$ and $\vec C$, arising for instance from mixed discretizations of the Stokes problem. If higher order discontinuous Galerkin time stepping is used, the resulting linear system matrix is a $(k+1)\times (k+1)$ block matrix with each of its blocks being of the form $\vec{\mathcal A}$. This imposes an additional facet of complexity on the geometric multigrid preconditioner. To the best of our knowledge,  analyses of GMG methods arising form higher order (variational) time discretizations of the coupled system \eqref{Eq:HPS} or related ones are still missing in the literature. 

GMG methods exploit different mesh levels of the underlying problem in order to reduce different frequencies of the error by employing a relatively cheap smoother on each grid level. Different iterative methods have been proposed in the literature as smoothing procedure; cf.~\cite{DJRWZ18} and the references therein.  They range from low-cost methods like Richardson, Jacobi, and SOR applied to the normal equation of the linear system to  collective smoothers, that are based on the solution of small local problems. Here, we use a Vanka-type smoother \cite{V86,M06,JT00} of the family of collective methods. Numerical computations show that an elementwise application of the Vanka smoother fails to reduce  the high frequencies of the error on the multigrid levels. The reason for this comes through spatial interelement couplings of degrees of freedom of the scalar variable $p$ in \eqref{Eq:HPS}. This holds for the discontinuous and continuous spatial approximation of $p$. For the spatially continuous approximation of $p$, the coupling of degrees of freedom is obvious. For the spatially discontinuous approximation of $p$, the interelement coupling of degrees of freedom of $p$ is due to the presence of the diffusion term in \eqref{Eq:HPS_2} and the surface integrals involving averages and jumps of the discrete solution over interelement faces resulting from its discretization. This is in contrast to the Navier--Stokes equations for which elementwise Vanka-type smoother are known to perform excellently if discontinuous pressure approximations are applied; cf.\ \cite{AB22,AB21,JT00}. The difference comes through the fact that in case of the Navier--Stokes system no interelement coupling of pressure degrees of freedom arises. As a remedy, we propose the application of the Vanka-type smoother on cell patches that are linked to the grid nodes and  built from four neighbored cells in two dimensions and eight neighbored elements in three dimensions, with appropriate adaptations for grid nodes close to or on the domain's boundary. Further we employ some relaxation strategy in the smoothing steps. Then an efficient damping of frequencies in the error on the multigrid hierarchy is obtained. Here, this Vanka-type smoother is presented in an mathematically exact way and its performance properties are investigated by numerical experiments. Our numerical experiments confirm that GMRES iterations that are preconditioned by the proposed multigrid method converge at a rate which is (nearly) independent of the mesh sizes in space and the time step size; cf.\ also \cite{AB22,AB21}. In the smoothing steps, the application of the Vanka-type operator requires the solution of patchwise linear systems of equations. For \eqref{Eq:HPS} and three space dimensions the size of the local systems no longer continues to remain negligible, in particular if a higher polynomial order in space and/or time is applied. Their solution deserves attention. In our approach we compute of inverse of the patchwise system matrix by LAPACK routines to solve to local system of equations. Details of an data-structure efficient MPI-based parallel implementation of the GMG algorithm in our in-house high-performance frontend solver for the deal.II finite element library \cite{Aetal21} can be found in \cite{AB21} where this is presented for the Navier--Stokes system. Multi-field formulations of the system \eqref{Eq:HPS}, that are not studied in this work, further increase the complexity of the arising linear systems; cf.~\cite[Subsec.~7.1]{FTW19}. However, the methodology of the linear solver, that is developed here, offers the potential of its generalization to more advanced problem formulations and their discretizations.

\subsection{Energy efficiency}

In the past, performance engineering and hardware engineering for large scale simulations of physical phenomena have been eclipsed by the longing for ever more performance where faster seemed to be the only paradigm. ''Classical'' performance engineering has been applied to enhance, firstly, the efficiency of the current method on the target hardware or to find numerical alternatives that might better fit to the hardware in use and/or, secondly, to develop other numerical methods can be found to improve the numerical efficiency. Tuning both simultaneously is called hardware-oriented numerics in the literature; cf.\ \cite{TBK06,TGBBW10}. Since recently, a growing awareness of energy consumption in computational science, particularly, in extreme scale computing with a view to exascale computing has raised; cf., e.g., \cite{R17}. It has been observed that as a consequence of decades of performance-centric hardware development there is a huge gap between pure performance and energy efficiency. An analysis of our algorithm's parallel scaling and energy consumption properties by performance models exceeds the scope of this work and would overburden it. However, since energy consumption of application codes on the available hardware is of growing awareness and a key for future improvements, we study the energy consumption and parallel scaling properties of our algorithm and its implementation by three-dimensional numerical experiments. The development of a proper model that quantifies performance and energy efficiency in some appropriate metric and can be used for a code optimization still deserves research and is left as a work for the future.

\subsection{Outline of the work}

This work is organized as follows. In Sec.\ \ref{Sec:Not} we introduce our notation. In Sec.\ \ref{Subsec:Disc} the space-time finite element approximation of arbitrary order of \eqref{Eq:HPS} is derived and well-posedness of the fully discrete problem is proved. Our GMRES--GMG solver is introduced in Sec.\ \ref{Sec:AlgSol}. In Sec.\ \ref{Sec:NumExp} our performed numerical computations for analyzing the performance properties of the overall approach are presented. In Sec.\ \ref{Sec:SumOut} we end with a summary and conclusions. In the appendix, supplementary results are summarized.     

\section{{Basic notation}}
\label{Sec:Not}

In this work, standard notation is used. We denote by $H^1(\Omega)$ the Sobolev space of $L^2(\Omega)$ functions with first-order weak derivatives in $L^2(\Omega)$. Further, $H^{-1}(\Omega)$ is the dual space of $H^1_{0}(\Omega)$, with the standard modification if the Dirichlet condition is prescribed on a part $\Gamma^D\subset \Gamma$ of the boundary $\partial \Omega$ only; cf.~\eqref{Eq:HPS}. The latter is not explicitly borne out by the notation of $H^{-1}(\Omega)$. It is always clear from the context. For vector-valued functions we write those spaces bold. By $\langle \cdot, \cdot \rangle_S$ we denote the $L^2(S)$ inner product for a domain $S$. For $S=\Omega$, we simply write $\langle \cdot, \cdot \rangle$. For the norms of the Sobolev spaces the notation is 
\begin{align*}
	\| \cdot \| := \| \cdot\|_{L^2}\,,\qquad 
	\| \cdot \|_1 := \| \cdot\|_{H^1}\,.
\end{align*} 
For short, we put 
\[
Q:=L^2(\Omega)\qquad \text{and} \qquad \vec V:= \left(H^1(\Omega)\right)^d\,.
\]
For a Banach space $B$, we let $L^2(0,T;B)$ be the Bochner space of $B$-valued functions, equipped with its natural norm. For a subinterval $J\subseteq [0,T]$, we will use the notation $L^2(J;B)$ for the corresponding Bochner space. In what follows, the constant $c$ is generic and indepedent of the size of the space and time meshes.

For the time discretization, we decompose the time interval $I:=(0,T]$ into $N$ subintervals $I_n=(t_{n-1},t_n]$, $n=1,\ldots,N$, where $0=t_0<t_1< \cdots < t_{N-1} < t_N = T$ such that $I=\bigcup_{n=1}^N I_n$. We put $\tau := \max_{n=1,\ldots, N} \tau_n$ with $\tau_n = t_n-t_{n-1}$. Further, the set $\mathcal{M}_\tau := \{I_1,\ldots,
I_N\}$ of time intervals is called the time mesh. For a Banach space $B$ and any $k\in
\N_0$, we let 
\begin{equation}
	\label{Def:Pk}
	\mathbb P_k(I_n;B) := \bigg\{w_\tau \,: \,  I_n \to B \,, \; w_\tau(t) = \sum_{j=0}^k 
	W^j t^j \; \forall t\in I_n\,, \; W^j \in B\; \forall j \bigg\}\,.
\end{equation}
For $k\in \N_0$ we define the space of piecewise polynomial functions in time with values in $B$ by 
\begin{equation}
	\label{Eq:DefYk}
	Y_\tau^{k} (B) := \left\{w_\tau: \overline I \rightarrow B  \mid w_\tau{}_{|I_n} \in
	\mathbb P_{k}(I_n;B)\; \forall I_n\in \mathcal{M}_\tau,\, w_\tau(0)\in B \right\}\subset L^2(I;B)\,.
\end{equation}
For any function $w: \overline I\to B$ that is piecewise sufficiently smooth with respect to the time mesh $\mathcal{M}_{\tau}$, for instance for $w\in Y^k_\tau (B)$, we define the right-hand sided and left-hand sided limit at a mesh point $t_n$ by
\begin{equation}
	\label{Eq:DefLim}
	w^+(t_n) := \lim_{t\to t_n+0} w(t) ,\quad\text{for}\; n<N\,,
	\qquad\text{and}\qquad
	w^-(t_n) := \lim_{t\to t_n-0} w(t) ,\quad\text{for}\; n>0\,.
\end{equation}
For the integration in time of a discontinuous Galerkin approach it is natural to use the right-sided $(k+1)$-point Gau{ss}--Radau quadrature formula. On the subinterval $I_n$, it reads as 
\begin{equation}
	\label{Eq:GF}
	Q_n(w) := \frac{\tau_n}{2}\sum_{\mu=1}^{k+1} \hat 
	\omega_\mu^{\operatorname{GR}} w(t_{n,\mu}^{\operatorname{GR}} ) \approx \int_{I_n} w(t) \ud t \,,
\end{equation}
where $t_{n,\mu}^{\operatorname{GR}}=T_n(\hat t_{\mu}^{\operatorname{GR}})$, for $\mu = 1,\ldots,k+1$, are the Gauss--Radau quadrature  points on $I_n$ and $\hat \omega_\mu^{\operatorname{GR}}$ the corresponding weights. Here, $T_n(\hat t):=(t_{n-1}+t_n)/2 + (\tau_n/2)\hat t$ is the affine transformation from $\hat I = [-1,1]$ to $I_n$ and $\hat t_{\mu}^{\operatorname{GR}}$ are the Gau{ss}--Radau quadrature points on $\hat I$. Formula \eqref{Eq:GF} is exact for all polynomials $w\in \mathbb P_{2k} (I_n;\R)$.  In particular, there holds that $t_{n,k+1}^{\operatorname{GR}}=t_n$.  

For the space discretization, let $\{\mathcal{T}_l\}_{l=0}^{L}$ be the decomposition on every multigrid level of $\Omega$ into (open) quadrilaterals or hexahedrals, with $\mathcal T_l = \{K_i\mid i=1,\ldots , N^{\text{el}}_l\}$, for $l=0,\ldots,L$. These element types are chosen for our implementation (cf.\ Sec.~\ref{Sec:NumExp}) that is based on the deal.II library \cite{Aetal21}.  The finest partition is $\mathcal T_h=\mathcal T_L$. We assume that all the partitions $\{\mathcal{T}_l\}_{l=0}^{L}$  are quasi-uniform with characteristic mesh size $h_l$ and $h_l=\gamma h_{l-1}$, $\gamma \in (0,1)$ and $h_0 = \mathcal O(1)$.  On the actual mesh level, the finite element spaces used for approximating the unknowns $\vec u$ and $p$ of \eqref{Eq:HPS} are of the form ($l\in \{0,\ldots,L\}$)
\begin{subequations}
	\label{Def:VhQh}
	\begin{alignat}{2}
		\label{Def:Vh}
		\vec V_{h_l}^l & := \{\vec v_h \in \vec V \cap C(\overline \Omega )^d\; : \; \vec v_{h}{}_{|K}\in 
		{\vec V(K)} \;\; \text{for all}\; K \in \mathcal{T}_l\}\,, \\[1ex]
		\label{Def:Qh}
		Q_{h_l}^{l,\text{cont}} & := \{\vec q_h \in Q\cap C(\overline \Omega ) \; : \; \vec q_{h}{}_{|K}\in {Q(K)} \;\; \text{for all}\; K \in \mathcal{T}_l\}\,, \\[1ex]
		\label{Def:Qh}
		Q_{h_l}^{l,\text{disc}} & := \{\vec q_h \in Q\; : \; \vec q_{h}{}_{|K}\in {Q(K)} \;\; \text{for all}\; K \in \mathcal{T}_l\}\,.
	\end{alignat}
\end{subequations}
By an abuse of notation, we skip the index $l$ of the mesh level when it is clear from the context and put
\begin{equation}
\label{Eq:DefVhQh}
\vec V_h := \vec V_h^l \quad \text{and} \quad Q_h := Q_h^l \;\; \text{with} \;\; Q_h^l\in \{Q_h^{l,\text{cont}}, Q_h^{l,\text{disc}}\}\,.
\end{equation}
 For the local spaces $\vec V(K)$ and $Q(K)$ we employ mapped versions of the pairs $\mathbb Q_r^d/\mathbb Q_{r-1}$ and $\mathbb Q_r^d/\mathbb P_{r-1}^\disc$, for $r\geq 2$. The pair $\mathbb Q_r^d/\mathbb Q_{r-1}$ with a (globally) continuous approximation of the scalar variable $p$ in $Q_h^{l,\text{cont}}$ is the well-known Taylor--Hood family of finite element spaces. The pair $\mathbb Q_r^d/\mathbb P_{r-1}^\disc$ comprises a discontinuous approximation of $p$ in the broken polynomial space $Q_h^{l,\text{disc}}$. For the Navier--Stokes equations, the multigrid method has shown to work best for higher-order finite element spaces with discontinuous discrete pressure; cf.\ \cite{JM01} and \cite{AB21}. For a further discussion of mapped and unmapped versions of the pair  $\mathbb Q_r^d/\mathbb P_{r-1}^\disc$ we refer to \cite[Subsec.~3.6.4]{J16}. For an analysis of stability properties of (spatial) discretizations for the quasi-static Biot system we refer to, e.g., \cite{RHOAGZ17}. In particular, both choices of the local finite element spaces, $\mathbb Q_r^d/\mathbb Q_{r-1}$ and $\mathbb Q_r^d/\mathbb P_{r-1}^\disc$, satisfy the inf--sup stability condition,
	\begin{equation}
	\label{Eq:InfSupCod}
	\inf_{q_h \in Q_h\backslash \{0\}} \sup_{\vec v_h\in \vec V_h\backslash \{\vec 0\}} \dfrac{b(\vec v_h,q_h)}{\|\vec v_h \|_1 \, \| q_h\|} \geq \beta > 0\,,
\end{equation}
for some constant $\beta$ independent of $h$; cf.~\cite{J16,MT02}. In \cite{A02,M01} optimal interpolation error estimates for mapped finite elements on quadrilaterals and hexahedra are studied. It turned out that the optimality is given for special families of triangulations. In two and three dimensions, families of meshes, which are obtained by a regular uniform refinement of an initial coarse grid, are among these special families. Such a regular refinement that is natural for the construction of the multigrid hierarchy is used in our computations.  Thus, for $\vec v\in \vec H^{r+1}(\Omega)$ and $q\in  H^r(\Omega)$ there exist approximations $i_h \vec v \in \vec V_h$ and $ j_h q \in Q_h$ such that, with some generic constant $c>0$ independent of $h$, 
	\begin{equation}
		\label{Eq:ApFES}
		\| \vec v - i_h \vec v \| + h \| \nabla (\vec v-i_h \vec v)\| \leq c h^{r+1} \qquad \text{and} \qquad \| q - j_h q\| \leq c  h^r\,.
	\end{equation}


\section{Space-time finite element approximation}
\label{Subsec:Disc}

For the discretization we rewrite \eqref{Eq:HPS} as a first-order in time system by introducing the new variable $\vec v:= \partial_t \vec u$. Then, we recover \eqref{Eq:HPS_1} and \eqref{Eq:HPS_2} as
\begin{subequations}
\label{Eq:HPS_11}
\begin{alignat}{2}
\label{Eq:HPS_8}
\partial_t \vec u - \vec v & = \vec 0\\[1ex]	
\label{Eq:HPS_9}
\rho \partial_t \vec v - \nabla \cdot (\vec C \vec \varepsilon (\vec u)) + \alpha \vec \nabla p & = \rho \vec f\,, \\[1ex]
\label{Eq:HPS_10}
c_0\partial_t p + \alpha \nabla \cdot \vec v  - \nabla \cdot (\vec K \vec \nabla p)  & = g
\end{alignat}
\end{subequations}
along with the initial and boundary conditions \ref{Eq:HPS_3} to \ref{Eq:HPS_7}. We employ discontinuous Galerkin methods (cf.\ \cite{T06}) for the discretization of the time variable and inf-sup stable pairs of finite elements (cf.\ Sec.~\ref{Sec:Not}) for the approximation of the space variables in \ref{Eq:HPS_11}. The derivation of the discrete scheme, presented below in Problem~\ref{Prob:DS}, is standard and not explicitly presented here. It follows the lines of \cite{BKR22}, where continuous in time Galerkin methods are applied to \eqref{Eq:HPS}, and \cite{AB22,AB21,HST11,HST13}, where discontinuous in time Galerkin methods are used to discretize the Navier--Stokes system. In contrast to \cite{BKR22}, Dirichlet boundary conditions are implemented here by Nitsche’s method \cite{B02,F78,N71}. This yields a strong link between two different families of inf-stable finite element pairs for the space discretization. The main reason for using Nitsche’s method here comes through our more general software framework. It captures problems on evolving domains solved on fixed computational background grids (cf.~\cite{AB21_2}), where Nitsche’s method is applied. We note that Nitsche’s method does not perturb the convergence behavior of the space-time discretization; cf.\ Sec.~\ref{Sec:NumExp}.

For the discrete scheme we need further notation. On the multigrid level $l$ with decomposition $\mathcal T_l$, for $\vec w_h, \vec \chi_h\in \vec V_h$ and $q_h, \psi_h \in Q_h$ we define
\begin{subequations}
	\label{Def:ACB}
	\begin{alignat}{2}
	\label{Def:A}
		A_\gamma (\vec w_{h},\vec \chi_h)	& := \langle \vec C \vec \varepsilon (\vec w_h),\vec \varepsilon (\vec \chi_h)\rangle - \langle \vec C\vec \varepsilon(\vec w_h) \vec n, \vec \chi_h\rangle_{\Gamma^D_{\vec u}}+ a_\gamma (\vec w_{h},\vec \chi_h)\,,\\[2ex]
	\label{Def:C}
		C (\vec \chi_{h},q_h)   & := -\alpha \langle \nabla \cdot \vec \chi_h, q_h\rangle +  \alpha \langle \vec \chi_h \cdot \vec n , q_h \rangle_{\Gamma^D_{\vec u}} \,, \\[2ex]
		B_\gamma (q_{h},\psi_h) & := \left\{\begin{array}{@{}ll}\langle \vec K \nabla q_h, \nabla \psi_h \rangle - \langle \vec K \nabla q_h \cdot \vec n, \psi_h \rangle_{\Gamma_p^D} + b_\gamma (q_h,\psi_h)\,,& \text{for}\;\; Q_h = Q_h^{l,\text{cont}}\,,\\[3ex]
	\label{Def:B}
		\begin{array}{@{}l} \displaystyle \sum_{K\in \mathcal T_l}\langle \vec K \nabla q_h, \nabla \psi_h \rangle_{K} - \sum_{F\in \mathcal F_h} \big (\langle \ldblbrace \vec K \nabla q_h \rdblbrace \cdot \vec n, \lsem \psi_h \rsem \rangle_{F}  \\[3ex] \displaystyle \quad + \langle \lsem q_h\rsem, \ldblbrace \vec K \nabla \psi_h \rdblbrace \cdot \vec n \rangle_{F}\big) + \sum_{F\in \mathcal F_h} \frac{\gamma}{h_F} \langle \lsem  q_h \rsem, \lsem \psi_h\rsem \rangle_F\,, \end{array} & \text{for}\;\; Q_h = Q_h^{l,\text{disc}}\,,\\[1ex]
		\end{array}\right.
	\end{alignat}
\end{subequations}
where, for $\vec w \in \vec H^{1/2}(\Gamma^D_{\vec u})$ and $q\in H^{1/2}(\Gamma^D_{p})$, 
\begin{subequations}
	\label{Def:ab}	
	\begin{alignat}{2}
		\label{Def:agam}
		a_\gamma (\vec w ,\vec \chi_h) & := - \langle \vec w, \vec C\vec \varepsilon(\vec \chi_h) \vec n \rangle_{\Gamma^D_{\vec u}}  + \frac{\gamma_a}{h_F}  \langle \vec w, \vec \chi_h \rangle_{\Gamma^D_{\vec u}}\,, \\[1ex]
		\label{Def:bgam}
		b_\gamma (q, \psi_h) & := - \langle q, \vec K \nabla \psi_h \cdot \vec n\rangle_{\Gamma_p^D} +  \frac{\gamma_b}{h_F}  \langle q, \psi_h \rangle_{\Gamma^D_{p}} \,.
	\end{alignat}
\end{subequations}
The second of the options in \eqref{Def:B} amounts to a symmetric interior penalty discontinuous Galerkin discretization of the scalar variable $p$; cf., e.g.,  \cite[Sec.~4.2]{PE12}. As usual, the average $\ldblbrace \cdot \rdblbrace$  and jump $\lsem \cdot \rsem$ of a function $w\in L^2(\Omega)$ on an interior face $F$ between two elements $K^+$ and $K^-$, such that $F=\partial K^+ \cap \partial K^-$, are 
\begin{equation*}
	\ldblbrace w \rdblbrace :=  \frac{1}{2} (w^{+}+ w^{-})\quad \text{and} \quad  \lsem w \rsem : =  w^{+} -  w^{-} \,.
\end{equation*}
For boundary faces $F \subset \partial K \cap \partial \Omega$, we set $\ldblbrace w\rdblbrace:= w_{|K}$ and $\lsem w \rsem := w_{|K}$. The set of all faces (interior and boundary faces) on the multigrid level $\mathcal T_l$ is denoted by $\mathcal F_h$. In the second of the options in \eqref{Def:B}, the parameter $\gamma$ of the last term has to be chosen sufficiently large, such that discrete coercivity on $Q_h$ of $B_\gamma$ is preserved. The local length $h_F$ is chosen {as $h_F = \ldblbrace h_F \rdblbrace := \frac{1}{2} (|K^+|_d + |K^-|_d)$} with Hausdorff measure $|\cdot |_d$; cf.~\cite[p.~125]{PE12}. {For boundary faces we set $h_F := |K|_d$.}
In \eqref{Def:ab}, the quantities $\gamma_a$ and $\gamma_b$ are the algorithmic parameters of the stabilization terms in the Nitsche formulation.
To ensure well-posedness of the discrete systems the parameters $\gamma_a$ and $\gamma_b$  have to be chosen  sufficiently large; cf.\ Appendix \ref{Sec:AppendixA}. {Based on our numerical experiments we choose the algorithmic parameter $\gamma$, $\gamma_a$ and $\gamma_b$ in \eqref{Def:B} and \eqref{Def:ab} as 
\begin{align*}
	\gamma_a = 5\cdot 10^4  \cdot r \cdot (r+1) \quad \text{and} \quad 
	\gamma = \gamma_b =  \frac{1}{2} \cdot r \cdot (r-1)\,,
\end{align*}
where $r$ is the polynomial degree of the finite element space \eqref{Def:Vh} for the displacement variable.}

Finally, for given $\vec f \in \vec H^{-1}(\Omega)$, $\vec u_D \in \vec H^{1/2}(\Gamma^D_{\vec u})$, $\vec t_N \in \vec H^{-1/2}(\Gamma^N_{\vec u})$ and $g\in H^{-1}(\Omega)$, $p_D \in H^{1/2}(\Gamma^D_{p})$, $p_N\in H^{-1/2}(\Gamma^N_{p})$ for $Q_h = Q_h^{l,\text{cont}}$, and properly fitted assumptions about the data for $Q_h = Q_h^{l,\text{disc}}$, we put 
\begin{subequations}
	\label{Def:LG}	
	\begin{alignat}{2}
	\label{Def:L}	
		F_\gamma (\vec \chi_h) & := \langle \vec f , \vec \chi_h\rangle - \langle \vec t_N,\vec \chi_h \rangle_{\Gamma^N_{\vec u}} + a_\gamma (\vec u_D ,\vec \chi_h)\,,\\[1ex]
	\label{Def:G}	
		G_\gamma (\psi_h) & := \left\{\begin{array}{@{}ll}
			\displaystyle \langle \langle g,\psi_h \rangle -  \alpha \langle \vec v_D \cdot \vec n , \psi_h \rangle_{\Gamma^D_{\vec u}}  - \langle p_N, \psi_h\rangle_{\Gamma^N_p} + b_\gamma(p_D, \psi_h)\,,& \text{for}\;\; Q_h = Q_h^{l,\text{cont}}\,,\\[3ex]
		\begin{array}{@{}l} \displaystyle \langle g,\psi_h \rangle -  \sum_{F\in \mathcal F_h^{D,\vec u}} \alpha \langle \vec v_D \cdot \vec n , \psi_h \rangle_{F} -  \sum_{F\in \mathcal F_h^{D,p}}\langle p_D, \ldblbrace \vec K \nabla \psi_h \rdblbrace \cdot \vec n \rangle_{F}\\[4ex] \displaystyle \;  + \sum_{F\in \mathcal F_h^{D,p}} \frac{\gamma}{h} \langle p_D, \lsem \psi_h\rsem \rangle_F - \sum_{F\subset \mathcal F_h^{N,p}} \langle p_N, \lsem \psi_h \rsem \rangle_F\,, \end{array} & \text{for}\;\; Q_h = Q_h^{l,\text{disc}}\,.
	   \end{array}\right.
	\end{alignat}
\end{subequations}
In the second of the options in \eqref{Def:G}, we denote by $\mathcal F_h^{D,p}\subset \mathcal F_h$ and $\mathcal F_h^{N,p}\subset \mathcal F_h$ the set of all element faces on the boundary parts $\Gamma_p^D$ and  $\Gamma_p^N$, respectively; cf.\ \eqref{Eq:HPS}. The second of the terms on the right-hand side of \eqref{Def:G}, with $\vec v_D = \partial_t \vec u_D$, is added to ensure consistency of the form \eqref{Def:C} in the fully discrete formulation \eqref{Eq:DPQ_A3} of \eqref{Eq:HPS_2}, i.e., that the discrete equation \eqref{Eq:DPQ_A3} is satisfied by the continuous solution to \eqref{Eq:HPS}.

We use a temporal test basis that is supported on the subintervals $I_n$; cf., e.g., \cite{AB21,HST13}. Then, a time marching process is obtained. In that, we assume that the trajectories $\vec u_{ \tau,h}$, $\vec v_{ \tau,h}$ and $p_{ \tau,h}$ have been computed before for all $t\in [0,t_{n-1}]$, starting with approximations  $\vec u_{\tau,h}(t_0) :=\vec u_{0,h}$, $\vec v_{\tau,h}(t_0) :=\vec u_{1,h}$ and $p_{\tau,h}(t_0) := p_{0,h}$ of the initial values $\vec u_0$, $\vec v_0$ and $p_{0,h}$. Then, we consider solving the following local problem on $I_n$.

\begin{prob}[Numerically integrated $I_n$-problem]
	\label{Prob:DSA}
	For given $\vec u_{h}^{n-1}:= \vec u_{\tau,h}(t_{n-1})\in \vec V_h$, $\vec v_{h}^{n-1}:=$ \linebreak[4]$  \vec v_{\tau,h}(t_{n-1})\in \vec V_h$,  and $p_{h}^{n-1}:= p_{\tau,h}(t_{n-1}) \in Q_h$ with  $\vec u_{\tau,h}(t_0) :=\vec u_{0,h}$, $\vec v_{\tau,h}(t_0) :=\vec u_{1,h}$ and $p_{\tau,h}(t_0) := p_{0,h}$, find $(\vec u_{\tau,h},\vec v_{\tau,h},p_{\tau,h}) \in \mathbb P_k (I_n;\vec V_h) \times \mathbb P_k (I_n;\vec V_h) \times \mathbb P_k (I_n;Q_h)$ such that
	\begin{subequations}
		\label{Eq:DPQ_A0}
		\begin{alignat}{2}
			\label{Eq:DPQ_A1}
			&\begin{aligned}
			& Q_n \big(\langle \partial_t \vec u_{\tau,h} , \vec \phi_{\tau,h} \rangle  - \langle \vec v_{\tau,h} , \vec \phi_{\tau,h} \rangle \big) + \langle \vec u^+_{\tau,h}(t_{n-1}), \vec \phi_{\tau,h}^+(t_{n-1})\rangle =  \langle \vec u_{h}^{n-1}, \vec \phi_{\tau,h}^+(t_{n-1})\rangle\,, \\[1ex]
			\end{aligned}\\
			\label{Eq:DPQ_A2}
			&\begin{aligned}
			& Q_n \Big(\langle \rho \partial_t \vec v_{\tau,h} , \vec \chi_{\tau,h} \rangle + A_\gamma(\vec u_{\tau,h}, \vec \chi_{\tau,h} ) + C(\vec \chi_{\tau,h},p_{\tau,h})\Big) + \langle \rho \vec v^+_{\tau,h}(t_{n-1}),  \vec \chi_{\tau,h}^+(t_{n-1})\rangle\\[1ex]  
			& \quad = Q_n \Big(F_\gamma(\vec \chi_{\tau,h})\Big) + \langle \rho \vec v_{h}^{n-1},  \chi_{\tau,h}^+(t_{n-1})\rangle \,,\\[1ex]
			\end{aligned}\\
			\label{Eq:DPQ_A3}
			&\begin{aligned}
			& Q_n \Big(\langle c_0 \partial_t p_{\tau,h},\psi_{\tau,h} \rangle  - C(\vec v_{\tau,h},\psi_{\tau,h})+ B_\gamma (p_{\tau,h}, \psi_{\tau,h})\Big) + \langle c_0 p^+_{\tau,h}(t_{n-1}), \psi_{\tau,h}^+(t_{n-1})\rangle \\[1ex] 
			&  \quad = Q_n \Big( G_\gamma(\psi_{\tau,h})\Big) + \langle c_0 p_{h}^{n-1}, \psi_{\tau,h}^+(t_{n-1})\rangle
			\end{aligned}
		\end{alignat}
	\end{subequations}
	for all $(\vec \phi_{\tau,h},\vec \chi_{\tau,h},\psi_{\tau,h})\in  \mathbb P_k (I_n;\vec V_h) \times \mathbb P_k (I_n;\vec V_h) \times \mathbb P_k (I_n;Q_h)$.
\end{prob}

The trajectories defined by Problem \ref{Prob:DSA}, for $n = 1,\ldots,N$, satisfy that $\vec u_{\tau,h},\vec v_{\tau,h}\in Y_\tau^k(\vec V_h)$ and $p_{\tau,h}\in Y_\tau^k(Q_h)$. The quadrature formulas on the left hand-side of \eqref{Eq:DPQ_A0} can be rewritten by time integrals since the Gauss--Radau formula \eqref{Eq:GF} is exact for all polynomials $w\in P_{2k}(I_n;\R)$. Well-posedness of Problem \ref{Prob:DSA} is ensured.

\begin{lem}[Existence and uniqueness of solutions to Problem \ref{Prob:DSA}]
\label{Lem:ExUniDS}
Problem \ref{Prob:DSA} admits a unique solution.
\end{lem}

\begin{mproof}
We prove Lem.\ \ref{Lem:ExUniDS} for $Q_h=Q_h^{l,\text{cont}}$ only, thus assuming the first of the options in \eqref{Def:B} and \eqref{Def:G}. For $Q_h=Q_h^{l,\text{disc}}$, the proof can be done similarly by using, in addition, standard techniques of error analysis for discontinuous Galerkin methods; cf., e.g., \cite[Sec.~4]{PE12}. Since Problem  \ref{Eq:DPQ_A0}  is finite-dimensional, it suffices to prove uniqueness of the solution. Existence of the solution then follows directly from its uniqueness. Let now $(\vec u^{(1)}_{\tau,h},\vec v^{(1)}_{\tau,h},p^{(1)}_{\tau,h}) \in \mathbb P_k (I_n;\vec V_h) \times \mathbb P_k (I_n;\vec V_h) \times \mathbb P_k (I_n;Q_h)$ and $(\vec u^{(2)}_{\tau,h},\vec v^{(2)}_{\tau,h},p^{(2)}_{\tau,h}) \in \mathbb P_k (I_n;\vec V_h) \times \mathbb P_k (I_n;\vec V_h) \times \mathbb P_k (I_n;Q_h)$  denote two triples of solutions to \eqref{Eq:DPQ_A1}. Their differences $(\vec u_{\tau,h},\vec v_{\tau,h},p_{\tau,h}) =(\vec u^{(1)}_{\tau,h},\vec v^{(1)}_{\tau,h},p^{(1)}_{\tau,h}) -(\vec u^{(2)}_{\tau,h},\vec v^{(2)}_{\tau,h},p^{(2)}_{\tau,h}) $ then satisfy the equations
	\begin{subequations}
	\label{Eq:DPQ_A6}
	\begin{alignat}{2}
		\label{Eq:DPQ_A7}
			& Q_n \big(\langle \partial_t \vec u_{\tau,h} , \vec \phi_{\tau,h} \rangle  - \langle \vec v_{\tau,h} , \vec \phi_{\tau,h} \rangle \big) + \langle \vec u^+_{\tau,h}(t_{n-1}), \vec \phi_{\tau,h}^+(t_{n-1})\rangle = 0\,, \\[1ex]
		\label{Eq:DPQ_A8}
			& Q_n \Big(\langle \rho \partial_t \vec v_{\tau,h} , \vec \chi_{\tau,h} \rangle + A_\gamma(\vec u_{\tau,h}, \vec \chi_{\tau,h} ) + C(\vec \chi_{\tau,h},p_{\tau,h})\Big) + \langle \rho \vec v^+_{\tau,h}(t_{n-1}), \vec \chi_{\tau,h}^+(t_{n-1})\rangle = 0 \,,\\[1ex]
		\label{Eq:DPQ_A9}
			& Q_n \Big(\langle c_0 \partial_t p_{\tau,h},\psi_{\tau,h} \rangle  - C(\vec v_{\tau,h},\psi_{\tau,h})+ B_\gamma (p_{\tau,h}, \psi_{\tau,h})\Big) + \langle c_0 p^+_{\tau,h}(t_{n-1}), \psi_{\tau,h}^+(t_{n-1})\rangle = 0 
	\end{alignat}
\end{subequations}
for all $(\vec \phi_{\tau,h},\vec \chi_{\tau,h},\psi_{\tau,h})\in  \mathbb P_k (I_n;\vec V_h) \times \mathbb P_k (I_n;\vec V_h) \times \mathbb P_k (I_n;Q_h)$. We let $\vec A_\gamma:  \vec V_h \mapsto \vec V_h$ be the discrete operator that is defined, for $\vec w_h \in  \vec V_h$ and all $\vec \phi_h\in \vec V_h$, by
\begin{equation}
	\label{Eq:DefOpAh}
	\langle \vec A_\gamma \vec w_h , \vec \phi_h \rangle = A_\gamma (\vec w_h,\vec \phi_h)\,.
\end{equation}
In \eqref{Eq:DPQ_A6} we choose $\vec \phi_{\tau,h}=\vec A_\gamma \vec u_{\tau,h}$, $\vec \chi_{\tau,h}=\vec v_{\tau,h} $ and $\psi_{\tau,h}=p_{\tau,h}$. Adding the resulting equations yields that 
\begin{equation}
\label{Eq:ExUniDS_1}
\begin{aligned}
& Q_n \big(\langle \partial_t \vec u_{\tau,h} , \vec A_\gamma \vec u_{\tau,h} \rangle  + \langle \rho \partial_t \vec v_{\tau,h} , \vec v_{\tau,h}  \rangle  + \langle c_0 \partial_t p_{\tau,h}, p_{\tau,h} \rangle + B_\gamma(p_{\tau,h},p_{\tau,h}) \big)   \\[1ex] 
& \quad + \langle \vec u^+_{\tau,h}(t_{n-1}), \vec A_\gamma u_{\tau,h}^+(t_{n-1})\rangle  + \langle \rho \vec v^+_{\tau,h}(t_{n-1}),  \vec v_{\tau,h}^+(t_{n-1})\rangle + \langle c_0 p^+_{\tau,h}(t_{n-1}), p_{\tau,h}^+(t_{n-1})\rangle  = 0 \,.
\end{aligned}
\end{equation}
Recalling  the exactness of the Gauss--Radau formula \eqref{Eq:GF} for $w\in \mathbb P_{2k}(I_n;\R)$, Eq.~\eqref{Eq:ExUniDS_1} yields that 
\begin{equation*}
	\label{Eq:ExUniDS_2}
	\begin{aligned}
		& \frac{1}{2} \int_{t_{n-1}}^{t_n} \frac{d}{dt} \big( \langle \vec A_\gamma \vec u_{\tau,h} , \vec u_{\tau,h} \rangle  + \langle \rho \vec v_{\tau,h} , \vec v_{\tau,h}  \rangle  + \langle c_0 p_{\tau,h}, p_{\tau,h} \rangle \big ) \ud t + Q_n\big(B_\gamma(p_{\tau,h},p_{\tau,h}) \big)  \\[1ex] 
		& \quad + \langle \vec u^+_{\tau,h}(t_{n-1}), \vec A_\gamma u_{\tau,h}^+(t_{n-1})\rangle  + \langle \rho \vec v^+_{\tau,h}(t_{n-1}),  \vec v_{\tau,h}^+(t_{n-1})\rangle + \langle c_0 p^+_{\tau,h}(t_{n-1}), p_{\tau,h}^+(t_{n-1})\rangle  = 0 \,.
	\end{aligned}
\end{equation*}
Using \eqref{Eq:DefOpAh}, this shows that 
\begin{equation}
	\label{Eq:ExUniDS_3}
	\begin{aligned}
		& A_\gamma (\vec u_{\tau,h}(t_n) , \vec u_{\tau,h}(t_n))  + \langle \rho \vec v_{\tau,h} (t_n), \vec v_{\tau,h}(t_n)  \rangle  + \langle c_0 p_{\tau,h}(t_n), p_{\tau,h}(t_n)\rangle +  2 Q_n \big(B_\gamma(p_{\tau,h},p_{\tau,h}) \big) \\[1ex]
		& \quad + A_\gamma(\vec u^+_{\tau,h}(t_{n-1}), \vec  u_{\tau,h}^+(t_{n-1}))  + \langle \rho \vec v^+_{\tau,h}(t_{n-1}),  \vec v_{\tau,h}^+(t_{n-1})\rangle + \langle c_0 p^+_{\tau,h}(t_{n-1}), p_{\tau,h}^+(t_{n-1})\rangle  = 0 \,.
	\end{aligned}
\end{equation}
From \eqref{Eq:ExUniDS_3} along with the discrete coercivity properties \eqref{Eq:CoercA} of $A_\gamma$ and \eqref{Eq:CoercB} of $B_\gamma$ we directly deduce that 
\begin{equation}
\label{Eq:ExUniDS_10}
\vec u_{\tau,h}(t_n)= \vec u_{\tau,h}^+(t_{n-1})=\vec 0\,, \quad \vec v_{\tau,h}(t_n)= \vec v_{\tau,h}^+(t_{n-1})=\vec 0 \quad \text{and} \quad p_{\tau,h}(t_n)=p_{\tau,h}^+(t_{n-1})= 0\,,
\end{equation}
as well as
\begin{equation}
\label{Eq:ExUniDS_11}
p_{\tau,h}\big(t_{n,\mu}^{\text{GR}}\big) = 0 \,, \quad \text{for}\;\; \mu = 1,\ldots, k+1\,.
\end{equation}
We note that $t_{n,k+1}^{\text{GR}}=t_n$. Relation \eqref{Eq:ExUniDS_11} now implies that $p_{\tau,h}\equiv 0$ on $I_n$. 

 To prove that $\vec u_{\tau,h}\equiv \vec 0$ and $\vec v_{\tau,h}\equiv \vec 0$,  by \eqref{Eq:ExUniDS_10} it is sufficient to show that $\vec u_{\tau,h}(t_{n,\mu}^{\text{G}})=\vec 0$ and $\vec v_{\tau,h}(t_{n,\mu}^{\text{G}}) = \vec 0$, for $\mu = 1,\ldots, k$, where $t_{n,\mu}^{\text{G}}$, for $\mu = 1,\ldots, k$, are the nodes of the $k$-point Gauss quadrature formula on $I_n$ that is exact for all polynomials in $\mathbb P_{2k-1}(I_n;\R)$. This is done now. Recalling \eqref{Eq:ExUniDS_10}, we conclude from \eqref{Eq:DPQ_A7} by a suitable choice of test functions that 
 \begin{equation}
 \label{Eq:ExUniDS_12}
\partial_t \vec u_{\tau,h}(t_{n,\mu}^{\text{GR}})=\vec v_{\tau,h}(t_{n,\mu}^{\text{GR}}) \,, \quad \text{for} \;\; \mu = 1,\ldots, k\,.	
\end{equation} 	
Next, choosing  $\chi_{\tau,h}=\vec v_{\tau,h}$ in \eqref{Eq:DPQ_A8} and recalling \eqref{Eq:ExUniDS_12} imply that 
\begin{equation}
\label{Eq:ExUniDS_13}
  Q_n \Big(\langle \rho \partial_t \vec v_{\tau,h} , \vec v_{\tau,h} \rangle + A_\gamma(\vec u_{\tau,h}, \partial_t \vec u_{\tau,h} ) \Big)= 0
\end{equation} 	
By the exactness of the Gauss--Radau formula \eqref{Eq:GF} for all $w\in \mathbb P_{2k}(I_n;\R)$ we have from \eqref{Eq:ExUniDS_13} that 
\begin{equation}
\label{Eq:ExUniDS_14}
\int_{t_{n-1}}^{t_n} \langle \rho \partial_t \vec v_{\tau,h} , \vec v_{\tau,h} \rangle  \ud t + \frac{1}{2} \int_{t_{n-1}}^{t_n} \frac{d}{dt} A_\gamma (\vec u_{\tau,h} , \vec u_{\tau,h} )  \ud t = 0 \,. 
\end{equation} 	
The second of the terms in \eqref{Eq:ExUniDS_14} vanishes by \eqref{Eq:ExUniDS_10}. The stability result of \cite[Lem.\ 2.1]{KM99} then implies that 
\begin{equation}	
\label{Eq:ExUniDS_15}
\vec v_{\tau,h}(t_{n,\mu}^{\text{G}}) = \vec 0 \,, \quad \text{for} \;\; \mu = 1,\ldots, k\,.
\end{equation}
From \eqref{Eq:ExUniDS_15} along with  \eqref{Eq:ExUniDS_10} we then deduce that $\vec v_{\tau,h}\equiv \vec 0$ on $I_n$.  Choosing the test function $\vec \phi_{\tau,h}=\vec u_{\tau_h}$ in  \eqref{Eq:DPQ_A6}, using $\vec v_{\tau,h}\equiv \vec 0$ and applying the stability result \cite[Lem.\ 2.1]{KM99}, it follows that 
\begin{equation}	
 	\label{Eq:ExUniDS_16}
 	\vec u_{\tau,h}(t_{n,\mu}^{\text{G}}) = \vec 0 \,, \quad \text{for} \;\; \mu = 1,\ldots, k\,.
 \end{equation}
 From \eqref{Eq:ExUniDS_16} along with  \eqref{Eq:ExUniDS_10} we then have that $\vec u_{\tau,h}\equiv \vec 0$ on $I_n$. Thus uniqueness of solutions to \eqref{Eq:DPQ_A0} and, thereby, well-posedness of Problem \ref{Prob:DSA}  is thus ensured.

\end{mproof}

In Appendix~\ref{Sec:AppendixB} an alternative formulation for the system \eqref{Eq:DPQ_A6} is still presented. It is based on using the time derivative $\partial_t \vec u_{\tau,h}$ of the primal variable $\vec u_{\tau,h}$ instead of the auxiliary variable $\vec v_{\tau,h}$ in \eqref{Eq:DPQ_A3}. In this case, an additional surface integral has to be included; cf.\  Eq.~\eqref{Eq:DPQ_3} .

\section{Algebraic solver by geometric multigrid preconditioned GMRES iterations}
\label{Sec:AlgSol}

On the algebraic level, the variational problem \eqref{Eq:DPQ_0} leads to linear systems of equations with complex block structure, in particular if higher order (piecewise) polynomial degrees $k$ for the approximation of the temporal variable are involved. This demands for a robust and efficient linear solver, in particular in the three-dimensional case $\Omega \subset \R^3$. For solving \eqref{Eq:DPQ_0} we consider using flexible GMRES iterations \cite{S03} that are preconditioned by a $V$-cycle geometric multigrid method (GMG) based on a local Vanka smoother. In \cite{AB22}, the combined GMG--GMRES approach is further embedded in a Newton iteration for solving space-time finite element discretizations of the Navier--Stokes system. Thus, nonlinear extensions of the prototype model \eqref{Eq:HPS} are expected to become feasible by our approach as well. 

To derive the algebraic form of \eqref{Eq:DPQ_A0}, the discrete functions $\vec u_{\tau,h}$, $\vec v_{\tau,h}$ and $p_{\tau,h}$ are represented in a Lagrangian basis $\{\chi_{n,m}\}_{m=1}^{k+1}\subset \mathbb P_k(I_n;\R)$ with respect to the $(k+1)$ Gauss--Radau quadrature points of $I_n$, such that 
\begin{subequations}
\label{Eq:TimeExp}
\begin{alignat}{2}
	\vec u_{\tau,h}{}_{|I_n}(\vec x,t) & = \sum_{m=1}^{k+1} \vec u_{n,m}(\vec x) \chi_{n,m}(t)\,, \quad \vec v_{\tau,h}{}_{|I_n}(\vec x,t) = \sum_{m=1}^{k+1} \vec v_{n,m}(\vec x) \chi_{n,m}(t)\,, \\[1ex]  p_{\tau,h}{}_{|I_n}(\vec x,t) & = \sum_{m=1}^{k+1} p_{n,m}(\vec x) \chi_{n,m}(t)\,.
\end{alignat}
\end{subequations}
The resulting coefficient functions $(\vec u_{n,m},\vec v_{n,m},p_{n,m})\in \vec V_h\times \vec V_h\times Q_h$, for $m=1,\ldots,k+1$, are developed in terms of the finite element basis of $\vec V_h$ and $Q_h$, respectively. Letting $\vec V_h= \operatorname{span}\{\vec \psi_1,\ldots,\vec \psi_R \}$ and $Q_h=\operatorname{span}\{\xi_1,\ldots,\xi_S \}$, we get that 
\begin{equation}
	\label{Eq:SpaceExp}
	\vec u_{n,m}(\vec x) = \sum_{r=1}^R u^{(r)}_{n,m} \vec \psi_r(\vec x)\,, \;\;
	\vec v_{n,m}(\vec x) = \sum_{r=1}^R v^{(r)}_{n,m} \vec \psi_r(\vec x) \quad \text{and} \quad 
	p_{n,m}(\vec x) = \sum_{s=1}^S p^{(s)}_{n,m}\,  \xi_s(\vec x)\,.
\end{equation}
For the coefficients of the expansions in \eqref{Eq:SpaceExp} we define the subvectors 
\begin{subequations}
\label{Eq:SpaceExp2}
\begin{alignat}{2}
\vec U_{n,m} & = \big(u^{(1)}_{n,m},\ldots,u^{(R)}_{n,m} \big)^\top\,, \;\; \vec V_{n,m}=  \big(v^{(1)}_{n,m},\ldots,v^{(R)}_{n,m} \big)^\top \quad \text{and} \\[1ex] 
\vec P_{n,m} & =  \big(p^{(1)}_{n,m},\ldots,p^{(S)}_{n,m} \big)^\top\,,  \quad \text{for}\;\;  m=1,\ldots,k+1\,,
\end{alignat}
\end{subequations}  
of the degrees of freedom for all Gauss--Radau quadrature points and the global solution vector on $I_n$ by
\begin{equation}
\label{Eq:DefXn}
\vec X_n^\top = \big((\vec V_{n,1})^\top,(\vec U_{n,1})^\top,(\vec P_{n,1})^\top,\ldots, (\vec V_{n,k+1})^\top,(\vec U_{n,k+1})^\top,(\vec P_{n,k+1})^\top \big)\,.
\end{equation}
We note that $\vec X_n$ comprises the (spatial) degrees of freedom for all $(k+1)$ Gauss--Radau nodes, representing the Lagrange interpolation points in time , of the subinterval $I_n$. The approximations at these time points will be computed simultaneously. Substituting \eqref{Eq:TimeExp} and \eqref{Eq:SpaceExp} into \eqref{Eq:DPQ_0} and choosing in \eqref{Eq:DPQ_0} the test basis $\{\chi_{n,m}\vec \psi_r,\chi_{n,m}\vec \psi_r,\chi_{n,m}\xi_s\}$, for $m=1,\ldots,k+1$, $r=1,\ldots, R$ and $s=1,\ldots,S$,  built from the trial basis in \eqref{Eq:TimeExp} and \eqref{Eq:SpaceExp}, we obtain the following algebraic system.
\begin{prob}[Algebraic $I_n$-problem]
\label{Prob:AlgSys}
For the vector $\vec X_n$, defined in \eqref{Eq:DefXn} along with \eqref{Eq:SpaceExp2}, of the coefficients of the expansions  \eqref{Eq:SpaceExp} solve 
\begin{equation}
\label{Eq:ALS_In}
	\vec A_n \vec X_n = \vec F_n\,,
\end{equation}
where the matrix $\vec A_n$ exhibits the $(k+1)\times (k+1)$ block structure 
\begin{equation}
\label{Eq:DefAn}
\vec A_n = \big(\vec A_{a,b}\big)_{a,b=1}^{k+1}
\end{equation}
with block submatrices $\vec A_{a,b}$, for $a,b=0,\ldots,k$, defined by
\begin{equation}
\label{Eq:DefAab}
\vec A_{a,b} = \left(\begin{array}{@{}ccc@{}} 
  - \vec M^{0,\vec V_h}_{{a}, {b}}
   & \vec M^{1,\vec V_h}_{{a},{b}}  
&  \vec 0 \\[1ex]
    \vec M^{1,\vec V_h}_{{a},{b}}  
& \vec A_{{a},{b}} + \vec N^{A}_{{a},{b}} 
& \vec C^\top_{{a},{b}}  \\[1ex]
    \vec 0
&  - \vec C_{{a},{b}}
& \vec M^{1,Q_h}_{a,b} + \vec B_{a,b} + \vec N^B_{a,b}
\end{array}\right)\,.	
\end{equation}
\end{prob}
For the choice $Q_h^l=Q_h^{l,\text{cont}}$ in \eqref{Eq:DefVhQh}, the explicit representation of the submatrices in \eqref{Eq:DefAab} reads as
\begin{equation*}
\begin{aligned}
\big(\vec M^{1,\vec V_h}_{{a},{b}}\big)_{i,j}  & := Q_n\big(\langle \rho \partial_t \chi_{n,a} \vec \psi_{j}, \chi_{n,b} \vec \psi_{i} \rangle \big) + \langle \rho \chi_{n,a}(t_{n-1}^+) \vec \psi_{j}, \chi_{n,b}(t_{n-1}^+) \vec \psi_{i} \rangle \,, \\[1ex]
 \big(\vec M^{0,\vec V_h}_{{a},{b}}\big)_{i,j} & := Q_n \big(\langle \rho \chi_{n,a} \vec \psi_{j}, \chi_{n,b} \vec \psi_{i} \rangle\big) \,, \\[1ex]
\big(\vec A_{{a},{b}}\big)_{i,j} & : = Q_n \big(\langle \vec C \vec \varepsilon (\chi_{n,a}\vec \psi_{j}), \vec \varepsilon (\chi_{n,b} \vec \psi_{i}) \rangle + \langle \vec C \vec \varepsilon (\chi_{n,a}\vec \psi_{j}) \vec n , \chi_{n,b} \vec \psi_{i} \rangle_{\Gamma^D_{\vec u}} \big)\,, \\[1ex]
\big(\vec N^{A}_{{a},\vec \psi_{b}}\big)_{i,j} 	 & :=  Q_n \big(a_\gamma(\chi_{n,a} \vec \psi_{i},\chi_{n,b} \vec \psi_{i})\big)
\end{aligned}
\end{equation*}
as well as 
\begin{align*}
\big(\vec C_{{a},{b}}\big)_{r,j} &  := Q_n \big(-\alpha\langle \nabla \cdot (\chi_{n,b} \vec \psi_{j}), \chi_{n,a} \xi_r \rangle + \alpha \langle \chi_{n,b} \vec \psi_{j} \cdot \vec n , \chi_{n,a} \xi_r \rangle_{\Gamma^D_{\vec u}}\big)\,,\\[1ex]
\big(\vec M^{1,Q_h}_{p_{a},p_{b}}\big)_{r,s}  & := Q_n\big(\langle c_0 \partial_t \chi_{n,a} \xi_s, \chi_{n,b} \xi_r \rangle \big) + \langle c_0 \chi_{n,a}(t_{n-1}^+) \xi_{s}, \chi_{n,b}(t_{n-1}^+) \xi_{r} \rangle\big) 
\end{align*}
and
\begin{align*}
\big(\vec B_{{a},{b}}\big)_{r,s}  &  := Q_n \big(\langle \vec K \nabla (\chi_{n,a}\xi_{s}), \nabla (\chi_{n,b} \xi_{r}) \rangle - \langle \vec K  \nabla (\chi_{n,a}\vec \xi_{s}) \cdot \vec n , \chi_{n,b} \xi_{r} \rangle_{\Gamma^D_{p}} \big)\,, \\[1ex]
\big(\vec N^B_{{a},{b}}\big)_{r,s}   & := Q_n \big(b_\gamma(\chi_{n,a} \xi_{s},\chi_{n,b} \xi_{r})\big)\,,
\end{align*}
with $a_\gamma(\cdot,\cdot)$ and $b_\gamma(\cdot,\cdot)$ being defined in \eqref{Def:ab} and $i,j=1,\ldots, R$ and $r,s=1,\ldots , S$. The right-hand side vector $\vec F_n$ in  \eqref{Eq:ALS_In} is defined similarly, according to \eqref{Eq:DPQ_0} along with \eqref{Def:LG}. Its definition is skipped here for brevity. We note that \eqref{Eq:HPS_8} is still multiplied by $\rho>0$ for the definition of the first row in \eqref{Eq:DefAab}. For the family of finite element pairs, corresponding to $Q_h^l=Q_h^{l,\text{disc}}$ in \eqref{Eq:DefVhQh}, the definition of the submatrix $B_{a,b}$ has to be adjusted to the second of the options in \eqref{Def:B}. Further, the contribution $\vec N^B_{a,b}$ has to be omitted. 

Increasing values of the piecewise polynomial degree in time $k$ enhance the complexity of the block structure of the system matrix $\vec A_n$ in \eqref{Eq:DefAn} along with \eqref{Eq:DefAab}. They impose an additional facet of challenge on the construction of efficient block preconditioners for \eqref{Eq:ALS_In}. We solve the linear system \eqref{Eq:ALS_In} for the unknown $\vec X_n$ on the subinterval $I_n$ by flexible GMRES iterations \cite{S03} that are preconditioned by a $V$-cycle geometric multigrid (GMG) algorithm \cite{TOS01}. Key ingredients of the GMRES--GMG approach are 
\begin{compactitem}
\item left preconditioning of each GMRES iteration by a single $V$-cycle GMG iteration,
\item prolongation to grids by interpolation and application of the transpose operator for restriction to grids,
\item $J_{\max}$ pre- and post-smoothing steps on each grid level $\{\mathcal T_l\}$ for $l=1,\ldots,L$ and
\item  application of a (parallel) direct solver on the coarsest mesh partition $\mathcal T_0$.
\end{compactitem}
In our computational studies presented in Sec.~\ref{Sec:NumExp}, we use $J_{\max}=4$ and the parallel direct solver \textit{SuperLU\_DIST} \cite{LD03}. For the implementation of the restriction and prolongation operators the deal.II classes \textit{MultiGrid} and MGTransferPrebuilt are used. For the deal.II finite element library we refer to \cite{Aetal21}. For details of the parallel implementation by the message passing interface (MPI) protocol of our GMG approach  we refer to \cite{AB21}.

The choice of the smoothing operator in the GMG method is extremely diverse; cf., e.g., \cite{DJRWZ18,JT00} for a further discussion. We use a collective smoother of Vanka type that is based on the solution of small local problems. Compared to the Navier--Stokes system \cite{AB21}, the dynamic Biot problem \eqref{Eq:HPS} requires adaptations in the construction of the local Vanka smoother. In inf-sup stable discretizations of the Navier--Stokes equations by the $\mathbb Q_r^d/\mathbb P_{r-1}^{\text{disc}}$, $r\geq 2$, family of finite element spaces no coupling between the pressure degrees of freedom is involved due to the discontinuous approximation of the pressure variable. This feature leads to excellent performance properties of the Vanka smoother \cite{AB21}. In contrast to this, the discretization \eqref{Eq:DPQ_A3} of \eqref{Eq:HPS_2} by $\mathbb P_{r-1}^{\text{disc}}$, $r\geq 2$, elements involves a coupling between degrees of freedom of the scalar variable $p$ due to the presence of the face integrals over the average and jump in the second of the options in \eqref{Def:B} for the definition of the bilinear form $B_\gamma$. This coupling reduces the smoothing properties of elementwise Vanka operator. For the $\mathbb Q_{r-1}$, $r\geq 2$, family of elements for the scalar variable $p$, leading to the first of the options in \eqref{Def:B}, the coupling of degrees of freedom of $p$ by its continuous in space approximation and the loss of smoothing properties  arise likewise. As a remedy, for both families of inf-sup stable approximation in \eqref{Eq:DPQ_A0} the local Vanka smoother is computed on overlapping patches of adjacent elements. In addition, an averaging of the patchwise updates is done after the Vanka smoother has been applied on all of them. 

To construct the patchwise Vanka, we let on the mesh partition $\mathcal T_l$, for $l=1,\ldots, L$, of the multigrid hierarchy the linear system, to be solved, be represented by
\begin{equation}
\label{Eq:Adbl}
\vec A_l \vec d_l = \vec b_l\,, \quad \text{for } \;\; l=1,\ldots, L\,.
\end{equation}
To each grid node $\vec \xi_l^m$, for $m = 1\,,\ldots \,, M_l$, where $M_l$ denotes the total number of grid nodes on the mesh partition $\mathcal T_l$,  we built a patch of adjacent elements such that 
\begin{equation}
\label{Eq:DefPlm}
P_l^m := \bigcup \{ K \in \mathcal T_l \mid \xi_l^m \in \overline K \}\,, \quad \text{for} \;\;m=1,\ldots, M_l \,.
\end{equation}
In two space dimensions $P_l^m$ is built from four elements, if $\xi_l^m \not\in \partial \Omega$. In three space dimensions, $P_l^m$ has eight elements, if $\xi_l^m \not\in \partial \Omega$.  If $\xi_l^m \in \partial \Omega$, patches of less elements are obtained.  On $\mathcal T_l$, let $Z_l$ denote the index set of all global degrees of freedom with cardinality $C_l$,
\begin{equation*}
C_l := \operatorname{card}(Z_l)\,.	
\end{equation*}
Let $Z_l(P_l^m)$ denote the subset of $Z_l$ all global degrees of freedom linked to the patch $P_l^m$, i.e.\ degrees of freedom of $\vec u_{\tau,h}$, $\vec v_{\tau,h}$ and $p_{\tau,h}$ for the $(k+1)$ Gauss--Radau time points of $I_n$. The cardinality of $Z_l(P_l^m)$ is denoted by $C_l^m$, 
\begin{equation*}
C_l^m : = \operatorname{card}(Z_l(P_m^l))\,, \quad \text{for} \;\; m=1,\ldots, M_l \,.
\end{equation*}
Further, we denote the index set of all local degrees of freedom on $P_l^m$  by $\hat Z_l(P_m^l) :=\{0,\ldots,C_l^m-1\}$. For the notation, we note that the set $\hat Z_l(P_m^l)$ depends on the cardinality of the patch $P_l^m$. For a given patch $P_l^m$ and a local degree of freedom with number $\hat \mu \in \hat Z_l(P_m^l) $ let the mapping 
\begin{equation}
\operatorname{dof}:\mathcal T_l\times \hat Z_l(P_m^l)  \rightarrow Z_l\,, \qquad \mu = \mathrm{dof}(P_l^m,\hat \mu)\in Z_l(P_l^m)\,,
\end{equation}
yield the uniquely defined global number $\mu \in Z_l$. Finally, we put $R = \operatorname{dim}\, \vec V_h$ and $S = \operatorname{dim}\, Q_h$. 

We are now in a position to define the local Vanka operator for the patch $P_l^m$ in an exact mathematical way; cf., e.g., \cite{AB21,HST14,J02,V86}.
\begin{defi}[Patchwise Vanka smoother]
\label{Def:VankaOP}
On patch $P_l^m$, for $m\in \{1,\ldots,M_l\}$, the local Vanka operator $\vec S_{P_l^m}: \R^{(k+1) \, \cdot \, (2R+S)} \rightarrow  \R^{C^l_m}$ is defined by  
\begin{equation}
\label{Eq:VankaOP0}
\vec S_{P_l^m}(\vec{d}) = \vec R_ {P_l^m} \vec d + \omega \, \vec A_{P_l^m}^{-1} \, \vec R_{P_l^m} \, (\vec b_l-\vec A_l \vec d)\,,
\end{equation}
with underrelaxation factor $\omega>0$ and  $P_l^m$-local restriction operator $\vec R_K:\R^{(k+1)\cdot (2R+S)} \rightarrow  \R^{C_l^m}$ defined by 
\begin{equation}
\label{Eq:DefR}
(\vec R_{P_l^m} \vec d)[\hat \mu] = \vec d[\operatorname{dof}(P_l^m,\hat \mu)]\,, \quad \text{for} \;\; \hat \mu \in \hat Z_l(P_l^m)\,.
\end{equation}
Further, the patch system matrix $\vec A_{P_l^m}\in \R^{C_l^m,C_l^m}$ in \eqref{Eq:VankaOP0} is defined by 
\begin{equation}
\vec A_{P_l^m}[\hat{\nu}][\hat{\mu}] := \vec A[\operatorname{dof}(P_l^m,\hat \nu)][\operatorname{dof}(P_l^m,\hat \mu)]\,, \quad \text{for}\;\; \hat \nu,\hat \mu \in \hat Z_l(P_l^m)\,.
\end{equation}
\end{defi}
In the numerical experiments of Sec.~\ref{Sec:NumExp} we choose $\omega = 0.7$.  The $P_l^m$-local restriction operator \eqref{Eq:DefR} assigns to a global defect vector $\vec d\in \R^{(k+1)\cdot (2R+S)} $ the local block vector $\vec R_K \vec d\in \R^{C_l^m}$ that contains all components of $\vec d$ that are associated with all degrees of freedom (for all $(k+1)$ Gauss--Radau points of $I_n$) belonging to the patch $P_l^m$.  For the computation of the inverse $(\vec A_{P_l^m})^{-1}$ of the patch system matrix$\vec A_{P_l^m}$ in \eqref{Eq:VankaOP0} we use LAPACK routines.

The application of the smoother for \eqref{Eq:Adbl} on the mesh partition $\mathcal T_l$ is summarized in Algorithm \ref{alg:smoother}.  

\begin{algorithm}[H]
	\SetAlgoLined
	Initialize solution $\vec d$ of \eqref{Eq:Adbl} with $\vec 0$\;
	\nllabel{alg:smooth_01}
	\For{$j=1,\ldots, J_{\max}$}{
    \nllabel{alg:smooth_02}
	   \For{($\mu =0$; $\mu < C_l$; $\mu\text{++})$}{
	   	  $p[\mu] = 0$\;
	      \nllabel{alg:smooth_03}
	   	  $z[\mu] = 0.0$\;
	   	  \nllabel{alg:smooth_04}
	   }
		\ForEach{patch $P_l^m$ on $\mathcal T_l$}{
	     \nllabel{alg:smooth_05}
			$\vec y = \vec S_{P_l^m}\vec(\vec d)$\; 
   	         \nllabel{alg:smooth_06}
   	         \For{ ($\hat \mu = 0$, $\hat \mu < C_l^m$; $\hat \mu\text{++}$)}{
   	         $z[\operatorname{dof}(P_l^m,\hat \mu)] \text{+=} (\vec E_{P_l^m} \vec y)[{\operatorname{dof}(P_l^m,\hat \mu)}]$\;
   	         \nllabel{alg:smooth_07}
   	         $p[\operatorname{dof}(P_l^m,\hat \mu)]\text{++}$ \;
   	         \nllabel{alg:smooth_08}
        }}
        \For{($\mu =0$; $\mu < C_l$; $\mu\text{++})$}{
        	$d[\mu] = z[\mu]/p[\mu]$\;
         \nllabel{alg:smooth_09}	
    }
}
\caption{Smoothing steps for \eqref{Eq:Adbl} on mesh partition $\mathcal T_l$.}
\label{alg:smoother}
\end{algorithm}

In line \ref{alg:smooth_01} of Alg.~\ref{alg:smoother} the defect and solution vector is pre-initialized with $\vec 0$. In line \ref{alg:smooth_02} the loop over all  $J_{\max}$ smoothing steps starts. In line \ref{alg:smooth_03}  and \ref{alg:smooth_04} the counter vector $\vec p$ for the number of updates of the degrees of freedom and the auxiliary vector $\vec z$  are initialized with $\vec 0$.	In line \ref{alg:smooth_05} the loop over all patches $P_l^m$, for $m=1,\ldots,M$, starts. In line \ref{alg:smooth_06} the local Vanka smoother is applied on patch $P_l^m$ to the current iterate $\vec d$ of the defect vector and the image is stored in the local patch vector $\vec y$. In line  \ref{alg:smooth_07} the local vector $\vec y$ is assigned to an auxiliary global  vector $\vec z$. In line \ref{alg:smooth_08} the components of the counter $\vec p$ are incremented for the (global) indices associated with the degrees of freedom of the patch $P_l^m$ processed in the loop. Finally, in line \ref{alg:smooth_09} the arithmetic mean of the local (patchwise) updates $\vec z$ is assigned to the update of the global defect vector $\vec d$. In line \ref{alg:smooth_07}, the $P_l^m$-dependent extension operator $\vec E_{P_l^m}:\R^{C_l^m} \rightarrow \R^{(k+1)\cdot (2R+S)}$ is defined by 
\begin{equation*}
(\vec E_{P_l^m} \vec y)[\mu] = \left\{\begin{array}{@{}ll} y[{\hat \mu}]\,, & \text{if}\; \exists \hat \mu \in \hat Z_l(P_l^m):\; \mu = \operatorname{dof}(K,\hat \mu)\,, \\[1ex] 0\,, & \text{if}\; \mu \not\in Z_l(P_l^m) \,.\end{array} \right.
\end{equation*}

Regarding the performance of Alg.~\ref{alg:smoother} and the overall GMRES--GMG approach we note the following. 
\begin{rem}
\begin{itemize}
\item Averaging of the patchwise updates implemented in line \ref{alg:smooth_07} and \ref{alg:smooth_09} of Alg.~\ref{alg:smoother}, that is used instead of overwriting successively the (global) degrees of freedom within the patch loop starting in line \ref{alg:smooth_05}, is essential and ensures the convergence and efficiency of the local Vanka smoother and, thereby, the desired performance of the overall GMRES--GMG linear solver. Without the averaging step we encountered convergence problems of the GMRES--GMG solver for the experiments  of Sec.~\ref{Sec:NumExp}.
\item The same holds for the application of the Vanka smoother on the patches \eqref{Eq:DefPlm} instead of using an elementwise Vanka smoother on the element $K$. The latter would lead to systems  \eqref{Eq:VankaOP0} of smaller dimension, however causes convergence problems for the GMRES--GMG solver. This is expected to be due coupling of the degrees of freedom of the scalar variable $p$ in the spatial discretizations used here,.  
\end{itemize} 
\end{rem}

\section{Numerical studies}
\label{Sec:NumExp}
 
In this section we study numerically the proposed space-time finite element and GMRES--GMG solver approach with respect to its computational and energy efficiency. Firstly, we demonstrate the accuracy of solutions in terms of convergence rates. Further, we illustrate the robustness of the GMRES--GMG solver with respect to the resolution in space and time and the polynomial approximation order. This is done for a challenging three-dimensional test problem of practical interest, in order to demonstrate the feasibility of large scale three-dimensional simulations. The numerical convergence of goal quantities of the benchmark is addressed. The parallel scaling properties of our implementation are also investigated. Beyond these studies of classical performance engineering, the energy efficiency of the approach is considered further.  

The implementation of the numerical scheme and the GMRES--GMG solver was done in an in-house high-performance frontend solver for the \texttt{deal.II} library \cite{Aetal21}. For details of the parallel implementation of the geometric multigrid solver we refer to \cite{AB21}. In all numerical experiments, the stopping criterion for the GMRES iterations is an absolute residual smaller than \mbox{1e-8}.
The computations were performed on a Linux cluster with 571 nodes, each of them with 2 CPUs and 36 cores per CPU. The CPUs are \emph{Intel Xeon Platinum 8360Y} with a base frequency of $\SI{2.4}{\giga\hertz}$, a maximum turbo frequency of $\SI{3.5}{\giga\hertz}$ and a level 3 cache of $\SI{54}{\mega\byte}$. Each node has $\SI{252}{\giga\byte}$ of main memory.

\subsection{Accuracy of the discretization:  Experimental order of convergence}
\label{Subsec:NumConv}
	
We study \eqref{Eq:HPS} for $\Omega=(0,1)^2$ and $I=(1,2]$ and the prescribed solution
\begin{equation}
	\label{Eq:givensolution}
	\boldsymbol u(\boldsymbol x, t) = \phi(\boldsymbol x, t) \boldsymbol E_2 \;\; \text{and}\;\;
	p(\boldsymbol x, t) = \phi(\boldsymbol x, t) \;\; \text{with}\;\; 
	\phi(\boldsymbol x, t) = \sin(\omega_1 t^2) \sin(\omega_2 x_1) \sin(\omega_2 x_2)
\end{equation}
and $\omega_1=\omega_2 = \pi$. We put $\rho=1.0$, $\alpha=0.9$, $c_0=0.01$ and $\boldsymbol K=\boldsymbol E_2$ with the identity $\vec E_2\in \R^{2,2}$. For the fourth order elasticity tensor $\boldsymbol C$, isotropic material properties with Young's modulus $E=100$ and Poisson's ratio $\nu=0.35$, corresponding to the Lam\'e parameters $\lambda = 86.4$ and $\mu = 37.0$, are chosen.  For an experiment with larger values of $\lambda$ and $\mu$ we refer to Table~\ref{Tab:App2} in the appendix. In our numerical experiments, the norm of $L^\infty(I;L^2)$ is approximated by computing the function's maximum value in the Gauss quadrature nodes $t_{n,m}$ of $I_n$, i.e., 
\begin{equation*}
	\label{Eq:DiscNorm}
	\| w\|_{L^\infty(I;L^2)} \approx \max \{ \| w_{|I_n}(t_{n,m})\|  \mid m=1,\ldots ,M\,, \; n=1,\ldots,N\}\,, \quad \text{with}\;\; M=100 \,.
\end{equation*}

We study the space-time convergence behavior of the scheme \eqref{Eq:DPQ_A0}. For this, the domain $\Omega$ is decomposed into a sequence of successively refined meshes of quadrilateral finite elements. The spatial and temporal mesh sizes are halved in each of the refinement steps. The step sizes of the coarsest space and time mesh are $h_0=1/(2\sqrt{2})$ and $\tau_0=0.1$. We choose the polynomial degree $k=2$ and {$r=3$}, such that discrete solutions $\vec u_{\tau,h}, \vec v_{\tau}\in Y_\tau^2(\vec V_h)$, $p_{\tau,h}\in Y_\tau^2(Q_h)$ with local spaces $\mathbb Q_3^2/ P_2^{\text{disc}}$ are obtained, as well as $k=3$ and {$r=4$} with $\vec u_{\tau,h}, \vec v_{\tau}\in Y_\tau^3(\vec V_h)$ and $p_{\tau,h}\in Y_\tau^3(Q_h)$ and local spaces $\mathbb Q_4^2/\mathbb P_3^{\text{disc}}$ ; cf.\ \eqref{Eq:DefYk}, \eqref{Def:VhQh} and \eqref{Eq:DefVhQh}. The calculated errors and corresponding experimental orders of convergence are summarized in Table~\ref{Tab:1} and \ref{Tab:2}, respectively. The error is measured in the quantities associated with the energy of the system \eqref{Eq:HPS}; cf.\ \cite[p.~15]{JR18} and \cite{BKR22}. Table~\ref{Tab:1} and \ref{Tab:2} nicely confirm the optimal rates of convergence with respect to the polynomial degrees in space and time of the overall approach. A notable difference in the convergence behavior between the pairs $\mathbb Q_r^2/\mathbb P_{r-1}^{\text{disc}}$ and $\mathbb Q_r^2/\mathbb Q_{r-1}$ of local spaces for the discretization of the spatial variables is not observed. For completeness, we summarize in Appendix~\ref{Sec:AppC} the convergence results obtained for the pair $\mathbb Q_r^2/\mathbb Q_{r-1}$ of spaces of the Taylor--Hood family. A minor superiority of the pair $\mathbb Q_r^2/\mathbb P_{r-1}^{\text{disc}}$ over the pair $\mathbb Q_r^2/\mathbb Q_{r-1}$ is only seen in the approximation of the scalar variable $p$. The coincidence of the convergence results also holds for the application of Problem~\ref{Prob:DS} instead of Problem~\ref{Prob:DSA}. The two discrete problems differ from each other by the discretization of the term $\alpha \nabla \cdot \partial_t \vec u$ in \eqref{Eq:HPS_2}. 

\begin{table}[h!t]
	\centering
	\begin{tabular}{l}
		\begin{tabular}{cccccccc}
			\toprule
			{$\tau$} & {$h$} &
			{ $\| \nabla (\vec u - \vec u_{\tau,h})  \|_{L^2(\vec L^2)} $ } & {EOC} &
			{ $\| \vec v - \vec v_{\tau,h}  \|_{L^2(\vec L^2)} $ } & {EOC} &
			{ $\| p-p_{\tau,h}  \|_{L^2(L^2)}  $ } & {EOC} 
			\\
			\cmidrule(r){1-2}
			\cmidrule(lr){3-8}
			$\tau_0/2^0$ & $h_0/2^0$ & 1.2544218392e-02 & {--} & 3.4897282317e-02  & {--} & 2.4070118274e-03 & {--} \\ 
			$\tau_0/2^1$ & $h_0/2^1$ & 1.5227995262e-03 & 3.04 & 3.9246006564e-03  & 3.15 & 2.8841669021e-04 & 3.06 \\
			$\tau_0/2^2$ & $h_0/2^2$ & 1.8904870171e-04 & 3.01 & 4.8175203148e-04  & 3.03 & 3.5986044195e-05 & 3.00 \\
			$\tau_0/2^3$ & $h_0/2^3$ & 2.3544020929e-05 & 3.01 & 5.9913957453e-05  & 3.01 & 4.4985485037e-06 & 3.00 \\
			$\tau_0/2^4$ & $h_0/2^4$ & 2.9374127099e-06 & 3.00 & 7.4788750858e-06  & 3.00 & 5.6221925444e-07 & 3.00 \\
			\bottomrule
		\end{tabular}\\
		\mbox{}\\
		\begin{tabular}{cccccccc}
			\toprule
			{$\tau$} & {$h$} &
			{ $\| \nabla (\vec u - \vec u_{\tau,h})   \|_{L^{\infty}(\vec L^2)} $ } & {EOC} &
			{ $\| \vec v - \vec v_{\tau,h}  \|_{L^{\infty}(\vec L^2)} $ } & {EOC} &
			{ $\| p-p_{\tau,h}  \|_{L^{\infty}(L^2)} $ } & {EOC} \\
			\cmidrule(r){1-2}
			\cmidrule(lr){3-8}
			$\tau_0/2^0$ & $h_0/2^0$ & 8.1391652440e-02 & {--} & 1.6415428887e-01 & {--} & 1.0828068294e-02 & {--}  \\ 
			$\tau_0/2^1$ & $h_0/2^1$ & 1.1474006637e-02 & 2.83 & 2.8042308570e-02 & 2.55 & 1.7292127803e-03 & 2.65  \\
			$\tau_0/2^2$ & $h_0/2^2$ & 1.5142126866e-03 & 2.92 & 3.6968566445e-03 & 2.92 & 2.3686467755e-04 & 2.87  \\
			$\tau_0/2^3$ & $h_0/2^3$ & 1.9337108028e-04 & 2.97 & 4.6572185327e-04 & 2.99 & 3.0596110688e-05 & 2.95  \\
			$\tau_0/2^4$ & $h_0/2^4$ & 2.4392244385e-05 & 2.99 & 5.8210228994e-05 & 3.00 & 3.8749300434e-06 & 2.98  \\
			\bottomrule
		\end{tabular}
	\end{tabular}
	\caption{%
		$L^2(L^2)$ and $L^\infty(L^2)$ errors and experimental orders of convergence (EOC) for \eqref{Eq:givensolution} with temporal polynomial degree $k=2$ and spatial degree $r=3$ for local spaces $\mathbb Q_r^2/\mathbb P_{r-1}^{\text{disc}}$.
	}
	\label{Tab:1}
\end{table}

	\begin{table}[h!t]
	\centering
	\begin{tabular}{l}
		\begin{tabular}{cccccccc}
			\toprule
			{$\tau$} & {$h$} &
			{ $\| \nabla (\vec u - \vec u_{\tau,h})  \|_{L^2(\vec L^2)} $ } & {EOC} &
			{ $\| \vec v - \vec v_{\tau,h}\|_{L^2(\vec L^2)}  $ } & {EOC} &
			{ $\| p - p_{\tau,h} \|_{L^2(L^2)}  $ } & {EOC}  \\
			\cmidrule(r){1-2}
			\cmidrule(lr){3-8}
			$\tau_0/2^0$ & $h_0/2^0$ & 2.3958455291e-03 & {--} & 1.7185653242e-02  & {--} & 7.0604463908e-04 & {--} \\ 
			$\tau_0/2^1$ & $h_0/2^1$ & 1.0529085600e-04 & 4.51 & 5.4558622568e-04  & 4.98 & 3.4799927094e-05 & 4.34 \\
			$\tau_0/2^2$ & $h_0/2^2$ & 5.9749777345e-06 & 4.14 & 1.9143193313e-05  & 4.83 & 1.2074793168e-06 & 4.85\\
			$\tau_0/2^3$ & $h_0/2^3$ & 3.6858164230e-07 & 4.02 & 9.6894400389e-07  & 4.30 & 6.1656313807e-08 & 4.29\\
			$\tau_0/2^4$ & $h_0/2^4$ & 2.2976498320e-08 & 4.00 & 5.6982903695e-08  & 4.09 & 3.6688337244e-09 & 4.07\\
			\bottomrule
		\end{tabular}\\
		\mbox{}\\
		\begin{tabular}{cccccccc}
			\toprule
			{$\tau$} & {$h$} &
			{ $\| \nabla (\vec u - \vec u_{\tau,h})  \|_{L^\infty(\vec L^2)} $ } & {EOC} &
			{ $\| \vec v - \vec v_{\tau,h}\|_{L^\infty(\vec L^2)}  $ } & {EOC} &
			{ $\| p - p_{\tau,h} \|_{L^\infty(L^2)}  $ } & {EOC}  \\
			\cmidrule(r){1-2}
			\cmidrule(lr){3-8}
			$\tau_0/2^0$ & $h_0/2^0$ & 1.3014842843e-02 & {--} & 1.1587851244e-01 & {--} & 4.9827180852e-03 & {--}  \\ 
			$\tau_0/2^1$ & $h_0/2^1$ & 8.7944275936e-04 & 3.89 & 4.1645282584e-03 & 4.80 & 2.6180886102e-04 & 4.25  \\
			$\tau_0/2^2$ & $h_0/2^2$ & 5.6086043381e-05 & 3.97 & 2.0121085807e-04 & 4.37 & 9.4019466765e-06 & 4.80  \\
			$\tau_0/2^3$ & $h_0/2^3$ & 3.5171465193e-06 & 4.00 & 1.1873567590e-05 & 4.08 & 5.0100645993e-07 & 4.23  \\
			$\tau_0/2^4$ & $h_0/2^4$ & 2.1959929657e-07 & 4.00 & 7.3034599295e-07 & 4.02 & 3.0246305328e-08 & 4.05  \\
			\bottomrule	
		\end{tabular}	
	\end{tabular}
	
	\caption{%
		$L^2(L^2)$ and $L^\infty(L^2)$  errors and experimental orders of convergence (EOC) with temporal polynomial degree $k=3$ and spatial degree $r=4$ for local spaces $\mathbb Q_r^2/\mathbb P_{r-1}^{\text{disc}}$.
	}
	\label{Tab:2}
\end{table}

Next, we show numerically that the time discretization is even superconvergent of order $2k+1$ in the discrete time nodes $t_n$, for $n=1,\ldots, N$. For this, we introduce the time mesh dependent norm
\begin{equation}
	\label{Eq:MDN}
	\| w\|_{l^\infty(L^2)} := \max\{ \| w(t_n)\| \mid n=1,\ldots, N\}\,.
\end{equation}
We prescribe the solution
\textcolor{black}{
	\begin{equation}
		\label{Eq:givensolution2}
		\boldsymbol u(\boldsymbol x, t) =  
		\begin{pmatrix}
			-2 (x-1)^2 x^2 (y-1) y (2 y-1) \sin (\omega_1  t) \\
			\phantom{-}2 (x-1) x (2 x-1) (y-1)^2 y^2 \sin (\omega_1  t)
		\end{pmatrix}
		\;\; \text{and}\;\;
		p(\boldsymbol x, t) = -2 (x-1)^2 x^2 (y-1) y (2 y-1) \sin (\omega_2  t) \;\;
	\end{equation}
}%
with $\omega_1=40\cdot\pi$ and $\omega_2 = 10\cdot\pi$. We put $\rho=1.0$, $\alpha=0.9$, $c_0=0.01$ and $\boldsymbol K=\boldsymbol E_2$ with the identity $\vec E_2\in \R^{2,2}$. For the elasticity tensor $\boldsymbol C$, isotropic material properties with Young's modulus $E=100$ and Poisson's ratio $\nu=0.35$ are used. For the local spaces we choose the pair $\mathbb Q_5^2/\mathbb Q_{4}$ such that, for any $t\in [0,T]$, the solution \eqref{Eq:givensolution2}   belongs to the discrete spaces $\vec V_h$ and $Q_h$, respectively,  and its spatial approximation is exact. This simplification is done here since we aim to study the convergence of the temporal discretization only. In the experiment, we choose $k=2$ such that discrete solutions $\vec u_{\tau,h}, \vec v_{\tau}\in Y_\tau^2(\vec V_h)$ and $p_{\tau,h}\in X_\tau^2(Q_h)$ are obtained.
\textcolor{black}{We use a spatial mesh that consists of 16 cells with $h=1/(2\sqrt{2})$ and set $\tau_0=0.02$.}
The calculated errors and corresponding experimental orders of convergence are summarized in Table~\ref{Tab:3}. Superconvergence of order $2k+1$ in the discrete time nodes is clearly observed in the second of the tables in Table~\ref{Tab:3}. 
\begin{table}[t]
	\centering
	\begin{tabular}{l}
		\begin{tabular}{cccccccc}
			\toprule
			{$\tau$} &
			{ $\| \nabla (\vec u - \vec u_{\tau,h})  \|_{L^2(\vec L^2)} $ } & {EOC} &
			{ $\| \vec v - \vec v_{\tau,h}\|_{L^2(\vec L^2)}  $ } & {EOC} &
			{ $\| p - p_{\tau,h} \|_{L^2(L^2)}  $ } & {EOC}  \\
			\cmidrule{1-1}
			\cmidrule(lr){2-7}
			$\tau_0/2^0$ & 3.3600958819e-03 & {--} & 5.2037881317e-02  & {--} & 8.5492528946e-05  & {--} \\ 
			$\tau_0/2^1$ & 3.9552579243e-04 & 3.09 & 6.6291344116e-03  & 2.97 & 4.8852611508e-06  & 4.13 \\
			$\tau_0/2^2$ & 4.8913600397e-05 & 3.02 & 8.3404665964e-04  & 2.99 & 1.7945937694e-07  & 4.77   \\
			$\tau_0/2^3$ & 6.1123116212e-06 & 3.00 & 1.0446618892e-04  & 3.00 & 1.0554654134e-08  & 4.09   \\
			$\tau_0/2^4$ & 7.6411713541e-07 & 3.00 & 1.3065222159e-05  & 3.00 & 1.1603768561e-09  & 3.19   \\
			$\tau_0/2^5$ & 9.5518196224e-08 & 3.00 & 1.6333729266e-06  & 3.00 & 1.4370494853e-10  & 3.01   \\
			$\tau_0/2^6$ & 1.1939894351e-08 & 3.00 & 2.0417851944e-07  & 3.00 & 1.7952564965e-11  & 3.00   \\
			\bottomrule	
		\end{tabular}\\
		\mbox{}\\
		
		\begin{tabular}{cccccccc}
			\toprule
			{$\tau$} &
			{ $\| \nabla (\vec u - \vec u_{\tau,h})  \|_{l^\infty(\vec L^2)} $ } & {EOC} &
			{ $\| \vec v - \vec v_{\tau,h}\|_{l^\infty(\vec L^2)}  $ } & {EOC} &
			{ $\| p - p_{\tau,h} \|_{l^\infty(L^2)}  $ } & {EOC}  \\
			\cmidrule{1-1}
			\cmidrule(lr){2-7}
			$\tau_0/2^0$ & 2.1302171198e-03 & {--} & 1.3451755096e-02 & {--} & 1.9152238265e-04 & {--}  \\ 
			$\tau_0/2^1$ & 1.1003144753e-04 & 4.28 & 6.9871664304e-04 & 4.27 & 1.0401014806e-05 & 4.20  \\
			$\tau_0/2^2$ & 3.6416424630e-06 & 4.92 & 2.3145015656e-05 & 4.92 & 3.5819518738e-07 & 4.86  \\
			$\tau_0/2^3$ & 1.1784707720e-07 & 4.95 & 7.3594357548e-07 & 4.97 & 1.1336271728e-08 & 4.98  \\
			$\tau_0/2^4$ & 3.6966657143e-09 & 4.99 & 2.3096831287e-08 & 4.99 & 3.5849645590e-10 & 4.98  \\
			$\tau_0/2^5$ & 1.1560104141e-10 & 5.00 & 7.2294482894e-10 & 5.00 & 1.1225015472e-11 & 5.00  \\
			$\tau_0/2^6$ & 3.5626018095e-12 & 5.02 & 2.2323871214e-11 & 5.02 & 3.4673336542e-13 & 5.02  \\
			\bottomrule	
		\end{tabular}
	\end{tabular}
	
	\caption{$L^2(L^2)$ and $l^\infty(L^2)$ errors and experimental orders of convergence (EOC) for \eqref{Eq:givensolution2} with temporal polynomial degree $k=2$ and spatial degree $r=5$ for local spaces $\mathbb Q_r^2/\mathbb Q_{r-1}$, showing superconvergence of order $2k+1$ in the discrete time nodes $t_n$, i.e., w.r.t.\ the norm $\|\cdot \|_{l^\infty(L^2)}$ defined in \eqref{Eq:MDN}.
	}
	\label{Tab:3}
\end{table}
	
\subsection{Computational efficiency:  Accuracy of goal quantities and convergence of the GMRES--GMG solver in a 3d benchmark}
\label{Subsec:CompEff}

\begin{figure}[t]
	\centering
	\subcaptionbox{Test setting: Cross-section plot of $\Omega$ in $(x_1,x_2)$-plane. The 3d domain extends the cross-section to $x_3\in [0,0.5]$. \label{fig:L_setting}}
	[0.34\columnwidth]
	{\includegraphics[width=0.3\textwidth,keepaspectratio]{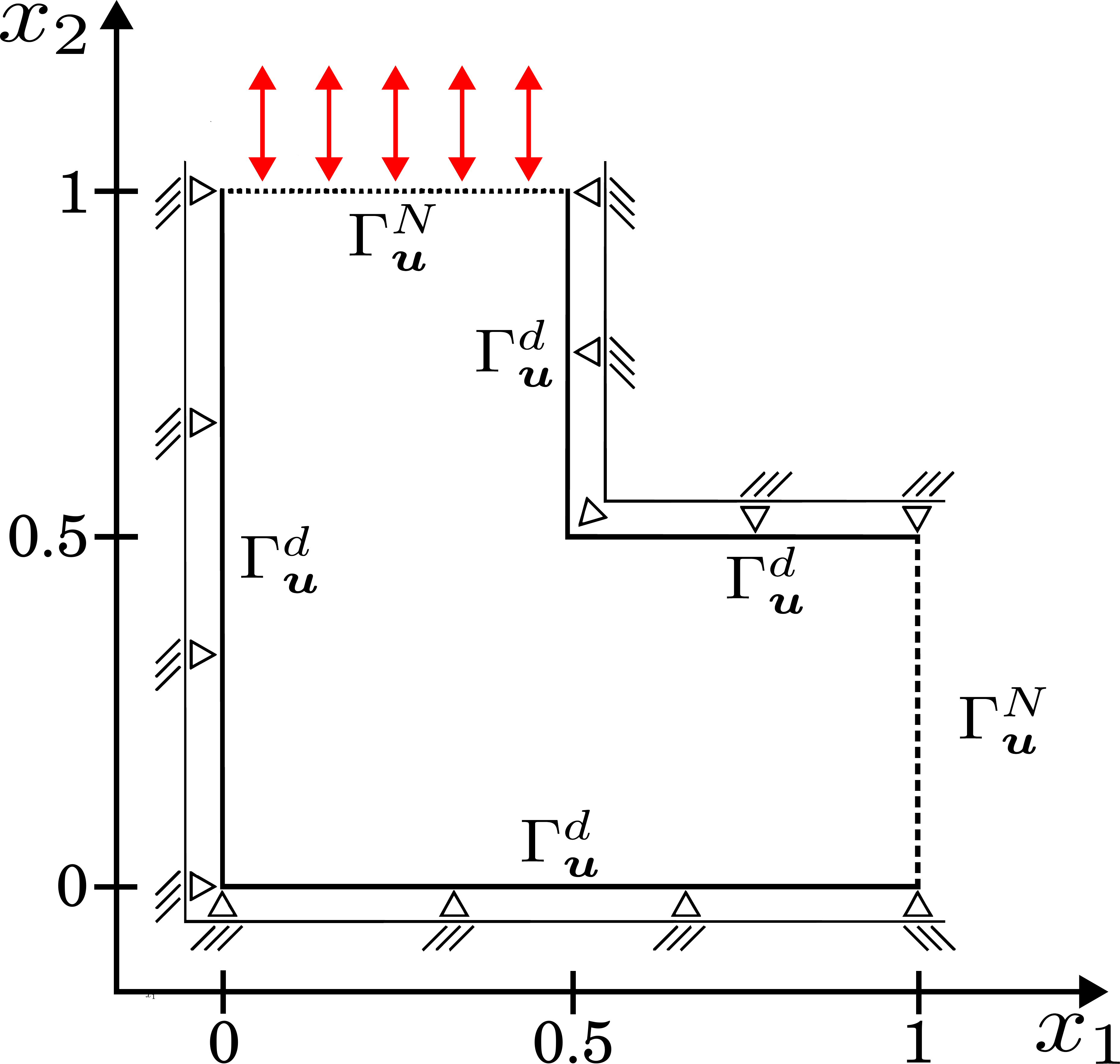}
	}
	\subcaptionbox{Modulus $|\vec u|$ of $\vec u$ at time t = 3.8. \label{fig:L_solution}}
	[0.32\columnwidth]
	{\includegraphics[width=0.38\textwidth,keepaspectratio]{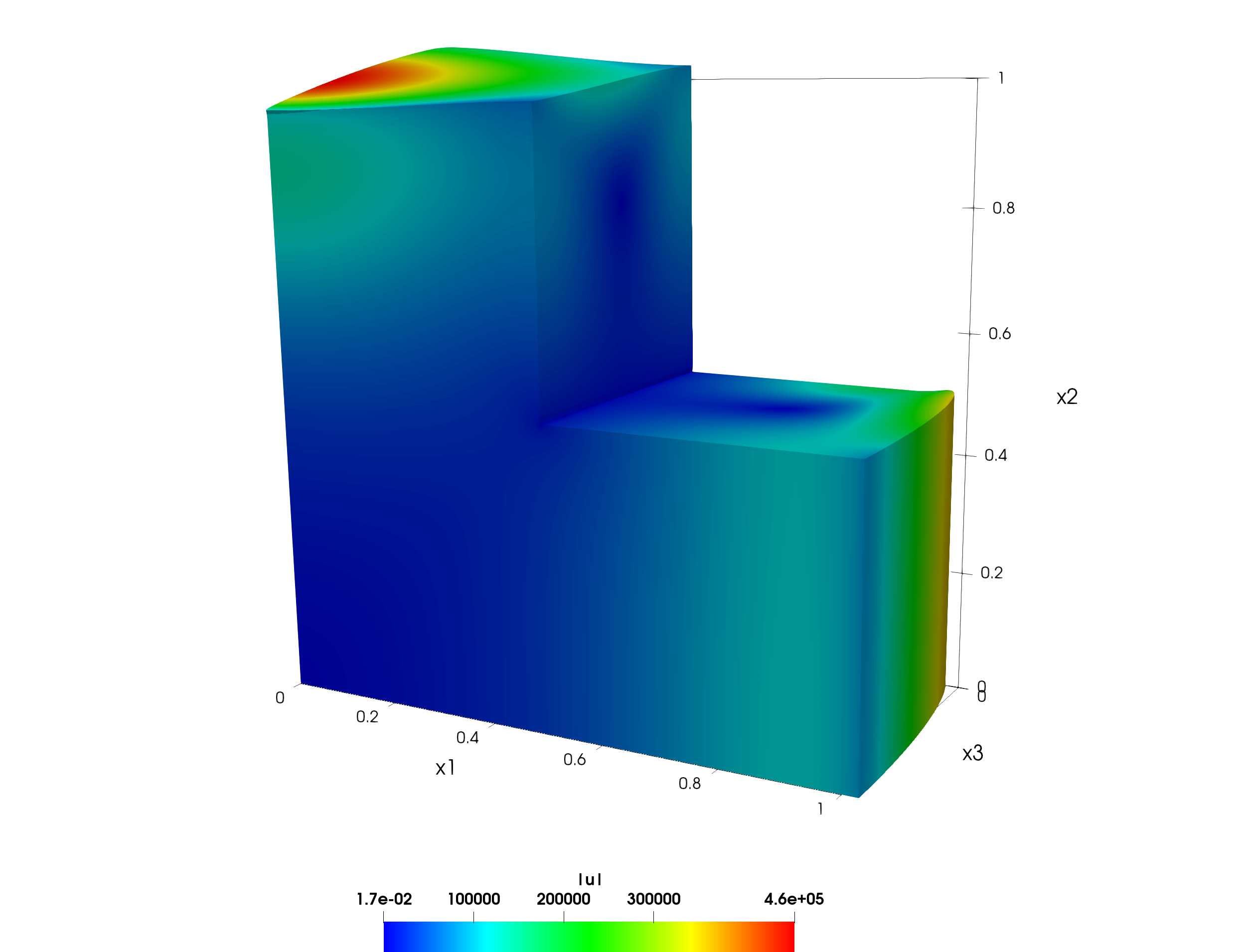}
	}
	\subcaptionbox{Variable $p$ at time t = 3.8. \label{fig:L_solution}}
	[0.32\columnwidth]
	{\includegraphics[width=0.38\textwidth,keepaspectratio]{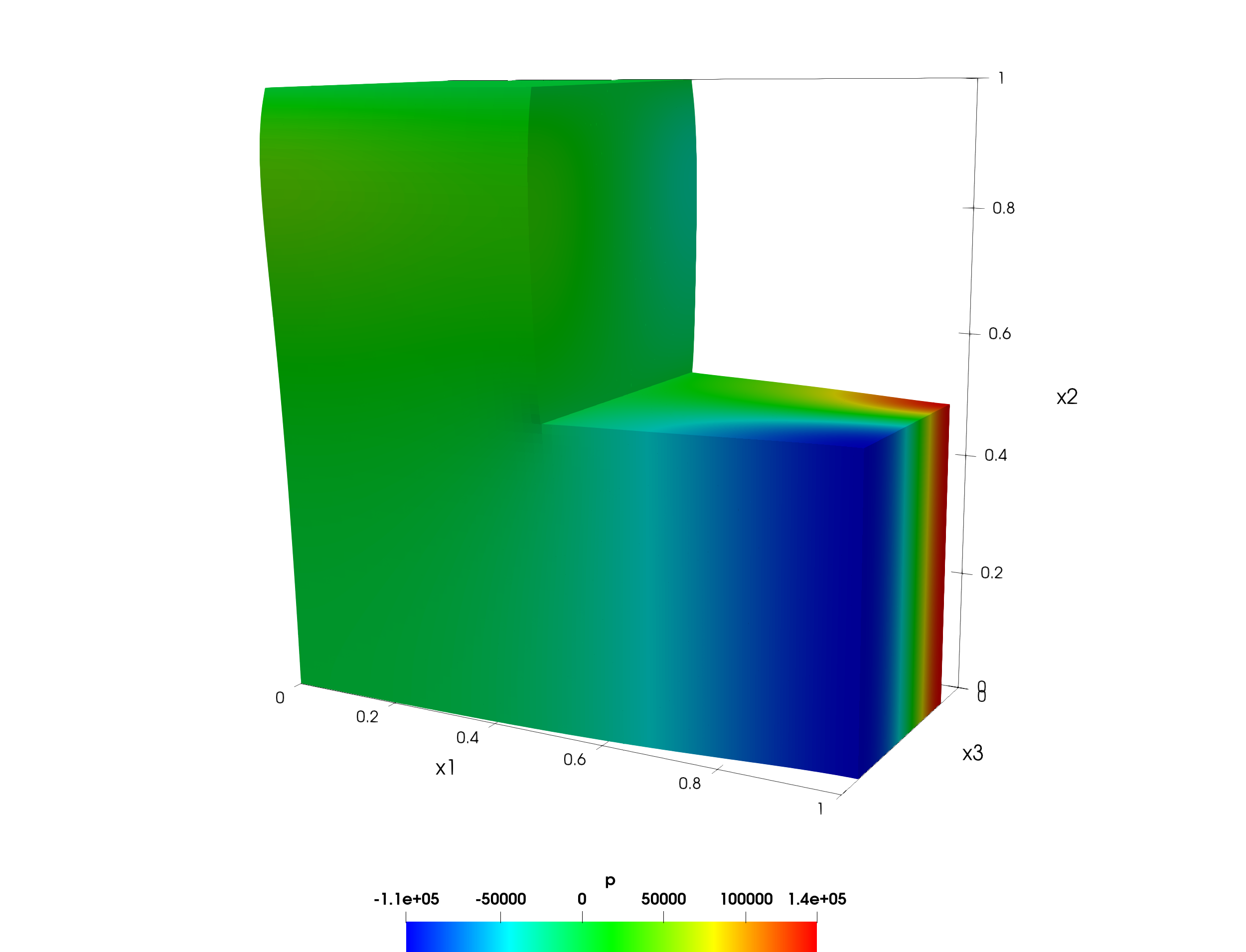}
	}
	\caption{Problem setting with boundary conditions for $\vec u$ and profile of the solution at time t = 3.8.}
	\label{fig:L_shaped_domain}
\end{figure}

For a three-dimensional problem configuration of practical interest (cf.~\cite[Sec.~5.1]{BRK17}), we study the potential of the proposed approach to compute reliably and efficiently goal quantities of physical interest that are defined in \eqref{Eq:DefGQ}. Concurrently, we analyze the convergence and performance properties of the GMRES--GMG solver of Sec.~\ref{Sec:AlgSol} for our families of space-time finite element methods. More precisely, we consider the L-shaped domain $\Omega\subset \R^3$ sketched in Fig.~\ref{fig:L_shaped_domain} along with the boundary conditions for $\vec u$ prescribed on the different parts of $\partial \Omega$. Beyond the boundary conditions \eqref{Eq:HPS_4} and \eqref{Eq:HPS_5}, we further apply (homogeneous) directional boundary conditions on parts of the boundary \textcolor{black}{for the variable $\vec{u}$}. This is done for the sake of physical realism. For the portion $\Gamma^d_{\vec u} \subset \partial \Omega$ they read as  
\begin{equation}
	\label{Eq:BCd}
	\vec u \cdot \vec n = 0 \qquad \text{and} \qquad  (\vec \sigma(\vec u)\vec n) \cdot \vec t_i  = 0\,, \quad \text{for}\;\; i=1,\ldots, d-1\,,  \quad \text{on}\;\; \Gamma^d_{\vec u} \times (0,T]\,,
\end{equation}
for the stress tensor $\vec \sigma(\vec u) = \vec C \vec \varepsilon(\vec u)$ and the unit basis vectors $\vec t_i$, for $i=1,\ldots,d-1$, of the tangent space at $\vec x \in \Gamma^d_{\vec u} $.  In the definition of $A_\gamma$ in \eqref{Def:A} and \eqref{Def:agam}, the conditions \eqref{Eq:BCd} still need to be implemented properly. For instance, for the left boundary of the L-shaped domain $\Omega$ the first of the equations in \eqref{Eq:BCd} implies that
\begin{equation*}
	\vec u_{\tau,h}\cdot \vec n = u_{\tau,h,1} = 0\,, \quad \text{for}\;\; \vec u_{\tau,h}=(u_{\tau,h,1},u_{\tau,h,2})^\top\,.
\end{equation*}
By the second of the conditions in \eqref{Eq:BCd} we get for the second of the terms on the right-hand of \eqref{Def:A} that 
\begin{equation}
	\langle \vec C \varepsilon(\vec u_{\tau,h})\vec n, \vec \chi_{\tau,h} \rangle_{\Gamma_{\vec u}^l}  = \langle \vec C \varepsilon(\vec u_{\tau,h}) \vec n \cdot \vec n,\vec n \cdot \vec \chi_{\tau,h} \rangle_{\Gamma_{\vec u}^d} \,. 
\end{equation}
For the boundary part ${\Gamma^d_{\vec u}}$ we then put, similarly to \eqref{Def:agam}, 
\begin{equation*}
	a^d_\gamma (\vec w,\vec \chi _{h}) := - \langle \vec w\cdot \vec n, \vec C \varepsilon(\vec \chi _{h}) \vec n \cdot \vec n\rangle_{{\Gamma^d_{\vec u}}} + \gamma_a\, h^{-1}  \langle \vec w\cdot \vec n, \vec \chi_h \cdot \vec n \rangle_{\Gamma^d_{\vec u}}\,, 	
\end{equation*}
while leaving $a_\gamma$ unmodified for the part ${\Gamma^D_{\vec u}}$ where Dirichlet boundary conditions are prescribed for $\vec u$. In its entirety, we thus have that $a_\gamma(\cdot,\cdot) := a_\gamma^D (\cdot,\cdot)+a_\gamma^d (\cdot,\cdot)$ with $a_\gamma^D (\cdot,\cdot)$ being defined by the right-hand side of \eqref{Def:agam}. By the arguments of Appendix~\ref{Sec:AppendixA}, the coercivity of the resulting form $A_\gamma$ is still ensured. On the right boundary we prescribe a homogeneous Neumann condition for $\vec u$ with  $\vec t_N=\vec 0$ in \eqref{Eq:HPS_5}. On the left part of the upper boundary, i.e.\ for $\Gamma^N_{\vec u}:=(0,0.5)\times \{1\} \times (0,0.5)$,  we impose the traction force (cf.\ \eqref{Eq:HPS_5})
\begin{equation*}
	\vec  t_N = 
	\begin{pmatrix}
		0 \\
		5 \cdot 10^{9} \cdot (32 x z-18 x-16 z+10) \sin (8 \pi  t) \\
		0
	\end{pmatrix}\,.
\end{equation*}
For the scalar variable $p$ we prescribe a homogeneous Dirichlet condition \eqref{Eq:HPS_6} on the left upper part of $\partial \Omega$, i.e., for $(x_1,x_2,x_3)\in \Gamma_p^D:=[0,0.5]\times \{1\}\times [0,0.5]$. On $\Gamma_p^N:=\partial \Omega \backslash \Gamma_p^D$ a homogeneous Neumann condition with $p_N=0$ in \eqref{Eq:HPS_7} is prescribed. We put $\rho=1.0$, $\alpha=0.9$, $c_0=0.01$ and $\boldsymbol K=\boldsymbol E_3$ with the identity $\vec E_3\in \R^{3,3}$. For the elasticity tensor $\boldsymbol C$, isotropic material properties with Young's modulus $E=20000$ and Poisson's ratio $\nu=0.3$, that model the material behavior of human bone (cf.\ \cite{WS00}), are used. We put $T=4$.


As goal quantities of this benchmark problem, we measure the displacement in the normal direction as well as the pressure at the right part of the boundary by 
\begin{equation}
	\label{Eq:DefGQ}
	b_{\vec{u}} = \int_{\Gamma_m} \vec{u} \cdot \vec{n} \, \ud o \qquad \text{and} \qquad 
	b_{p} = \int_{\Gamma_m} p \, \ud o \qquad \text{for}\;\; \Gamma_m:= \{1\}\times (0,0.5)\times (0,0.5)\,.
\end{equation} 

\begin{figure}[h!t]
	\centering
	\subcaptionbox{Benchmark quantity $b_{\vec{u}}$ for refinement levels  ref$_4$ to ref$_6$ and $t\in [0, T]$. \label{fig:L_bu_full}}
	[0.65\columnwidth]
	{\includegraphics[width=0.65\textwidth,keepaspectratio]{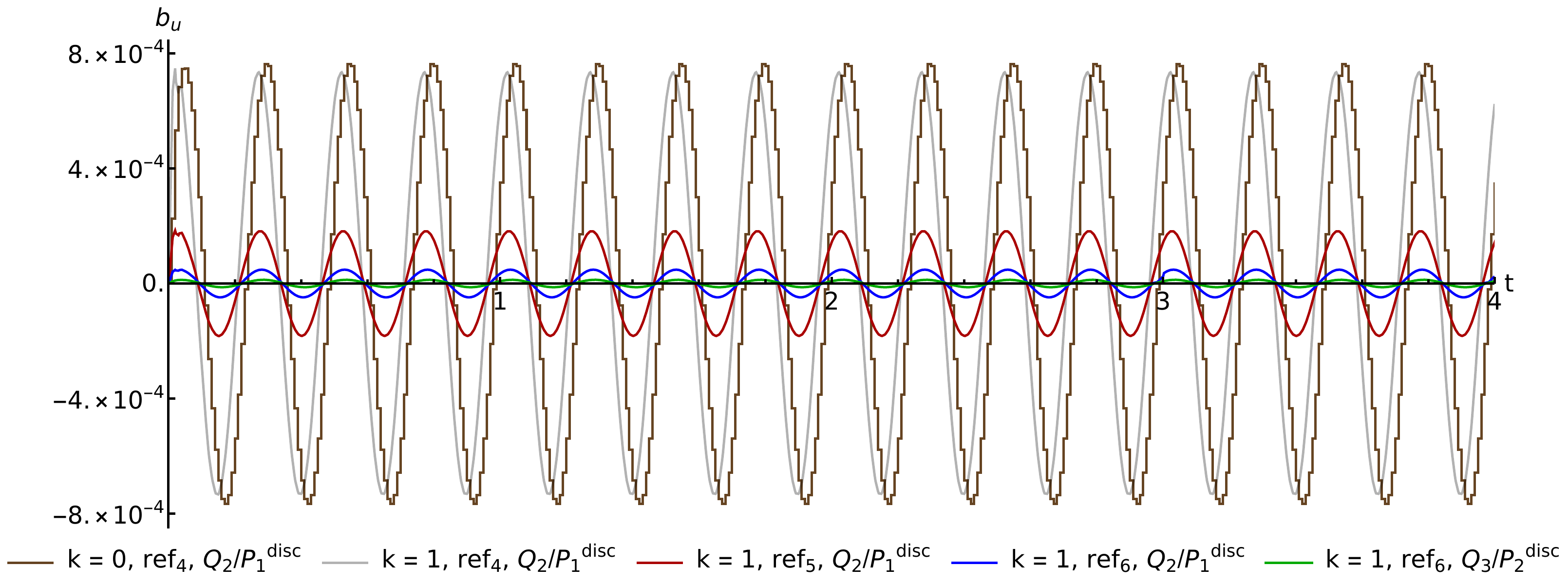}
	}\\
	\subcaptionbox{Benchmark quantity $b_{p}$ for refinement levels  ref$_4$ to ref$_6$ and $t\in [0, T]$. \label{fig:L_bp_full}}
	[0.65\columnwidth]
	{\includegraphics[width=0.65\textwidth,keepaspectratio]{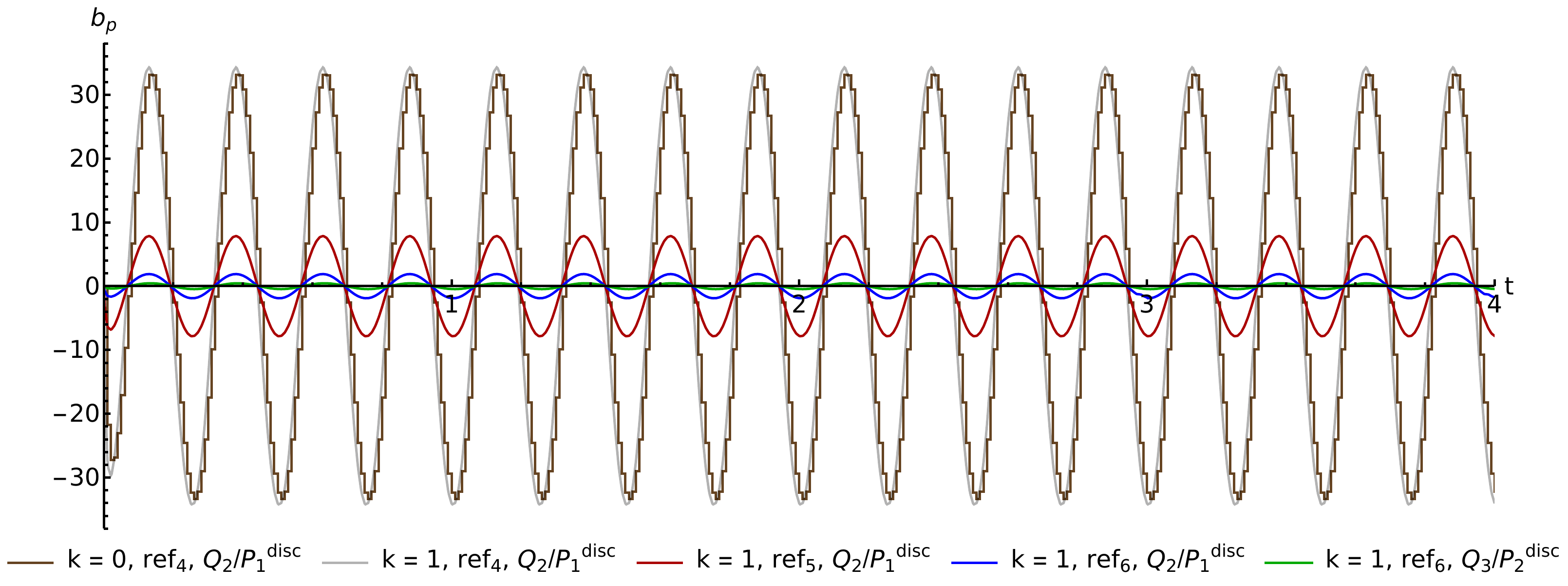}
	}
	\\[2ex]
	\subcaptionbox{Benchmark quantity $b_{\vec{u}}$ for refinement levels ref$_5$ and ref$_6$ and $t\in [0, 1]$. \label{fig:L_bu}}
	[0.65\columnwidth]
	{\includegraphics[width=0.65\textwidth,keepaspectratio]{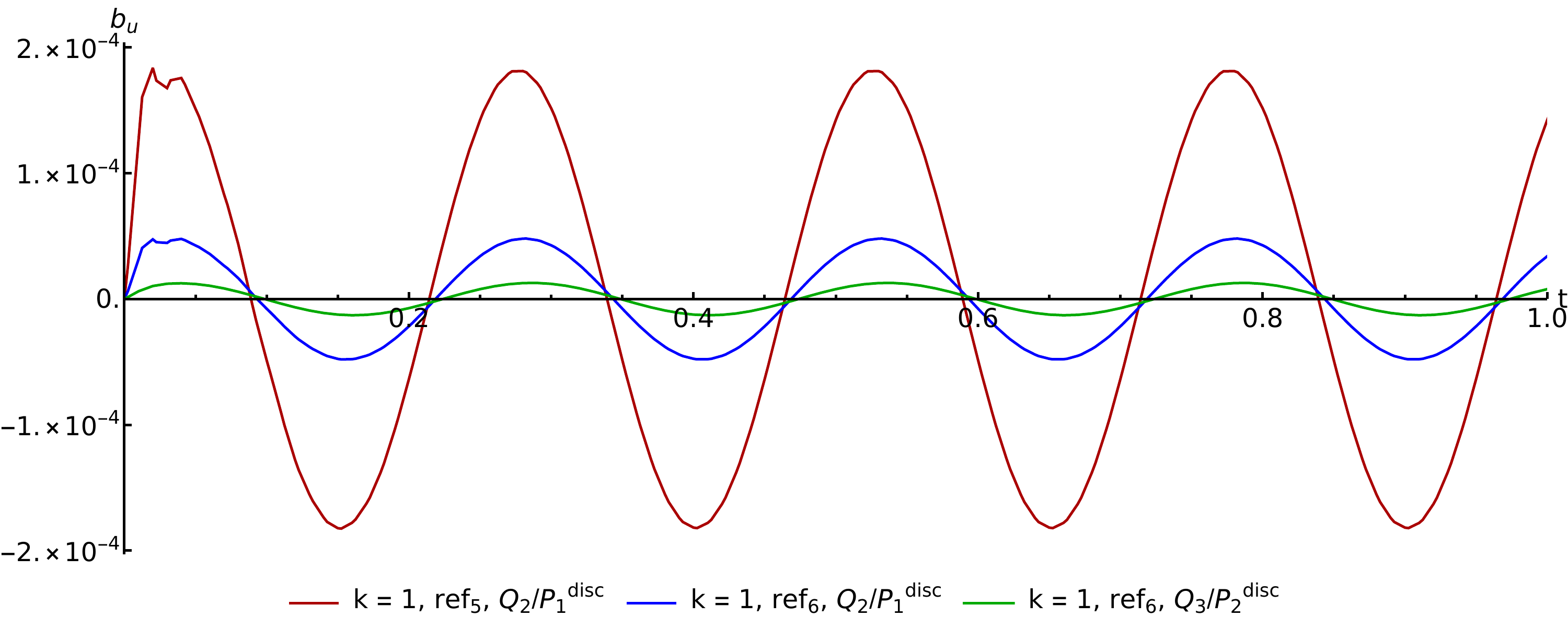}
	}\\
	\subcaptionbox{Benchmark quantity $b_{p}$ for refinement levels  ref$_5$ and ref$_6$ and $t\in [0,1]$. \label{fig:L_bp}}
	[0.65\columnwidth]
	{\includegraphics[width=0.65\textwidth,keepaspectratio]{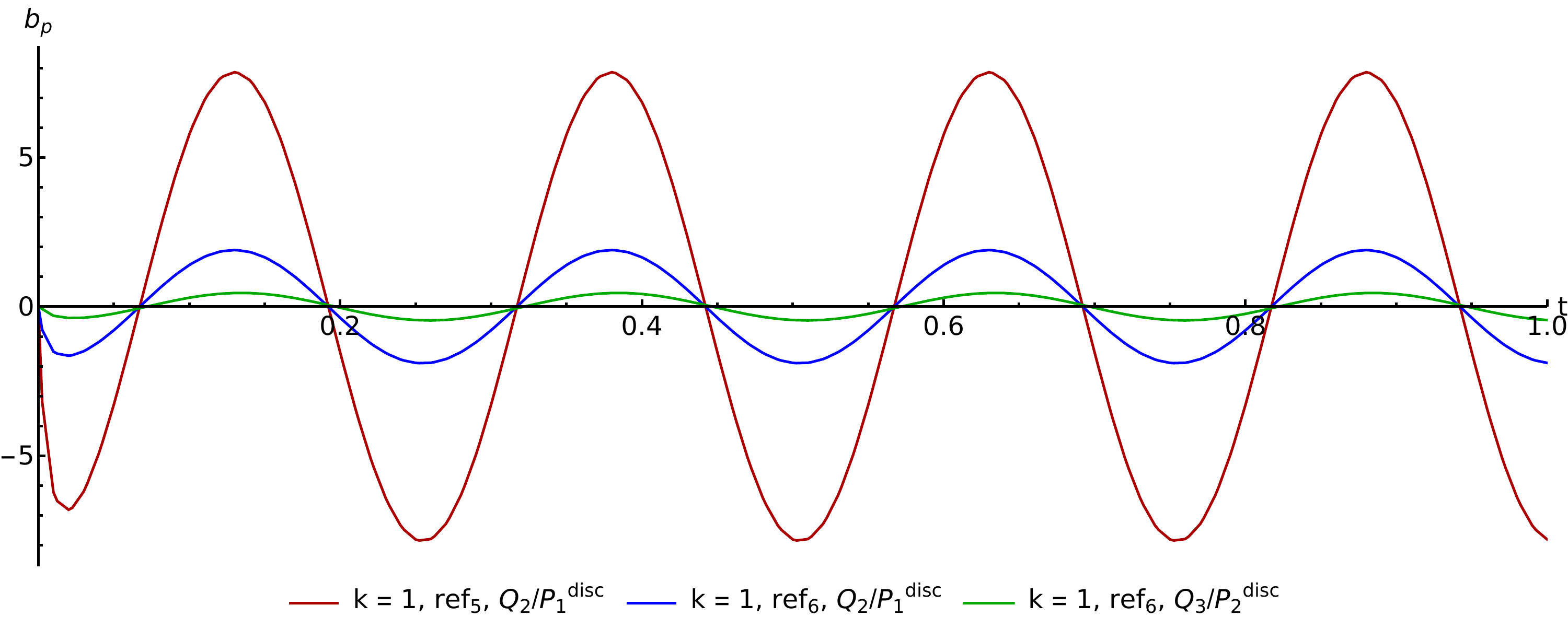}
	}
	\caption{Goal quantities $b_{\vec{u}}$ and $b_{p}$ of \eqref{Eq:DefGQ} for different discretizations (spatial mesh resolution and polynomial degrees).}
	\label{fig:L_shaped_domain_plots}
\end{figure}

The time step size of the temporal discretization is chosen as $\tau = 0.01$. The calculated profile for $t=3.8$ of the modulus of the vectorial variable $\vec u$ and the scalar variable $p$ are illustrated in Fig.~\ref{fig:L_shaped_domain}. 
In Table \ref{tab:3d_L_results} we summarize characteristic quantities and results of our simulations for various spatial resolutions and different temporal and spatial polynomial degrees of the STFEMs. In Fig.~\ref{fig:L_shaped_domain_plots} we visualize the computed benchmark quantities \eqref{Eq:DefGQ} of some of the simulations over the temporal axis. For the temporal discretization with polynomial degree $k = 0$ and $k=1$ the computed quantities differ on the spatial refinement level 4 only by a small shift of the phases of the goal quantifies. This indicates that for the current benchmark setting (choice of the model parameter and boundary conditions) the spatial discretization error dominates the temporal one. Therefore, in order to study numerically the convergence and accuracy of the goal quantities, we computed the benchmark problem for higher spatial resolutions (spatial mesh and approximation order; cf.~Table \ref{tab:3d_L_results}. The lowest row in Table \ref{tab:3d_L_results} contains the results of the finest simulation that we could run on our hardware. Fig.~\ref{fig:L_shaped_domain_plots} shows that the solution (i.e., the goal quantities) is not fully converged yet. This underlines the complexity of the benchmark. The final column of Table \ref{tab:3d_L_results} shows the convergence statistics of the proposed GMRES--GMG solver. In terms of the average number of iterations per time step the solver is (almost) grid independent. This underlines its capability and robustness for solving efficiently the complex systems arising from space-time finite element discretizations of partial differential equations.   

{
	\sisetup{scientific-notation = false,
		round-mode=places,
		round-precision=2,
		output-exponent-marker=\ensuremath{\mathrm{e}},
		table-figures-integer=4, 
		table-figures-decimal=2, 
		table-figures-exponent=1, 
		table-sign-mantissa = false, 
		table-sign-exponent = false, 
		table-number-alignment=right} 

\begin{table}[h!t]
	\caption{Computed goal quantities \eqref{Eq:DefGQ} for different temporal and spatial approximations. Here, $k$ and  $r$ are the polynomial degrees (cf.~Sec.~\ref{Sec:Not}), ref$_n$ is the number $n$ of refinements of the initial spatial grid (with $(k+1)\cdot 390$ degrees of freedom for $r=2$) and DoF$_I{_n}$ is the total number of degrees of freedom on the fine grid level $L$ of GMG per subinterval $I_n$ (i.e.\ for all $k+1$ Gauss--Radau quadrature points on $I_n$). The coarse grid level of GMG is ref$_1$ (with $(k+1) \cdot 2013$ degrees of freedom for $r=2$).  Further, the minimum and maximum of the benchmark quantities $b_{\vec u}$ and $b_p$ as well as the average number $\bar{n}_{\text{GMRES}}$ of GMRES iterations per time step are summarized.}
	\centering
	%
	\begin{tabular}{*{3}{c@{\hskip 4ex}}  r *{4}{S@{\hskip 1ex}} S}
		\toprule
		{$k$} & {$r$} & ref$_n$ &{DoF$_{I_n}$} & {$b_{\vec{u}_{\text{min}}}$} &{$b_{\vec{u}_{\text{max}}}$} & {$b_{p_{\text{min}}}$} & {$b_{p_{\text{max}}}$} & {$\bar{n}_{\text{GMRES}}$}  \\
		\cmidrule(lr){1-4} \cmidrule(r){5-6} \cmidrule(r){7-8} \cmidrule(r){9-9}
		{0} & 2 & 4 & \num{682950}   & -7.652e-04 & 7.627e-04 & -33.330 & 33.125 & \num{9} \\
		\cmidrule(lr){1-9} 
		{1} & 2 & 2 & \num{25836}   & -1.674e-02 & 1.676e-02 & -942.600 & 945.243 & \num{11} \\
		{1} & 2 & 3 & \num{182220}  & -3.263e-03 & 3.251e-03 & -165.600 & 166.085 & \num{12} \\
		{1} & 2 & 4 & \num{1365900} & -7.305e-04 & 7.335e-04 & -34.200 & 34.292 & \num{14} \\
		{1} & 2 & 5 & \num{10571532}& -1.812e-04 & 1.812e-04 &  -7.836 & 7.857 & \num{12} \\
		{1} & 2 & 6 & \num{83172876}& -4.789e-05 & 4.806e-05 &  -1.892 & 1.897 & \num{12} \\
		\cmidrule(lr){1-9} 
	    {1} & 3 & 5& \num{34596492}&  -4.518e-05 & 4.518e-05 & -1.874  & 1.879  & \num{11} \\
		{1} & 3 & 6 & \num{273640716}& -1.286e-05 & 1.280e-05 &  -0.461 & 0.459 & \num{8} \\
		\bottomrule
	\end{tabular}
	\label{tab:3d_L_results}
\end{table}
}

\subsection{Parallel scaling and energy efficiency}

Here, we study briefly the parallel scaling and energy efficiency of our solver. By studying energy consumption, we'd like to draw attention to this emerging dimension in the tuning of algorithms. Energy efficiency broadens the classical hardware-oriented numerics that is applied to enhance the performance of the current method on the target hardware and/or to find other numerical methods to improve the numerical efficiency. For the longer term, energy and power consumption needs to be mapped into a rigorous performance model. Here, we restrict ourselves to illustrate numerically the parallel scaling and energy consumption properties of our implementation that uses Message Passing Interface(MPI) libraries and multi threading parallelism. 

\begin{figure}[!h]
	\centering
	\subcaptionbox{Wall time and ideal wall time for increased number of nodes. \label{fig:wall_time}}
	[0.49\columnwidth]
	{\includegraphics[width=0.48\textwidth,keepaspectratio]{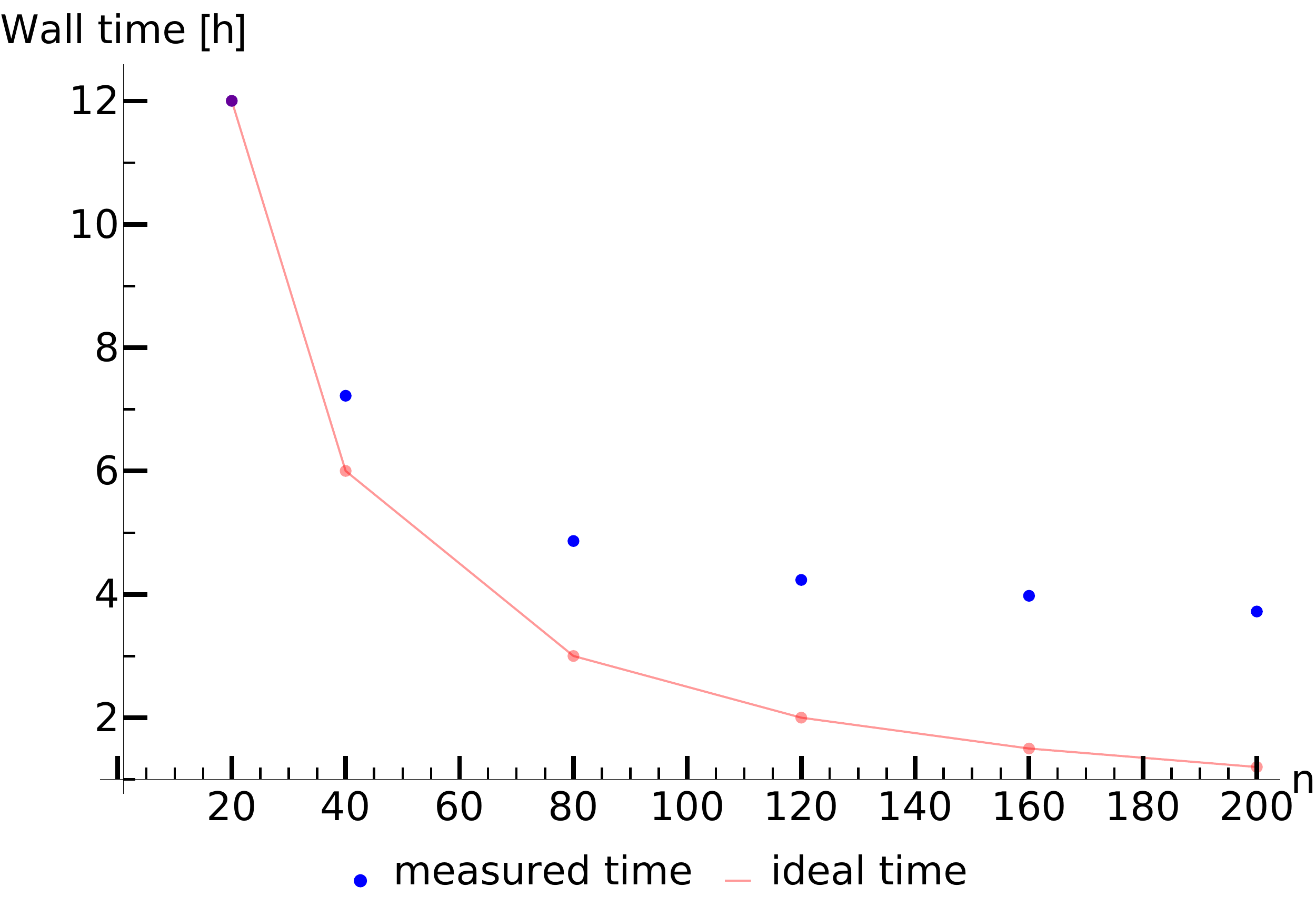}
	}
	\subcaptionbox{Parallel speedup $S$ and energy ratio $R$ for increased number of nodes. \label{fig:speedup}}
	[0.49\columnwidth]
	{\includegraphics[width=0.48\textwidth,keepaspectratio]{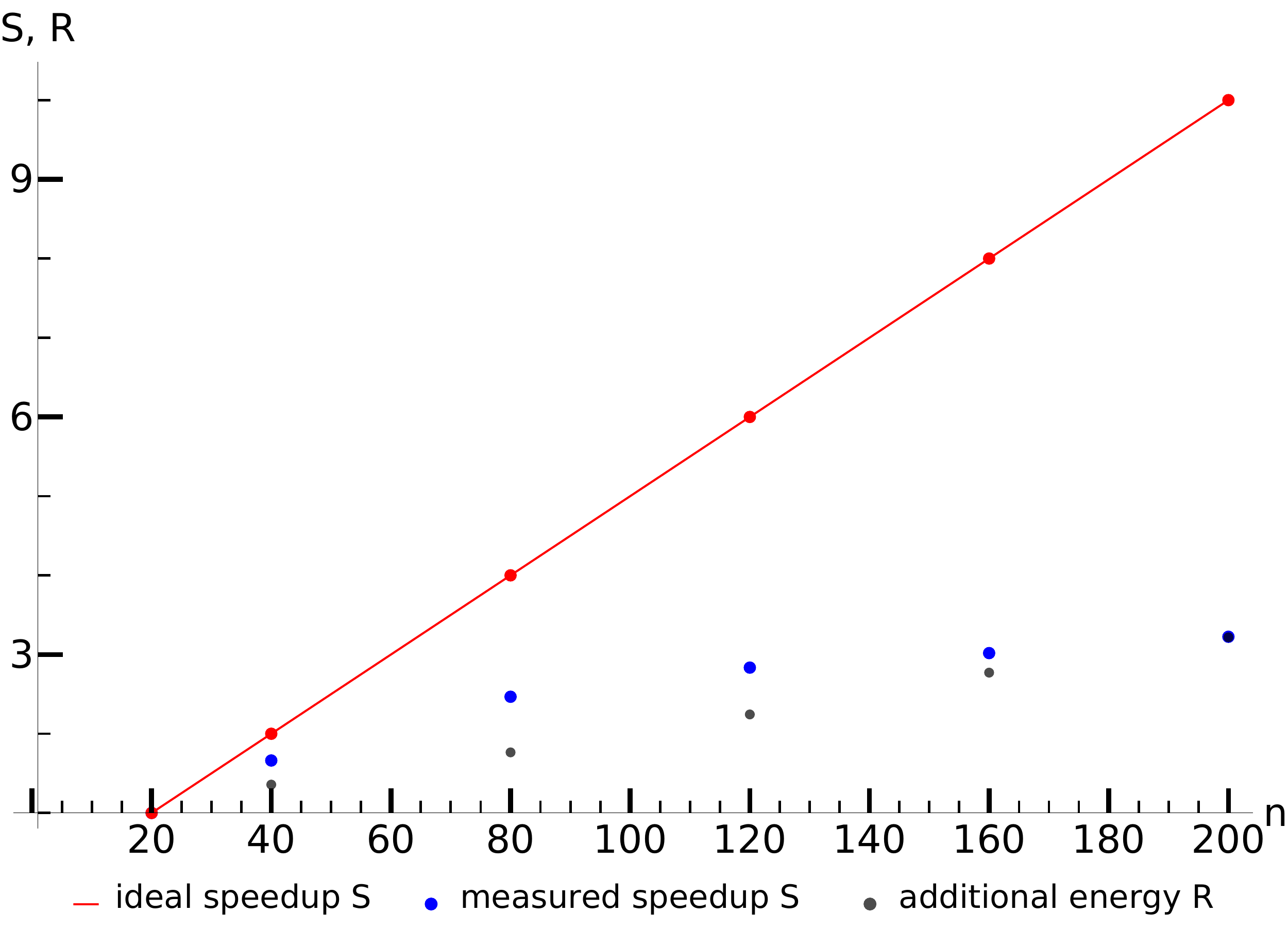}
	}
	\caption{Results of the strong scaling and energy consumption benchmark.}
	\label{fig:strong_scaling}
\end{figure}

We perform a strong scaling benchmark for the test problem of Subsec.~\ref{Subsec:CompEff} with $k = 1$, $r = 2$ and ref$_n$ = 5, with \num{10571532} degrees of freedom in each subinterval $I_n$ on the fine level $\mathcal T_L$ (cf.\ Table \ref{tab:3d_L_results}) and \num{98304} mesh cells on $\mathcal T_L$. Throughout, we assign 18 MPI processes to each of the nodes used for the computations and vary the number $n$ of nodes from $n=20$ to $n=200$. For the evaluation of the parallel scaling properties, we compute the parallel speedup of the code (cf.~\cite{A67}) that is approximated by
\begin{equation}
\label{Eq:DefCS}
S = \frac{t_{\text{wall}}(n=n_{\min})}{t_{\text{wall}}(n)}\,,
\end{equation}
where $t_{\text{wall}}(n)$ denotes the wall time of the simulation of fixed size on $n$ compute nodes and $t_{\text{wall}}(n=n_{\min})$ is the wall time of the simulation on the minimum number of nodes involved in the scaling experiment. Secondly, we compute similarly the energy ratio by means of 
\begin{equation}
\label{Eq:DefPSup}
R = \frac{E(n)}{E(n=n_{\min})}\,,
\end{equation}
where $E(n)$ measures the total energy consumption of the simulation on $n$ nodes. The energy consumption is determined by the Linux cluster workload manager slurm \cite{S23}. The energy consumption data is collected from hardware sensors using Intel's Running Average Power Limit (RAPL) mechanism. It measures the energy consumption of the processor and memory.  On our system, the sampling interval of energy consumption is determined by the value of 30 seconds. Fig.~\ref{fig:strong_scaling} illustrates the results of the performance test.

\begin{figure}[h!t]
	\begin{center}
		\includegraphics[width=0.45 \textwidth]{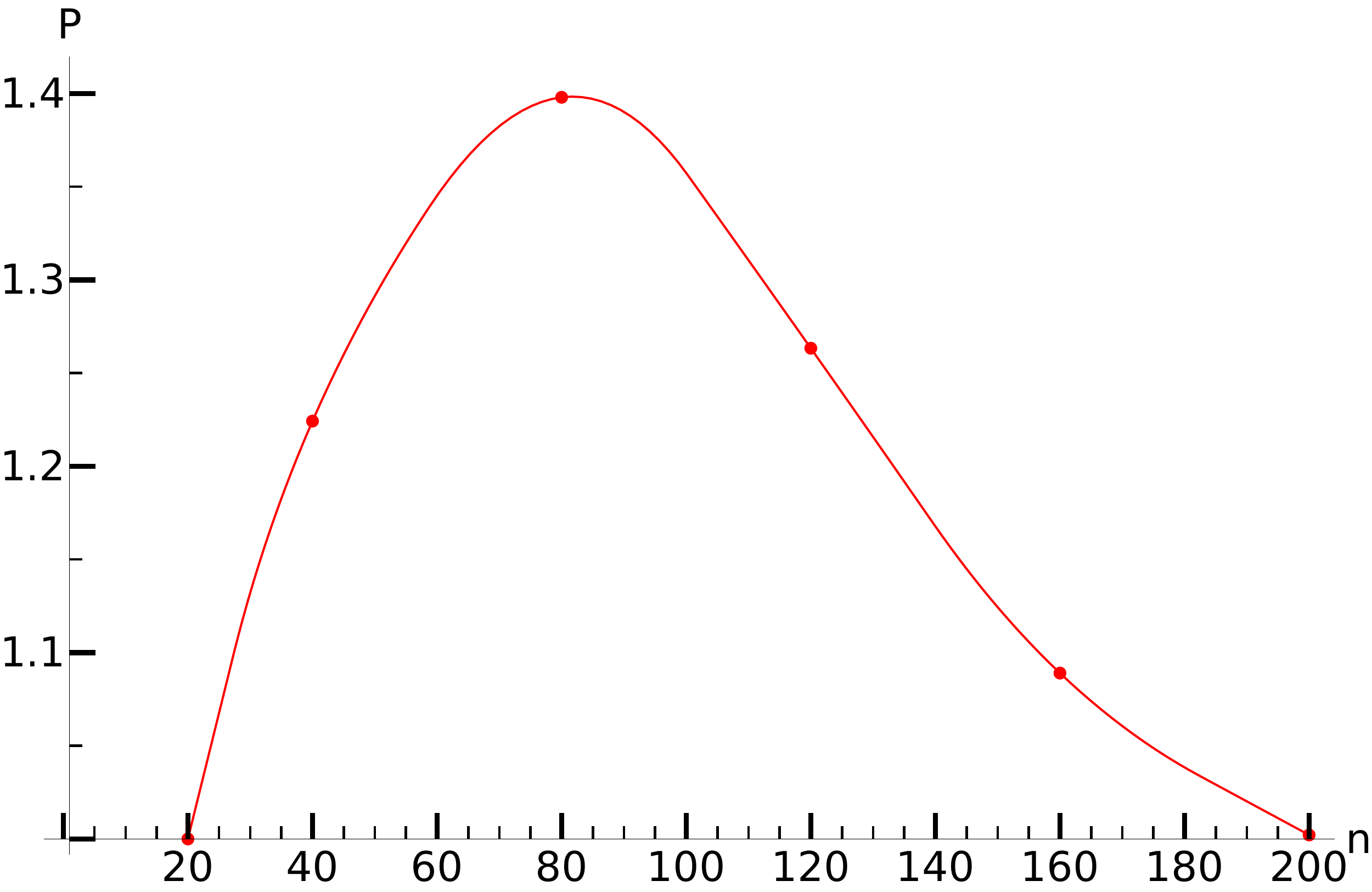}
	\end{center}
	\caption{Piecewise quadratic interpolation of the productivity function \eqref{Eq:productivity}.}	
	\label{Fig:productivity}
\end{figure}

To quantify and evaluate the productivity or resource costs of the algorithm and its implementation, we use a simple model for the $\text{Productivity }P = \frac{\text{Output }O}{\text{Input } I}$; cf.\ \cite{C00}. In an economic sense, all outputs should be the desired ones. Therefore, we use the reciprocal of the wall time $t_{\text{wall}}$ as the output $O = \frac{1}{t_{\text{wall}}}$, such that a decrease in $t_{\text{wall}}$ represents an increase of the (abstract) output. As the input $I$ we use the total energy consumption $E$ of the simulation. We scale the result by multiplying $P$ with the constant factor $E(n_{\min}) \cdot t_{\text{wall}}(n_{\min})$ such that the computation with $n = n_{\min}$ has a productivity of $P=1.0$:
\begin{equation}
	P = \frac{\frac{1}{t_{\text{wall}}(n)}}{E(n)} \cdot E(n_{\min}) \cdot t_{\text{wall}}(n_{\min}) = \frac{S}{R}\,,
	\label{Eq:productivity}
\end{equation}
with $S$ and $R$ being defined in \eqref{Eq:DefCS} and \eqref{Eq:DefPSup}, respectively. The resulting productivity curve of our computations is presented in Fig.\ \ref{Fig:productivity}. In our model the simulation on 80 compute nodes is the most productive one, that is, the ratio of output (low wall-time) to input (energy) is best. The characteristic quantities of our performance study are also summarized in Table~\ref{tab:strong_scaling}.

{
	\sisetup{scientific-notation = false,
		round-mode=places,
		round-precision=2,
		output-exponent-marker=\ensuremath{\mathrm{e}},
		table-figures-integer=3, 
		table-figures-decimal=3, 
		table-figures-exponent=0, 
		table-sign-mantissa = false, 
		table-sign-exponent = false, 
		table-number-alignment=right} 

\begin{table}[h!t]
	\caption{Computed quantities of the scaling and performance benchmark.}
	\centering
	\begin{tabular}{c *{2}{S} *{3}{c}  }
		\toprule
		{$n$} & {$t_{\text{wall}}$ [h]} & {E [kWh]} & {$S$} &{$R$} & {$P$}  \\
		\cmidrule(r){1-1} \cmidrule(r){2-2} \cmidrule(r){3-3} \cmidrule(r){4-4} \cmidrule(r){5-5} \cmidrule{6-6}
		{20} & 12.00 & 98.31 & 1.00   & 1.00 & 1.00  \\
		{40} & 7.22 & 133.49 & 1.66   & 1.36 & 1.22  \\
		{80} & 4.87 & 173.50 & 2.47   & 1.76 & 1.40  \\
		{120} & 4.23 & 220.61 & 2.83  & 2.24 & 1.26  \\
		{160} & 3.98 & 272.48 & 3.02  & 2.77 & 1.09  \\
		{200} & 3.72 & 316.34 & 3.22  & 3.22 & 1.00  \\
		\bottomrule
	\end{tabular}
	\label{tab:strong_scaling}
\end{table}
}

\section{Summary and outlook}
\label{Sec:SumOut}

In this work we presented and analyzed families of space-time finite element discretizations of the coupled hyperbolic-parabolic system \eqref{Eq:HPS} of poro- or thermoelasticity. The time  discretization uses the discontinuous Galerkin method. The space discretization is based on inf-sup stable pairs of finite element spaces with continuous and discontinuous approximation of the scalar variable $p$. Well-posedness of the discrete problems is proved. For efficiently solving the arising algebraic systems with complex block structure in the case of increasing polynomial degrees of the time discretization a geometric multigrid preconditioner with a local Vanka smoother on patches of finite elements is proposed and studied. The overall approach is carefully evaluated for a sequence of numerical experiments, in particular for a challenging three-dimensional test setting of practical interest. Parallel scaling and energy consumption is also considered in the experiments. The benchmark problem of Subsec.~\ref{Subsec:CompEff} deserves further studies. In particular, soft materials with a dominating impact of the hyperbolic equation \eqref{Eq:HPS_1} to the overall system  \eqref{Eq:HPS} need to be studied as well. The definition of the goal quantities of the proposed benchmark might be reconsidered and revised with regard to either an improved choice of such quantities or their computation with higher rate of convergence. This evaluation remains as a work for the future. Finally, multi-field formulations of \eqref{Eq:HPS} with an explicit approximation of the stress tensor $\vec \sigma = \vec C \vec \varepsilon$ and the flux vector $\vec q=-\vec K \nabla p$ might be advantageous for applications of \eqref{Eq:HPS} in that these quantities are of interest or desired. By using the concept of the local (patchwise) Vanka smoother GMG preconditioner, that was presented here, multi-field formulations of \eqref{Eq:HPS} are expected to be become feasible as well. This also remains as a work for the future. 

\section*{Acknowledgements}

Computational resources (HPC-cluster HSUper) have been provided by the project hpc.bw, funded by dtec.bw — Digitalization and Technology Research Center of the Bundeswehr.

\appendix
\def\appendixname{}

\section*{Appendix}
\section{Discrete coercivity}
\label{Sec:AppendixA}

Here, we prove the discrete coercivity of the bilinear form $A_\gamma$ defined in \eqref{Def:A} along with \eqref{Def:agam}. For this we introduce the mesh- and parameter-dependent norm
\begin{equation}
\label{Eq:Mnorm}
\| \vec v_h \|_{h,\widetilde \gamma_a} := \Big(\| \varepsilon(\vec v_h) \|^2 + \widetilde \gamma_a \, h^{-1}\|\vec v_h \|_{\Gamma^D_{\vec u}}^2\Big)^{1/2}\,, \quad \text{for}\;\; \vec v_h \in \vec V_h\,,
\end{equation}
with some parameter $\widetilde \gamma_a> 0$. Its choice is addressed below. The norm property of \eqref{Eq:Mnorm} is ensured by a variant of Korn's inequality; cf.~\cite[Eq.~(1.19)]{B03}. We recall the well-known inverse inequality (cf.~\cite[p.~28]{T06})
\begin{equation}
\label{Eq:Trace}
\| \vec \varepsilon(\vec v_h) \vec n \|_{\Gamma^D_{\vec u}} \leq c h^{-1/2} \| \vec \varepsilon (\vec v_h) \|\,, \quad \text{for}\;\; \vec v_h \in \vec V_h\,.	
\end{equation}
From \eqref{Def:A} along with \eqref{Def:agam}, it follows by the positive definiteness \eqref{Eq:PosDefC} of $\vec C$, the inequality of Cauchy--Young with a sufficiently small constant $\delta >0$ and the trace inequality \eqref{Eq:Trace} that, for all $\vec v_h \in \vec V_h$, 
\begin{align*}
	\nonumber
	A_\gamma( \vec v_h, \vec v_h) & = \langle \vec C \vec \varepsilon (\vec v_h),\vec \varepsilon (\vec v_h)\rangle - \langle \vec C\vec \varepsilon(\vec v_h) \vec n, \vec v_h\rangle_{\Gamma^D_{\vec u}} - \langle \vec v_h, \vec C\vec \varepsilon(\vec v_h) \vec n \rangle_{\Gamma^D_{\vec u}}  + \gamma_a \,h^{-1} \langle \vec v_h, \vec v_h \rangle_{\Gamma^D_{\vec u}}\\[1ex]
	\nonumber
	& \geq c \|\vec \varepsilon(\vec v_h) \|^2 - c \| \vec \varepsilon(\vec v_h) \vec n \|_{\Gamma^D_{\vec u}} \| \vec v_h\|_{\Gamma^D_{\vec u}}+ \gamma_a \,h^{-1} \| \vec v_h \|^2_{\Gamma^D_{\vec u}}\\[1ex]
	\nonumber
	& \geq c \| \vec \varepsilon(\vec v_h) \|^2  - c \delta \, h \| \vec \varepsilon(\vec v_h) \vec n \|^2_{\Gamma^D_{\vec u}}  - c \delta^{-1} h^{-1}\| \vec v_h\|^2_{\Gamma^D_{\vec u}} + \gamma_a \,h^{-1} \| \vec v_h \|^2_{\Gamma^D_{\vec u}}\\[1ex]
	\nonumber
	& \geq c \| \vec \varepsilon(\vec v_h) \|^2  - c \delta \| \vec \varepsilon (\vec v_h) \|^2 - c \delta^{-1} h^{-1}\| \vec v_h\|^2_{\Gamma^D_{\vec u}} + \gamma_a \,h^{-1} \| \vec v_h \|^2_{\Gamma^D_{\vec u}}\\[1ex]
	\label{Eq:ExUniDS_4}
	& = c\Big( \| \vec \varepsilon (\vec v_h) \|^2 + \widetilde \gamma_a \,h^{-1} \| \vec v_h \|^2_{\Gamma^D_{\vec u}}\Big)\,,
\end{align*}	
with some constant $\widetilde \gamma_a >0$ for a sufficiently large choice of the algorithmic parameter $\gamma_a$ in \eqref{Def:agam}, such that $\widetilde \gamma_a :=\gamma_a -c\delta^{-1}>0$. Thus, there holds for some constant $c>0$ that 
\begin{equation}
\label{Eq:CoercA}
A_\gamma( \vec v_h, \vec v_h) \geq c \| \vec v_h \|_{h,\widetilde \gamma_a}\,, \quad \text{for all}\;\; \vec v_h \in \vec V_h\,.	
\end{equation}

Secondly, we prove the discrete coercivity of the bilinear form $B_\gamma$ defined in \eqref{Def:B} along with \eqref{Def:bgam}. For brevity, we study the case $Q_h= Q_h^{l,\text{cont}}$ only, such that $Q_h\subset H^1(\Omega)$ is satisfied. The case $Q_h= Q_h^{l,\text{disc}}$ can be captured similarly. We introduce the mesh- and parameter-dependent norm
\begin{equation*}
	\|  q_h \|_{h,\widetilde \gamma_b} := \Big(\|\nabla q_h \|^2 + \widetilde \gamma_b \, h^{-1}\| q_h \|_{\Gamma^D_{p}}^2\Big)^{1/2}\,, \quad \text{for}\;\; q_h \in Q_h\subset H^1(\Omega)\,,
\end{equation*}
with some parameter $\widetilde \gamma_b> 0$, where $Q_h=Q_h^{l,\text{cont}}$; cf.~\eqref{Def:VhQh}. Further, we recall the well-known inverse inequality
\begin{equation}
	\label{Eq:Trace2}
	\| \nabla q_h \cdot \vec n \|_{\Gamma^D_{p}} \leq c h^{-1/2} \| \nabla q_h \|\,, \quad \text{for}\;\; q_h \in Q_h\,.	
\end{equation}
From \eqref{Def:B} along with \eqref{Def:bgam}, it follows by the positive definiteness \eqref{Eq:PosDefK} of $\vec K$, the inequality of Cauchy--Young with a sufficiently small constant $\delta >0$ and the trace inequality \eqref{Eq:Trace2} that, for all $q_h \in Q_h$, 
\begin{align*}
	\nonumber
	B_\gamma( q_h, q_h) & = \langle \vec K \nabla q_h ,\nabla q_h \rangle - \langle \vec K\vec \nabla q_h \cdot \vec n, q_h\rangle_{\Gamma^D_{p}} - \langle q_h, \vec K\nabla q_h \cdot \vec n \rangle_{\Gamma^D_{p}}  + \gamma_b \,h^{-1} \langle q_h, q_h \rangle_{\Gamma^D_{p}}\\[1ex]
	\nonumber
	& \geq c \|\nabla q_h) \|^2 - c \| \nabla q_h \cdot \vec n \|_{\Gamma^D_{p}} \| q_h\|_{\Gamma^D_{p}}+ \gamma_b \,h^{-1} \| q_h \|^2_{\Gamma^D_{p}}\\[1ex]
	\nonumber
	& \geq c \| \nabla q_h \|^2  - c \delta \, h \| \nabla q_h \cdot \vec n \|^2_{\Gamma^D_{p}}  - c \delta^{-1} h^{-1}\| q_h\|^2_{\Gamma^D_{p}} + \gamma_b \,h^{-1} \| q_h \|^2_{\Gamma^D_{p}}\\[1ex]
	\nonumber
	& \geq c \| \nabla q_h \|^2  - c \delta \| \nabla q_h \|^2 - c \delta^{-1} h^{-1}\| \vec v_h\|^2_{\Gamma^D_{p}} + \gamma_b \,h^{-1} \| q_h \|^2_{\Gamma^D_{p}}\\[1ex]
	\label{Eq:ExUniDS_5}
	& = c\Big( \| \nabla q_h \|^2 + \widetilde \gamma_b \,h^{-1} \| q_h \|^2_{\Gamma^D_{p}}\Big)\,,
\end{align*}	
with some constant $\widetilde \gamma_b >0$ for a sufficiently large choice of the algorithmic parameter $\gamma_b$ in \eqref{Def:bgam}, such that $\widetilde \gamma_b :=\gamma_a -c\delta^{-1}>0$. Thus, there holds for some constant $c>0$ that 
\begin{equation}
	\label{Eq:CoercB}
	B_\gamma( q_h, q_h) \geq c \| q_h \|_{h,\widetilde \gamma_b}\,, \quad \text{for all}\;\; q_h \in Q_h\,.	
\end{equation}

\section{Alternative  formulation of the fully discrete problem}
\label{Sec:AppendixB}

Here, we present an alternative formulation for the fully discrete system to \eqref{Eq:HPS}. The difference to Problem~\ref{Prob:DSA} comes through using $\partial_t \vec u_{\tau,h}$ instead of $\vec v_{\tau,h}$ in \eqref{Eq:DPQ_A3} by means of \eqref{Eq:HPS_8} and \eqref{Eq:DPQ_A1}, respectively. However, this modification requires the inclusion of an additional boundary integral in the resulting equation \eqref{Eq:DPQ_3}. This non-obvious adaptation results from the proof of well-posedness of the discrete problem. The following problem is thus considered. 

\begin{prob}[Numerically integrated $I_n$-problem with $\partial_t \vec u_{\tau,h}$]
	\label{Prob:DS}
	For given $\vec u_{h}^{n-1}:= \vec u_{\tau,h}(t_{n-1})\in \vec V_h$, $\vec v_{h}^{n-1}:= \vec v_{\tau,h}(t_{n-1})\in \vec V_h$,  and $p_{h}^{n-1}:= p_{\tau,h}(t_{n-1}) \in Q_h$ with  $\vec u_{\tau,h}(t_0) :=\vec u_{0,h}$, $\vec v_{\tau,h}(t_0) :=\vec u_{1,h}$ and $p_{\tau,h}(t_0) := p_{0,h}$, find $(\vec u_{\tau,h},\vec v_{\tau,h},p_{\tau,h}) \in \mathbb P_k (I_n;\vec V_h) \times \mathbb P_k (I_n;\vec V_h) \times \mathbb P_k (I_n;Q_h)$ such that
	\begin{subequations}
		\label{Eq:DPQ_0}
		\begin{alignat}{2}
			\label{Eq:DPQ_1}
			&\begin{aligned}
				&Q_n \big(\langle \partial_t \vec u_{\tau,h} , \vec \phi_{\tau,h} \rangle  - \langle \vec v_{\tau,h} , \vec \phi_{\tau,h} \rangle \big) + \langle \vec u^+_{\tau,h}(t_{n-1}), \vec \phi_{\tau,h}^+(t_{n-1})\rangle =  \langle \vec u_{h}^{n-1}, \vec \phi_{\tau,h}^+(t_{n-1})\rangle  \,,\\[1ex]
			\end{aligned}\\
			\label{Eq:DPQ_2}
			&\begin{aligned}
				& Q_n \Big(\langle \rho \partial_t \vec v_{\tau,h} , \vec \chi_{\tau,h} \rangle + A_\gamma(\vec u_{\tau,h}, \vec \chi_{\tau,h} ) + C(\vec \chi_{\tau,h},p_{\tau,h})\Big) + \langle \rho \vec v^+_{\tau,h}(t_{n-1}), \vec \chi_{\tau,h}^+(t_{n-1})\rangle \\[1ex]
				& \qquad = Q_n \Big(F_\gamma(\vec \chi_{\tau,h})\Big) + \langle \rho \vec v_{h}^{n-1}, \vec \chi_{\tau,h}^+(t_{n-1})\rangle \,,
			\end{aligned}\\[1ex]
			\label{Eq:DPQ_3}
			&\begin{aligned}
				&Q_n \Big(\langle c_0 \partial_t p_{\tau,h},\psi_{\tau,h} \rangle  - C(\partial_t \vec u_{\tau,h},\psi_{\tau,h})+ B_\gamma (p_{\tau,h}, \psi_{\tau,h})\Big)\\[1ex]
				& \qquad \quad + \langle c_0 p^+_{\tau,h}(t_{n-1}) +\alpha \nabla \cdot \vec u^+_{\tau,h}(t_{n-1}),  \psi_{\tau,h}^+(t_{n-1})\rangle - \alpha \langle \vec u^+_{\tau,h}(t_{n-1}) \cdot \vec n ,\psi_{\tau,h}^+(t_{n-1}) \rangle_{\Gamma^D_{\vec u}} \\[1ex] 
				& \qquad = Q_n \Big( G_\gamma(\psi_{\tau,h})\Big) + \langle c_0 p_{h}^{n-1} + \alpha \nabla \cdot  \vec u_{h}^{n-1}, \psi_{\tau,h}^+(t_{n-1})\rangle - \alpha \langle \vec u_D(t_{n-1}) \cdot \vec n ,\psi_{\tau,h}^+(t_{n-1}) \rangle_{\Gamma^D_{\vec u}} 
			\end{aligned}
		\end{alignat}
	\end{subequations}
	for all $(\vec \phi_{\tau,h},\vec \chi_{\tau,h},\psi_{\tau,h})\in  \mathbb P_k (I_n;\vec V_h) \times \mathbb P_k (I_n;\vec V_h) \times \mathbb P_k (I_n;Q_h)$.
\end{prob}

\begin{lem}[Existence and uniqueness of solutions to Problem \ref{Prob:DS}]
	\label{Lem:ExUniDS2}
	Problem \ref{Eq:DPQ_0} admits a unique solution.
\end{lem}

\begin{mproof}
The proof follows basically the lines of the proof of Lem.~\ref{Lem:ExUniDS} and is kept short. Only differences to Lem.~\ref{Lem:ExUniDS} are depicted. For the differences $(\vec u_{\tau,h},\vec v_{\tau,h},p_{\tau,h})$ of two triples satisfying \eqref{Eq:DPQ_0} and the test functions $\vec \phi_{\tau,h}= A_\gamma \vec u_{\tau,h}$, $\vec \chi_{\tau,h}=\vec v_{\tau,h}$ and $\psi_{˝\tau,h}=p_{\tau}$ there holds that 
\begin{equation}
\label{Eq:ExUniDS_Z1}
\begin{aligned}
& \frac{1}{2} \int_{t_{n-1}}^{t_n} \frac{d}{dt} \big( A_\gamma (\vec u_{\tau,h} , \vec u_{\tau,h}) + \langle \rho \vec v_{\tau,h} , \vec v_{\tau,h}  \rangle  + \langle c_0 p_{\tau,h}, p_{\tau,h} \rangle \big ) \ud t + Q_n\big(B_\gamma(p_{\tau,h},p_{\tau,h}) \big)  \\[1ex] 
& \quad     + Q_n\big(C_\gamma (\vec v_{\tau,h} - \partial_t \vec u_{\tau,h},p_{\tau,h}) \big) + \alpha \langle \nabla \cdot \vec u^+_{\tau,h}(t_{n-1}),  p_{\tau,h}^+(t_{n-1})\rangle {\color{black} - \alpha \langle \vec u^+_{\tau,h}(t_{n-1}) \cdot \vec n ,p_{\tau,h}^+(t_{n-1}) \rangle_{\Gamma^D_{\vec u}}}\\[1ex] 
& \quad + A_\gamma (\vec u^+_{\tau,h}(t_{n-1}), u_{\tau,h}^+(t_{n-1}))  + \langle \rho \vec v^+_{\tau,h}(t_{n-1}),  \chi_{\tau,h}^+(t_{n-1})\rangle + \langle c_0 p^+_{\tau,h}(t_{n-1}), p_{\tau,h}^+(t_{n-1})\rangle  = 0 \,.
\end{aligned}
\end{equation}
Now, let $l\in \{1,\ldots,k+1\}$ be arbitrary but fixed and $\vec \phi_{\tau,h}\in \mathbb P_{k}(I_n;\vec V_{h})$ be chosen as  
\begin{equation*}
\vec \phi_{\tau,h}(t) := \xi_{n,l}(t) \vec \phi_h \quad \text{with}\quad \xi_{n,l}(t) := \left(\prod_{i=1 \atop i \neq l}^{k+1}\big(t-t_{n,i}^{\text{GR}}\big)\right) \left(\prod_{i=1 \atop i \neq l}^{k+1}\big(t_{n,l}^{\text{GR}} -t_{n,i}^{\text{GR}}\big)\right)^{-1}  \in \mathbb P_{k}(I_n;\R)\,, \;\; \vec \phi_h \in \vec V_h\,, 
\end{equation*}
and the Gauss--Radau quadrature nodes $t_{n,\mu}^{\text GR}$, for $\mu = 1,\ldots , k+1$; cf.~\eqref{Eq:GF}. By the exactness of the Gauss--Radau quadrature formula \eqref{Eq:GF} for all polynomials in $\mathbb P_{2k}(I_n;\R)$ we deduce from \eqref{Eq:DPQ_1} that
\begin{align*}
0 & = \frac{\tau_n}{2}\sum_{\mu=1}^{k+1} \hat \omega_\mu^{\text{GR}}(\langle \partial_t \vec u_{\tau,h}(t_{n,\mu}^{\text{GR}}), \vec \phi_{\tau,h}(t_{n,\mu}^{\text{GR}})  \rangle - \langle \vec v_{\tau,h}(t_{n,\mu}^{\text{GR}}), \vec \phi_{\tau,h}(t_{n,\mu}^{\text{GR}}) \rangle) + \langle \vec u^+_{\tau,h}(t_{n-1}), \vec \phi_{\tau,h}^+(t_{n-1})\rangle \\[1ex]
& = \frac{\tau_n}{2} \hat \omega_l^{\text{GR}}  \langle \partial_t \vec u_{\tau,h}(t_{n,l}^{\text{GR}}) - \vec v_{\tau,h}(t_{n,l}^{\text{GR}}), \vec \phi_h \rangle + \langle \vec u^+_{\tau,h}(t_{n-1}),  \xi^+_{n,l}(t_{n-1}) \vec \phi_h \rangle\,.
\end{align*}
Thus, we have that ($l=1,\ldots, k+1$)
\begin{equation}
\label{Eq:ExUniDS_Z2}
 \vec v_{\tau,h}(t_{n,l}^{\text{GR}}) - \partial_t \vec u_{\tau,h}(t_{n,l}^{\text {GR}})  = c_{n,l} \, \vec u^+_{\tau,h}(t_{n-1})\quad \text{with}\;\; c_{n,l} = 2 \tau_n^{-1}\,  \left(\hat \omega_l^{\text{GR}}\right)^{-1}\, \xi^+_{n,l}(t_{n-1})\,.
\end{equation}
Substituting \eqref{Eq:ExUniDS_Z2} into the third term on the left-hand side of \eqref{Eq:ExUniDS_Z1}, we get that 
\begin{align}
\nonumber
Q_n\big(C_\gamma (\vec v_{\tau,h} - \partial_t \vec u_{\tau,h},p_{\tau,h}) \big) & = \sum_{\mu=1}^{k+1}  C_\gamma (\vec u_{\tau_h}^+(t_{n-1}),\xi^+_{n,\mu}(t_{n-1}) p_{\tau,h}(t_{n,\mu}^{\text{GR}}) ) = C_\gamma (\vec u_{\tau_h}^+(t_{n-1}), p_{\tau,h}^+(t_{n-1}))\\[1ex]
& = - \alpha \langle \nabla \cdot \vec u_{\tau_h}^+(t_{n-1}), p_{\tau,h}^+(t_{n-1})\rangle + \alpha \langle \vec u_{\tau_h}^+(t_{n-1}) \cdot \vec n ,  p_{\tau,h}^+(t_{n-1}) \rangle_{\Gamma^D_{\vec u}} \,.
\label{Eq:ExUniDS_Z3}
\end{align}
Combining \eqref{Eq:ExUniDS_Z1} with \eqref{Eq:ExUniDS_Z3} then implies that 
\begin{equation}
	\label{Eq:ExUniDS_Z4}
		\begin{aligned}
			& A_\gamma (\vec u_{\tau,h}(t_n) , \vec u_{\tau,h}(t_n))  + \langle \rho \vec v_{\tau,h} (t_n), \vec v_{\tau,h}(t_n)  \rangle  + \langle c_0 p_{\tau,h}(t_n), p_{\tau,h}(t_n)\rangle +  2 Q_n \big(B_\gamma(p_{\tau,h},p_{\tau,h}) \big) \\[1ex]
			& \quad + A_\gamma(\vec u^+_{\tau,h}(t_{n-1}), \vec  u_{\tau,h}^+(t_{n-1}))  + \langle \rho \vec v^+_{\tau,h}(t_{n-1}),  \chi_{\tau,h}^+(t_{n-1})\rangle + \langle c_0 p^+_{\tau,h}(t_{n-1}), p_{\tau,h}^+(t_{n-1})\rangle  = 0 \,.
		\end{aligned}
\end{equation}
From \eqref{Eq:ExUniDS_Z4} along with the discrete coercivity properties \eqref{Eq:CoercA} and \eqref{Eq:CoercB} we deduce that 
	\begin{equation*}
	\label{Eq:ExUniDS_Z5}
		\vec u_{\tau,h}(t_n)= \vec u_{\tau,h}^+(t_{n-1})=\vec 0\,, \quad \vec v_{\tau,h}(t_n)= \vec v_{\tau,h}^+(t_{n-1})=\vec 0 \quad \text{and} \quad p_{\tau,h}(t_n)=p_{\tau,h}^+(t_{n-1})= 0
	\end{equation*}
as well as 
\begin{equation*}
	p_{\tau,h}\big(t_{n,\mu}^{\text{GR}}\big) = 0 \,, \quad \text{for}\;\; \mu = 1,\ldots, k+1\,.
\end{equation*}
The rest then follows as in the proof of Lem.~\ref{Lem:ExUniDS}.	
\end{mproof}

\section{Additional numerical experiments}
\label{Sec:AppC}

Here, we present some additional results of our numerical experiments. Table~\ref{Tab:App1} shows the results of the convergence study introduced in Subsec.~\ref{Subsec:NumConv} for the solution \eqref{Eq:givensolution} of \eqref{Eq:HPS}. For the computations of Table~\ref{Tab:App1} the Taylor-Hood pair of finite element spaces $\mathbb Q_r^2/\mathbb Q_{r-1}$, with $r=4$, is used instead of the pair $\mathbb Q_r^2/\mathbb P_{r-1}^{\text{disc}}$ chosen for the results of Table~\ref{Tab:2}. For the pressure approximation marginally smaller errors are observed. In Table~\ref{Tab:App2} the computed errors of the approximation of \eqref{Eq:givensolution} for larger values of the Lam\'e parameters $\lambda$ and $\mu$ in the elasticity tensor $\vec C$ are summarized. No significant increase of the errors is observed, indicating the independence of the error constant on the magnitude of $\vec C$.

	\begin{table}[H]
	\centering
	\begin{tabular}{l}
		\begin{tabular}{cccccccc}
			\toprule
			{$\tau$} & {$h$} &
			{ $\| \nabla (\vec u - \vec u_{\tau,h})  \|_{L^2(\vec L^2)} $ } & {EOC} &
			{ $\| \vec v - \vec v_{\tau,h}\|_{L^2(\vec L^2)}  $ } & {EOC} &
			{ $\| p - p_{\tau,h} \|_{L^2(L^2)}  $ } & {EOC}  \\
			\cmidrule(r){1-2}
			\cmidrule(lr){3-8}
			$\tau_0/2^0$ & $h_0/2^0$ & 2.3958497242e-03 & {--} & 1.7185669625e-02  & {--} & 7.6634974254e-04 & {--} \\ 
			$\tau_0/2^1$ & $h_0/2^1$ & 1.0529091363e-04 & 4.51 & 5.4558642696e-04  & 4.98 & 3.8993048591e-05 & 4.30 \\
			$\tau_0/2^2$ & $h_0/2^2$ & 5.9749805593e-06 & 4.14 & 1.9143193940e-05  & 4.83 & 1.6140682324e-06 & 4.59\\
			$\tau_0/2^3$ & $h_0/2^3$ & 3.6858181764e-07 & 4.02 & 9.6894400439e-07  & 4.30 & 9.0805047399e-08 & 4.15\\
			$\tau_0/2^4$ & $h_0/2^4$ & 2.2976509245e-08 & 4.00 & 5.6982904446e-08  & 4.09 & 5.5524640600e-09 & 4.03\\
			\bottomrule
		\end{tabular}\\
		\mbox{}\\
		\begin{tabular}{cccccccc}
			\toprule
			{$\tau$} & {$h$} &
			{ $\| \nabla (\vec u - \vec u_{\tau,h})  \|_{L^\infty(\vec L^2)} $ } & {EOC} &
			{ $\| \vec v - \vec v_{\tau,h}\|_{L^\infty(\vec L^2)}  $ } & {EOC} &
			{ $\| p - p_{\tau,h} \|_{L^\infty(L^2)}  $ } & {EOC}  \\
			\cmidrule(r){1-2}
			\cmidrule(lr){3-8}
			$\tau_0/2^0$ & $h_0/2^0$ & 1.3014883491e-02 & {--} & 1.1587851233e-01 & {--} & 4.9980604461e-03 & {--}  \\ 
			$\tau_0/2^1$ & $h_0/2^1$ & 8.7944378641e-04 & 3.89 & 4.1645287352e-03 & 4.80 & 2.6289118267e-04 & 4.25  \\
			$\tau_0/2^2$ & $h_0/2^2$ & 5.6086076018e-05 & 3.97 & 2.0121086363e-04 & 4.37 & 9.4136464267e-06 & 4.80  \\
			$\tau_0/2^3$ & $h_0/2^3$ & 3.5171476997e-06 & 4.00 & 1.1873567381e-05 & 4.08 & 5.0122396011e-07 & 4.23  \\
			$\tau_0/2^4$ & $h_0/2^4$ & 2.1959938532e-07 & 4.00 & 7.3034502924e-07 & 4.02 & 2.8669553969e-08 & 4.13  \\
			\bottomrule	
		\end{tabular}	
	\end{tabular}
	
	\caption{%
		$L^2(L^2)$ and $L^\infty(L^2)$  errors and experimental orders of convergence (EOC) with temporal polynomial degree $k=3$ and spatial degree $r=4$ for local spaces $\mathbb Q_r^2/\mathbb Q_{r-1}$.
	}
	\label{Tab:App1}
\end{table}

\begin{table}[H]
	\centering
	\begin{tabular}{l}
		\begin{tabular}{cccccccc}
			\toprule
			{$\tau$} & {$h$} &
			{ $\| \nabla (\vec u - \vec u_{\tau,h})  \|_{L^2(\vec L^2)} $ } & {EOC} &
			{ $\| \vec v - \vec v_{\tau,h}\|_{L^2(\vec L^2)}  $ } & {EOC} &
			{ $\| p - p_{\tau,h} \|_{L^2(L^2)}  $ } & {EOC}  \\
			\cmidrule(r){1-2}
			\cmidrule(lr){3-8}
			$\tau_0/2^0$ & $h_0/2^0$ & 1.1835824122e-02 & {--} & 2.8040747896e-02  & {--} & 2.9610298753e-03 & {--} \\ 
			$\tau_0/2^1$ & $h_0/2^1$ & 1.5714797742e-03 & 2.91 & 5.7609393975e-03  & 2.28 & 3.9787786632e-04 & 2.90 \\
			$\tau_0/2^2$ & $h_0/2^2$ & 1.8935735637e-04 & 3.05 & 5.1606189537e-04  & 3.48 & 3.7075698634e-05 & 3.42\\
			$\tau_0/2^3$ & $h_0/2^3$ & 2.3548556912e-05 & 3.01 & 6.0883693658e-05  & 3.08 & 4.5049814701e-06 & 3.04\\
			$\tau_0/2^4$ & $h_0/2^4$ & 2.9374871960e-06 & 3.00 & 7.5088512252e-06  & 3.02 & 5.6225700500e-07 & 3.00\\
			\bottomrule
		\end{tabular}\\
		\mbox{}\\
		\begin{tabular}{cccccccc}
			\toprule
			{$\tau$} & {$h$} &
			{ $\| \nabla (\vec u - \vec u_{\tau,h})  \|_{L^\infty(\vec L^2)} $ } & {EOC} &
			{ $\| \vec v - \vec v_{\tau,h}\|_{L^\infty(\vec L^2)}  $ } & {EOC} &
			{ $\| p - p_{\tau,h} \|_{L^\infty(L^2)}  $ } & {EOC}  \\
			\cmidrule(r){1-2}
			\cmidrule(lr){3-8}
			$\tau_0/2^0$ & $h_0/2^0$ & 6.5064235342e-02 & {--} & 1.3580389516e-01 & {--} & 1.4838995481e-02 & {--}  \\ 
			$\tau_0/2^1$ & $h_0/2^1$ & 1.1842822460e-02 & 2.46 & 1.8023863764e-02 & 2.91 & 2.1492360068e-03 & 2.79  \\
			$\tau_0/2^2$ & $h_0/2^2$ & 1.5191266076e-03 & 2.96 & 3.4010640997e-03 & 2.41 & 2.3825320454e-04 & 3.17  \\
			$\tau_0/2^3$ & $h_0/2^3$ & 1.9340613477e-04 & 2.97 & 4.6154471881e-04 & 2.88 & 3.0597316621e-05 & 2.96  \\
			$\tau_0/2^4$ & $h_0/2^4$ & 2.4393096482e-05 & 2.99 & 5.8122744877e-05 & 2.99 & 3.8749301950e-06 & 2.98  \\
			\bottomrule	
		\end{tabular}	
	\end{tabular}

	\caption{%
	$L^2(L^2)$ and $L^\infty(L^2)$  errors and experimental orders of convergence (EOC) with temporal polynomial degree $k=2$ and spatial degree $r=3$ for local spaces $\mathbb Q_r^2/\mathbb P_{r-1}^{\text{disc}}$ for Young's modulus $E=10000$ and Poisson's ratio $\nu=0.35$, corresponding to the Lam\'e parameters $\lambda = 8642$ and $\mu = 3704$.}

	\label{Tab:App2}
\end{table}


\begin{thebibliography}{99}
	
\bibitem{A67}	
G.\ M.\ Amdahl, \textit{Validity of the single processor approach to achieving large scale computing capabilities}, in Proceedings of the April 18-20, 1967, Spring Joint Computer Conference on --- AFIPS ’67 (Spring), Atlantic City, New Jersey: ACM Press, 1967, p.~483.	
	
\bibitem{AB21}
M.\ Anselmann, M.\ Bause, \textit{A geometric multigrid method for space-time finite element discretizations of the Navier–Stokes equations and its application to 3d flow simulation}, ACM Trans.\ Math.\ Softw., \textbf{accepted} (2023), https://doi.org/10.1145/3582492, pp.\ 1--27; arXiv:2107.10561.

\bibitem{AB22}
M.\ Anselmann, M.\ Bause, \textit{Efficiency of local Vanka smoother geometric multigrid preconditioning for space-time finite element methods to the Navier-Stokes equations}, PAMM Proc.\ Appl.\ Math.\ Mech., \textbf{accepted} (2022), pp.\ 1--6; arXiv:2210.02690.
	
\bibitem{AB21_2} 
M.\ Anselmann, M.\ Bause, \textit{CutFEM and ghost stabilization techniques for higher order space-time discretizations of the Navier--Stokes equations}, Int.\ J.\ Numer.\ Meth.\ Fluids, \textbf{94} (2022), pp.~775--802.

\bibitem{ABBM20} 
M.\ Anselmann, M.\ Bause, S.\ Becher, G.\ Matthies, \textit{Galerkin-collocation approximation in time for the wave equation and its post-processing}, ESAIM Math.\ Model.\ Numer.\ Anal., \textbf{54} (2020), pp.\ 2099--2123.

\bibitem{Aetal21}
D.\ Arndt, W.\ Bangerth, B.\ Blais, M.\ Fehling, R.\ Gassmöller, T.\ Heister, L.\ Heltai, U.\ Köcher, M.\ Kronbichler, M.\ Maier, P.\ Munch, J.-P.\ Pelteret, S.\ Proell, K.\ Simon, B.\ Turcksin, D.\ Wells, J.\ Zhang, \textit{The deal.II Library, Version 9.3}, J.\ Numer.\ Math., 29 (2021), pp.~171--186.

\bibitem{A02}
D.\ N.\ Arnold, D.\ Boffi, R.\ S.\ Falk, \textit{Approximation by quadrilateral finite elements}, Math.\ Comp., \textbf{71} (2002), pp.\ 909--922.

\bibitem{BGR10} 
W.\ Bangerth, M.\ Geiger, R.\ Rannacher, \textit{Adaptive {G}alerkin finite element methods for the wave equation},
Comput.\ Meth.\ Appl.\ Math., \textbf{10} (2010), pp.~3--48.
 
\bibitem{BR03}
W.\ Bangerth, R.\ Rannacher, \textit{Adaptive Finite Element Methods for Differential Equations}, Birkhäuser, Basel, 2003.

\bibitem{BKR22}
M.\ Bause, M.\ Anselmann, U.\ Köcher, F.\ A.\ Radu, \textit{Convergence of a continuous Galerkin method for hyperbolic-parabolic systems}, Comput.\ Math.\ with Appl., \textbf{submitted} (2022), pp.\ 1--24;  arXiv:2201.12014.

\bibitem{BRK17} 
M.\ Bause, R.\ Radu, U.\ Köcher, \textit{Space-time finite element approximation of the Biot poroelasticity system with iterative coupling}, Comput.\ Methods Appl.\ Mech.\ Engrg., \textbf{320} (2017), pp. 745--768.

\bibitem{B02}
R.\ Becker, \textit{Mesh adaptation for Dirirchlet flow control via Nitsche's method}, Commun.\ Numer.\ Meth.\ Engrg., \textbf{18} (2002), pp.\ 669--680.
	
\bibitem{B41}
M.\ Biot, \textit{General theory of three-dimensional consolidation}, J.\ Appl.\ Phys., \textbf{12} (1941), pp.~155--164.

\bibitem{B55}
M.\ Biot, \textit{Theory of elasticity and consolidation for a porous anisotropic solid}, J.\ Appl.\ Phys., \textbf{26} (1955), pp.~182--185.

\bibitem{B72}
M.\ Biot, \textit{Theory of finite deformations of porous solids, Indiana Univ.\ Math.\ J.}, \textbf{21} (1972), pp.~597--620.
		
\bibitem{BKNR22}
J.\ W.\ Both, N.\ A.\ Barnafi, F.\ A.\ Radu, P.\ Zunino, A.\ Quarteroni, \textit{Iterative splitting schemes for a soft material poromechanics model}, Comput.\ Methods Appl.\ Mech.\ Engrg., \textbf{388} (2022), 114183.

\bibitem{BL11}
A.\ Brandt, O.\ E.\ Livne, \textit{Multigrid Techniques—1984 Guide with Applications to Fluid Dynamics}, SIAM, Philadelphia, 2011. 

\bibitem{B03}
S.\ C.\ Brenner, \textit{Korn's inequalities for piecewise $\vec H^1$ vector fields},  Math.\ Comp., \textbf{73} (2003), pp.~1067--1087.
	
\bibitem{BKB22}
M.\ P.\ Bruchh\"auser,  U.\ Köcher, M.\ Bause, \textit{On the implementation of an adaptive multirate framework for coupled transport and flow}, J.\ Sci.\  Comput.,  \textbf{93:59} (2022), https://doi.org/10.1007/s10915-022-02026-z, pp.\ 1--29. 

\bibitem{C72}
D.\ E.\ Carlson, \textit{Linear thermoelasticity}, Handbuch der Physik  V Ia/2, Springer, Berlin, 1972.

\bibitem{C00}
W.\ W.\ Cooper, L.\ M.\ Seiford, K.\ Tone \textit{Data Envelopment Analysis}, Kluwer Academic Publishers, Dordrecht, 2000.

\bibitem{PE12}
D.\ A.\ Di Pietro, A.\ Ern, \textit{Mathematical Aspects of Discontinuous Galerkin Methods}, Springer, Berlin, 2012. 

\bibitem{DJRWZ18}
D.\ Drzisga, L.\ John, U.\ Rüde, B.\ Wohlmuth, W.\ Zulehner, \textit{On the analysis of block smoothers for saddle point problems}, SIAM J.\ Matrix Anal.\ Appl., \textbf{39} (2018), pp.~932--960.

\bibitem{F78}
G.\ Fairweather, \textit{Finite Element Galerkin Methods for Differential Equations}, Lecture Notes in Pure and Applied Mathematics 34, Marcel Dekker Inc., New York, 1978.

\bibitem{FTW19}
S.\ Franz, S.\ Trostorff, M.\ Waurick, \textit{Numerical methods for changing type systems}, IMA J.\ Numer.\ Anal., \textbf{39} (2019), pp.~1009--1038.

\bibitem{GHJRW16}
B.\ Gmeiner, M.\ Huber, L.\ John, U.\ Rüde, B.\ Wohlmuth, \textit{A quantitative performance study for Stokes solvers at the extreme scale}, J.\ Comput.\ Sci., \textbf{17} (2016), pp.~509--521.

\bibitem{GRSWW15}
B.\ Gmeiner, U.\ Rüde, H.\ Stengel, C.\ Waluga, B.\ Wohlmuth, \textit{Performance and scalability of hierarchical hybrid multigrid solvers for Stokes systems}, SIAM J.\ Sci.\ Comput., \textbf{37} (2015), pp.~C143--C168.

\bibitem{GRSWW15_2}
B.\ Gmeiner, U.\ Rüde, H.\ Stengel, C.\ Waluga, B. \ Wohlmuth, \textit{Towards textbook efficiency for parallel multigrid}, Numer.\ Math.\ Theory Methods Appl., \textbf{8} (2015), pp.~22--46.

\bibitem{H85}
W.\ Hackbusch, \textit{Multigrid methods and applications}, Springer, Berlin, 1985.

\bibitem{HK18}
Q.\ Hong, J.\ Kraus, \emph{Parameter-robust stability of classical three-field formulation of Biot's consolidation model}, Electron.\ Trans.\ Numer.\ Anal., \textbf{48} (2018), pp.~202--226.

\bibitem{HKXZ16}
Q.\ Hong, J.\ Kraus, J.\ Xu, L.\ Zikatanov, \textit{A robust multigrid method for discontinuous Galerkin discretizations of Stokes and linear elasticity equations}, Numer.\ Math., \textbf{132} (2016), pp.~23--49.

\bibitem{HST14}
S.\ Hussain, F.\ Schieweck, S.\ Turek, \textit{Efficient Newton-multigrid solution techniques for higher order space-time Galerkin discretizations of incompressible flow}, Appl.\ Numer.\ Math., \textbf{83} (2014),  pp.~51--71.

\bibitem{HST13}
S.\ Hussain, F.\ Schieweck, S.\ Turek, \textit{An efficient and stable finite element solver of higher order
in space and time for nonstationary incompressible flow}, Internat.\ J. Numer.\ Methods Fluids, \textbf{73} (2013), pp.\ 927--952.
		
\bibitem{HST11}
S.\ Hussain, F.\ Schieweck, S.\ Turek, \textit{Higher order Galerkin time discretizations and fast multigrid solvers for the heat equation}, J.\ Numer.\ Math., \textbf{19} (2011), pp.\ 41--61.

\bibitem{JR18}
S.\ Jiang, R.\ Racke, \textit{Evolution equations in thermoelasticity}, CRC Press, Boca Raton, 2018.
	
\bibitem{J16}
V.\ John, \textit{Finite Element Methods for Incompressible Flow Problems}, Springer, Cham, 2016.
	
\bibitem{J02}
V. \ John, \textit{Higher order finite element methods and multigrid solvers in a benchmark problem for the 3D Navier-Stokes equations}, Int. J. Numer. Meth. Fluids, \textbf{40} (2002), pp.~775--798.

\bibitem{JM01}	
V.\ John, G.\ Matthies,	\textit{Higher-order finite element discretizations in a benchmark problem for incompressible flows}, Int.\ J.\ Numer.\ Meth.\ Fluids, \textbf{37} (2001), pp.~885--903.	
	
\bibitem{JT00}
V.\ John, L. \ Tobiska, \textit{Numerical performance of smoothers in coupled multigrid methods for the parallel solution of the incompressible Navier–Stokes equations}, Int.\ J.\ Numer.\ Meth.\ Fluids, \textbf{33} (2000),  pp.~453--473. 

\bibitem{KR18}
 G.\ Kanschat, B.\ Riviere, \textit{A finite element method with strong mass conservation for Biot's linear consolidation model}, J.\ Sci.\ Comput., \textbf{77} (2018), pp.~1762--1779.

\bibitem{KM04}
O.\ Karakashian, C.\ Makridakis, \textit{Convergence of a continuous Galerkin method with mesh modification for nonlinear wave equations}, Math.\ Comp., \textbf{74} (2004), pp.\ 85--102.

\bibitem{KM99}
O.\ Karakashian, C.\ Makridakis, \textit{A space-time finite element method for the nonlinear Schrödinger equation: The continuous Galerkin method}, SIAM J.\ Numer.\ Anal., \textbf{36} (1999), pp.\ 1779--1807.

\bibitem{L86}
R.\ Leis, \textit{Initial boundary value problems in mathematical physics}, Teubner, Stuttgart, John Wiley \& Sons, Chichester, 1986.

\bibitem{LD03}
X.\ S.\ Li, J.\ W.\ Demmel, \textit{SuperLU\_DIST: A scalable distributed-memory sparse direct solver for unsymmetric linear systems}, ACM Trans.\ Math.\ Softw., \textbf{29} (2003), pp.\ 110--140.
\bibitem{LLRS94}
J.\ Linden, G.\ Lonsdale, H.\ Ritzdorf, A. Schüller, \textit{Scalability aspects of parallel multigrid}, Future Gener.\ Comput.\ Syst., \textbf{10} (1994), pp.~429--439.

 \bibitem{M06}
S.\ Manservisi, \textit{Numerical analysis of Vanka-type solvers for steady Stokes and Navier--Stokes flows},  SIAM.\ J.\ Numer.\ Anal., \textbf{44} (2006), pp.~2025--2056.

\bibitem{M01}
G.\ Matthies, \textit{Mapped finite elements on hexahedra. Necessary and sufficient conditions for optimal interpolation errors}, Numer.\ Algor., \textbf{27} (2001), pp.\ 317--327.

\bibitem{MT02}
G.\ Matthies, L.\ Tobiska, \textit{The inf-sup condition for the mapped $\mathbb Q_k^d/P_{k-1}^\disc$ element in arbitrary space dimensions}, Computing, \textbf{69} (2002), pp.\ 119--139. 
		
\bibitem{MW12}
A.\ Mikeli\'{c}, M.\ F.\ Wheeler, {\em Theory of the dynamic Biot--Allard equations and their link to the quasi-static Biot system}, J.\ Math.\ Phys., \textbf{53} (2012), 123702:1--15.

\bibitem{N71}
J.\ Nitsche, \emph{Über ein Variationsprinzip zur Lösung von Dirichlet Problemen bei Verwendung von Teilräumen, die keinen Randbedingungen unterworfen sind (in German)},  Abh.\ Math.\ Sem.\ Univ.\ Hamburg, \textbf{36} (1971), pp.\ 9--15.

\bibitem{PW07}
P.\ J.\ Philips, M.\ F.\ Wheeler, {\em A coupling of mixed and continuous Galerkin finite 
	element methods for poroelasticity I, II}, Comput.\ Geosci., \textbf{11} (2007), 
131--158.

\bibitem{PW08}
P.\ J.\ Philips, M.\ F.\ Wheeler, {\em A coupling of mixed and discontinuous Galerkin 
	finite 
	element methods for poroelasticity }, Comput.\ Geosci., \textbf{12} 
(2008), 417--435.

\bibitem{RHOAGZ17}
C.\ Rodrigo, X.\ Hu, P.\ Ohm, J.\ H.\ Adler, F.\ J.\ Gaspar, L.\ T.\ Zikatanov, \textit{New stabilized discretizations for poroelasticity and the Stokes' equations}, Comput.\ Methods Appl.\ Mech.\ Engrg., \textbf{341} (2018), pp.\ 467--484.

\bibitem{R17}
U.\ Rüde, \textit{Algorithmic Efficiency and the Energy Wall}, 2017. DOI: 10.13140/RG.2.2.33914.18881.

\bibitem{STW22}
C.\ Seifert, S.\ Trostorff, M.\ Waurick, \textit{Evolutionary Equations: Picard's Theorem for Partial Differential Equations, and Applications}, Birkhäuser, Cham, 2022.

\bibitem{S03}
Y.\ Saad, \textit{Iterative Methods for Sparse Linear Systems. 2nd edition}. SIAM, Philadelphia, 2003.

\bibitem{S00}
R.\ Showalter, \textit{Diffusion in poro-elastic media}, J.\ Math.\ Anal.\ Appl.,  \textbf{251} (2000), pp.\ 310--340.

\bibitem{S89}
M.\ Slodi\v{c}ka, \textit{Application of Rothe's method to integrodifferential equation}, Comment.\ Math.\ Univ.\ Carolinae, \textbf{30} (1989), pp.~57--70.

\bibitem{S23}
Slurm Workload Manager, Version 22.05, \url{https://github.com/SchedMD/slurm/tree/master/src/plugins/acct_gather_energy/rapl}, 2023. 


\bibitem{SZ22}
O.\ Steinbach, M.\ Zank, \textit{A generalized inf–sup stable variational formulation for the wave equation}, J.\ Math.\ Anal.\ Appl.,  \textbf{505} (2022), 125457. 

\bibitem{T06}
V.\ Thome\'e, \textit{Galerkin Finite Element Methods for Parabolic Problems}, Springer, Berlin, 2006.

\bibitem{TOS01}
U.\ Trottenberg, C.\ W.\ Oosterlee, A.\ Schüller, Multigrid, Academic Press, San Diego, CA, 2001.

\bibitem{T99}
S.\ Turek, \textit{Efficient Solvers for Incompressible Flow Problems}, Springer,  Berlin, 1999.

\bibitem{TBK06}
S.\ Turek, C.\ Becker, D.\ Kilian, \textit{Hardware-oriented numerics and concepts for PDE software}, Future Generation Computer Systems, \textbf{22} (2006), pp.~217--238.
	
\bibitem{TGBBW10}
S.\ Turek, D.\ Göddecke, C.\ Becker, S.\ Buijssen, H.\ Wobker, \textit{FEAST --- realisation of hardware-oriented numerics for HPC simulations with finite elements}, Concurrency and Computation: Practice and Experience 6
(May), pp.\ 2247--2265. Special Issue Proceedings of ISC 2008, doi:10.1002/cpe.1584.
	
\bibitem{V86}
S.\ Vanka, \textit{Block-implicit multigrid solution of Navier-–Stokes equations in primitive variables},
J.\ Comput.\ Phys., \textbf{65} (1986), pp.~138--158.

\bibitem{WS00}
D.\ C.\ Wirtz, N.\ Schiffers, T.\ Pandorf, K.\ Rademacher, D.\ Weichert, R.\ Forst, \textit{Critical evaluation of known bone material properties to realize anisotropic FE-simulation of the proximal femur},
J.\ Biomech., \textbf{33} (2000), pp.~1325--1330.

\end{thebibliography}
\end{document}